\documentclass[11pt]{amsart}

\usepackage{amsmath, amsthm, amssymb}

\usepackage[dvipsnames]{xcolor}
\definecolor{darkgreen}{RGB}{0,150,0}
\definecolor{darkred}{RGB}{200,0,0}
\definecolor{darkblue}{RGB}{0,50,205}
\definecolor{darkpurple}{RGB}{100,15,200}
\definecolor{gold}{RGB}{255,218,9}
\usepackage[colorlinks, linkcolor=darkblue, citecolor=darkgreen, urlcolor=darkgreen]{hyperref}
\hypersetup{hypertexnames=false}

\usepackage{newpxtext}
\usepackage[bigdelims,vvarbb]{newpxmath}  
\usepackage[T1]{fontenc}  

\usepackage{setspace}
\usepackage[margin=1in]{geometry}
             
\usepackage{comment}
\usepackage{thmtools, thm-restate}
\usepackage{overpic}
\usepackage{subcaption}

\usepackage{aliascnt}
\usepackage[noabbrev,capitalize,nameinlink]{cleveref}

\newtheorem{theorem}{Theorem}[section]
\newtheorem{proposition}[theorem]{Proposition}
\newtheorem{lemma}[theorem]{Lemma}
\newtheorem{corollary}[theorem]{Corollary}
\newtheorem{conjecture}[theorem]{Conjecture}
\newtheorem{mainthm}{Theorem}

\newtheorem{maincor}[mainthm]{Corollary}
\theoremstyle{definition}
\newtheorem{definition}[theorem]{Definition}

\theoremstyle{remark}
\newtheorem{remark}[theorem]{Remark}
\numberwithin{equation}{section}

\usepackage{mleftright}
\newcommand{\set}[2]{\mleft\{\,#1~\middle|~#2\,\mright\}}
\newcommand{\R}{\mathbb{R}}
\newcommand{\C}{\mathbb{C}}

\newcommand{\Z}{\mathbb{Z}}

\newcommand{\ve}{\varepsilon}
\newcommand{\st}{\mathrm{st}}
\newcommand{\tb}{\mathrm{tb}}
\newcommand{\rot}{\mathrm{rot}}

\newcommand{\RNum}[1]{\uppercase\expandafter{\romannumeral #1\relax}}

\makeatletter 
 \def\l@subsection{\@tocline{2}{0pt}{2pc}{6pc}{}} \makeatother

\begin{document}

\title{Convex ribbon disks, quasipositivity, and Lagrangian sliceness}

\author{Joseph Breen}
\address{University of Alabama, Tuscaloosa, AL}
\email{jjbreen@ua.edu} \urladdr{https://sites.google.com/view/joseph-breen}

\author{Alexander Zupan}
\address{University of Nebraska - Lincoln, Lincoln, NE}
\email{zupan@unl.edu} \urladdr{http://www.math.unl.edu/~azupan2}

\thanks{JB was partially supported by NSF grant DMS-2038103 and an AMS-Simons Travel Grant. AZ was partially supported by NSF grant DMS-2405301 and a Simons Travel Support award.}

\begin{abstract}
We introduce a version of convex surface theory for immersed ribbon surfaces with Legendrian boundary. As an application, we establish a special case of K3: Problem 1.70 on sharpness in the slice-Bennequin inequality implying quasipositivity, and a special case of a conjecture of Cornwell, Ng, and Sivek on characterizing Lagrangian sliceness. Specifically, we prove that if $K$ is a knot that bounds a "star-shaped" ribbon disk and satisfies $\overline{\mathrm{sl}}(K) = -1$ (resp.\ $\overline{\tb}(K) = -1$), then $K$ is quasipositive (resp.\ Lagrangian slice). In fact, we prove these statements for any fixed transverse (resp.\ Legendrian) representative.
\end{abstract}

\maketitle

\tableofcontents

\section{Introduction}\label{sec:intro}

Given a knot $K\subset S^3$, the \emph{slice-Bennequin inequality}, established by Rudolph \cite{rudolph1993quasipositivity} as a consequence of the local Thom conjecture \cite{kronheimer1994genus}, is
\begin{equation}\label{eq:slice-bennequin-inequality}
\overline{\mathrm{sl}}(K) \leq 2g_4(K) - 1.    
\end{equation}
Here, $\overline{\mathrm{sl}}(K)$ denotes the maximum self-linking number among all transverse representatives of $K$ in $(S^3, \xi_{\st})$, and $g_4(K)$ is the smooth $4$-genus. Quasipositivity of $K$ implies equality in \eqref{eq:slice-bennequin-inequality}, and it is an open question as to whether the converse holds; see \cite[Question 7.3]{etnyre2015monoids}, \cite[Question 1.12]{hamer2022defect}, or the recent K3 version \cite[Problem 1.70]{K3}.

Rudolph \cite{rudolph1983seifertribbons} and Boileau and Orevkov \cite{boileau2001quasipositive} established an equivalence between quasipositive knots and knots with transverse representatives bounding symplectic surfaces in $(B^4, \omega_{\st})$. Consequently, if a knot $K$ admits a transverse representative which is symplectically slice, it is quasipositive and thus attains equality in \eqref{eq:slice-bennequin-inequality} with $g_4(K) = 0$, implying $\overline{\mathrm{sl}}(K) = -1$.  Conversely, the $g_4 = 0$ version of \cite[Problem 1.70]{K3}, which we pose as a conjecture, is: 

\begin{conjecture}\label{q:sl=-1}
If $K\subset S^3$ is smoothly slice and $\overline{\mathrm{sl}}(K) = -1$, then $K$ is quasipositive.
\end{conjecture}

There is a similar story for Legendrian knots bounding Lagrangian disks. Work of Chantraine \cite{chantraine2010concordance} implies that if a topological knot type $K$ admits a Legendrian representative which is Lagrangian slice (that is, $K$ is \emph{Lagrangian slice}), then the maximum Thurston-Bennequin invariant of $K$ is $\overline{\mathrm{tb}}(K) = -1$. Cornwell, Ng, and Sivek conjectured, and verified through $14$ crossings, that smooth sliceness and $\overline{\mathrm{tb}}$ alone are sufficient to conclude Lagrangian sliceness:

\begin{conjecture}[\cite{cornwell2016obstructions}]\label{conj:CNS}
If $K\subset S^3$ is smoothly slice and $\overline{\mathrm{tb}}(K) = -1$, then $K$ is Lagrangian slice.     
\end{conjecture}

A natural strengthening of these conjectures concerns not only the topological knot type, but any fixed transverse or Legendrian representative; here, a transverse knot is \emph{quasipositive} if it is transversely isotopic to the closure of a quasipositive braid.

\begin{conjecture}\label{conj:QP_fixed_rep}
If $\mathcal{T}\subset (S^3, \xi_{\st})$ is transverse, smoothly slice, and $\mathrm{sl}(\mathcal{T}) = -1$, then $\mathcal{T}$ is quasipositive.     
\end{conjecture}

\begin{conjecture}\label{conj:CNS_fixed_rep}
If $\Lambda\subset (S^3, \xi_{\st})$ is Legendrian, smoothly slice, and $\mathrm{tb}(\Lambda) = -1$, then $\Lambda$ is Lagrangian slice.
\end{conjecture}

Beyond analogy, there is a strong connection between \cref{conj:QP_fixed_rep} and \cref{conj:CNS_fixed_rep}. First, results of Eliashberg \cite{eliashberg1995pushoff} imply that if a Legendrian $\Lambda$ is Lagrangian slice, then its transverse pushoff is symplectically slice and hence quasipositive; second, Mark and Tosun \cite[Corollary 2.10]{mark2024fillable} proved that if a smoothly slice Legendrian knot $\Lambda$ with $\mathrm{tb}(\Lambda) = -1$ has a quasipositive transverse pushoff, then $\Lambda$ is Lagrangian slice. Thus, \cref{conj:QP_fixed_rep} implies \cref{conj:CNS_fixed_rep}. Clearly, \cref{conj:QP_fixed_rep} implies \cref{q:sl=-1}, and \cref{conj:CNS_fixed_rep} implies \cref{conj:CNS}, so \cref{conj:QP_fixed_rep} implies all four of the above conjectures.

Note that a quasipositive slice knot is automatically ribbon, so ribbonness is a necessary conclusion in a verification of \cref{q:sl=-1} and \cref{conj:CNS}. Therefore, it is natural to modulate by the slice-ribbon conjecture and instead pose these questions under the assumption of ribbonness, rather than sliceness. 

\begin{figure}[ht]
	\centering
	\begin{overpic}[scale=.4]{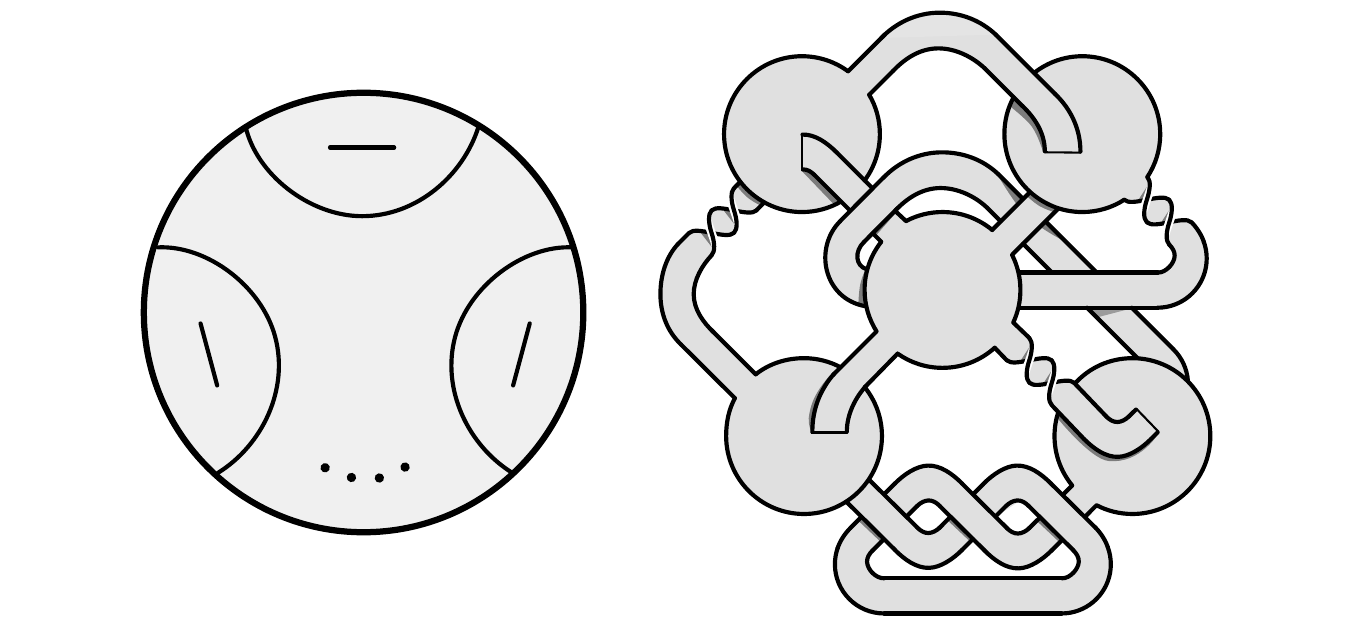}
       
	\end{overpic}
	\caption{An abstract star-shaped ribbon disk (left); a concrete example with $4$ ribbon singularities (right).}
	\label{fig:starshaped1}
\end{figure}

With this in mind, the main result of this article is a verification of \cref{conj:QP_fixed_rep} (thus all four conjectures above) for knots bounding \emph{star-shaped ribbon disks}; see \cref{fig:starshaped1}, and \cref{sec:ribbon_disk_analysis} for a precise definition. For clarity, we state both the transverse and Legendrian versions separately.

\begin{mainthm}\label{cor:slice-benn-sharp}
Let $\mathcal{T}\subset (S^3, \xi_{\st})$ be transverse with $\mathrm{sl}(\mathcal{T}) = -1$. If $\mathcal{T}$ bounds a star-shaped ribbon disk, then $\mathcal{T}$ is quasipositive.
\end{mainthm}

\begin{mainthm}\label{cor:CNSrib2}
Let $\Lambda\subset (S^3, \xi_{\st})$ be Legendrian with $\mathrm{tb}(\Lambda) = -1$. If $\Lambda$ bounds a star-shaped ribbon disk, then $\Lambda$ is Lagrangian slice.
\end{mainthm}

\begin{maincor}
\cref{q:sl=-1} and \cref{conj:CNS} hold for knots bounding star-shaped ribbon disks.
\end{maincor}
The star-shaped assumption in these theorems is in place for combinatorial reasons and is non-trivial; for instance, if $K$ bounds a star-shaped ribbon disk, $\det(K)$ is the square of products of integers of the form $2^n -1$. Nevertheless, ribbon knots bounding star-shaped disks and satisfying $\overline{\tb} = -1$ or $\overline{\mathrm{sl}} = -1$ appear broadly and in infinitude.

The first author, together with Roy and Wang \cite{BRW_PANLF}, recently proved that every Lagrangian disk with Legendrian boundary in a $4$-dimensional Weinstein domain is regular in the sense of \cite{eliashberg2020flexible}, up to Weinstein deformation equivalence. It follows that:

\begin{corollary}
If a Legendrian knot $\Lambda\subset (S^3, \xi_{\st})$ is star-shaped ribbon and has $\tb(\Lambda) = -1$, then $\Lambda$ is regularly Lagrangian slice.  
\end{corollary}

There are three ingredients in the proofs of \cref{cor:slice-benn-sharp} and \cref{cor:CNSrib2}. The first, and the main technical contribution of this article, is a version of convex surface theory for ribbon disks. The second is a result from the predecessor of this work, \cite{BZ0forthcoming}, which we will summarize in \cref{subsec:BZ1intro} and review as required in \cref{sec:BZ1review}.  The third is an intricate combinatorial argument limiting the ways in which dividing curves can interact with star-shaped ribbon disks, appearing in \cref{sec:ribbon_disk_analysis}.

\subsection{Convex ribbon disks}
Convex surface theory was introduced by Giroux \cite{giroux1991convexite}, further developed by the likes of Kanda \cite{kanda1997classification,kanda1998legendrian} and Honda \cite{honda2000classification}, and has since become a permanent fixture of $21$st-century $3$-dimensional contact topology; see \cite{etnyre2002tight,etnyre2005cabling,ghiggini2006classification,colin2011sutures}, or \cite{dimitroglou2024c0limits,licata2024heegaard,etnyre2026pnastori,fernandez2026spaces} for some recent examples. Convex surfaces are $C^{\infty}$-dense in contact manifolds, and come equipped with a multicurve isotopy class, the dividing set, which encodes all contact topology in a neighborhood of the surface. The power of convexity lies in this reduction of semi-local contact structure to the combinatorics of curves on surfaces. 

\begin{figure}[ht]
	\centering
	\begin{overpic}[scale=.3]{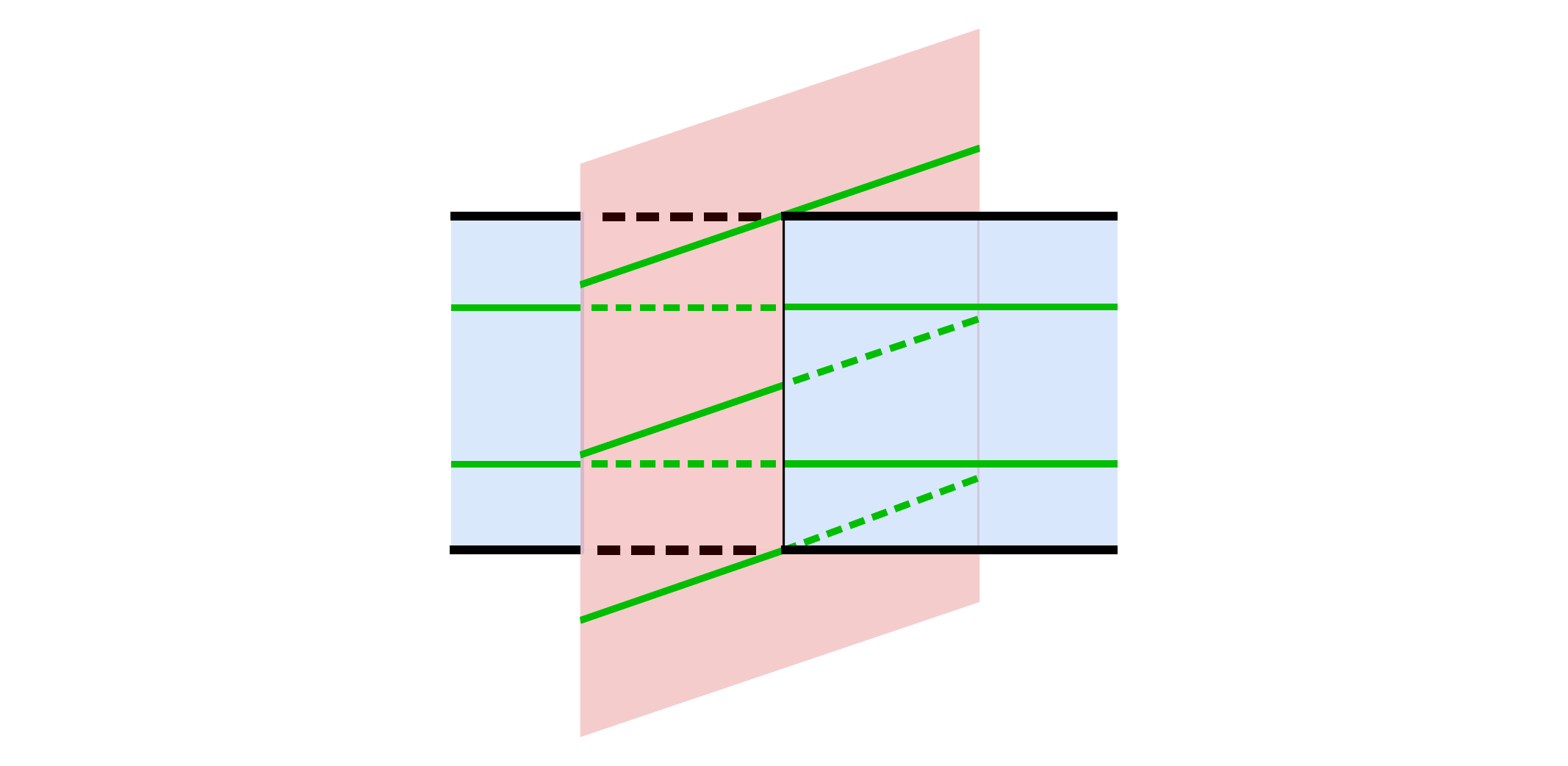}
       
	\end{overpic}
	\caption{A convex ribbon disk near a Legendrian ribbon intersection.}
	\label{fig:normalized-intro}
\end{figure}

By definition, convex surfaces are embedded. In this article we define a notion of convexity for immersed ribbon surfaces with Legendrian boundary and Legendrian ribbon intersections, where the dividing curves of the intersecting sheets have a standard "interleaving" model as in \cref{fig:normalized-intro}. While there has been past experimentation with this idea among experts, this is (to our knowledge) the first instance of such a framework being put to practical use.

This aspect of our work is summarized by the following informal theorem. The definition and existence statement is formulated precisely in \cref{sec:ribbon}, and the manipulations come as a suite of lemmas collected in \cref{sec:manipulation}.

\begin{mainthm}\label{thm:main-ribbon-disk-convex}
Let $D\subset (S^3, \xi_{\st})$ be a ribbon disk with Legendrian boundary $\Lambda$. There is an isotopy relative to $\Lambda$ so that $D$ is "convex" in the sense suggested by \cref{fig:normalized-intro}. Moreover, there are several semi-local "bypass manipulations" of such a convex ribbon disk which modify the dividing curves. 
\end{mainthm}

We expect this notion of ribbon convexity to be of general interest beyond our own applications.

\subsection{Derivative links}\label{subsec:BZ1intro}

Given a knot $K\subset S^3$ and a genus-$g$ Seifert surface $F$ bounding $K$, a \emph{derivative of $K$} is $g$-component link $J = J_1 \cup \cdots \cup J_g \subset F$ such that $F - J$ is connected and planar, and such that the Seifert form of $F$ vanishes on the summand of $H_1(F; \Z)$ generated by $J$. Since the work of J.\ Levine in the 1960s \cite{levine1969knotcobordismcodim2}, derivative links have been fruitfully used to study sliceness. For example, $K$ is ribbon if and only if $K$ admits a Seifert surface with an unlink derivative \cite{cochran2015davis,miller2023handle}.

In \cite{BZ0forthcoming}, the authors initiated a systematic study of Lagrangian and symplectic sliceness via derivative links. Several characterizations of different types of sliceness were established in terms of the existence of derivative links.  For our purpose here, we only need the following implication of one of those results, whose proof will be reviewed for completeness in \cref{sec:BZ1review}: 

\begin{restatable}{theorem}{bz}(\cite{BZ0forthcoming})\label{thm:BZ1review}
Let $\mathcal{T}\subset (S^3, \xi_{\st})$ be a transverse knot with $\mathrm{sl}(\mathcal{T}) = -1$. If $\mathcal{T}$ admits a Seifert surface with a max-sl transverse unlink derivative, then $\mathcal{T}$ is symplectically slice, hence quasipositive. 
\end{restatable}

\begin{figure}[ht]
	\centering
	\begin{overpic}[scale=.3]{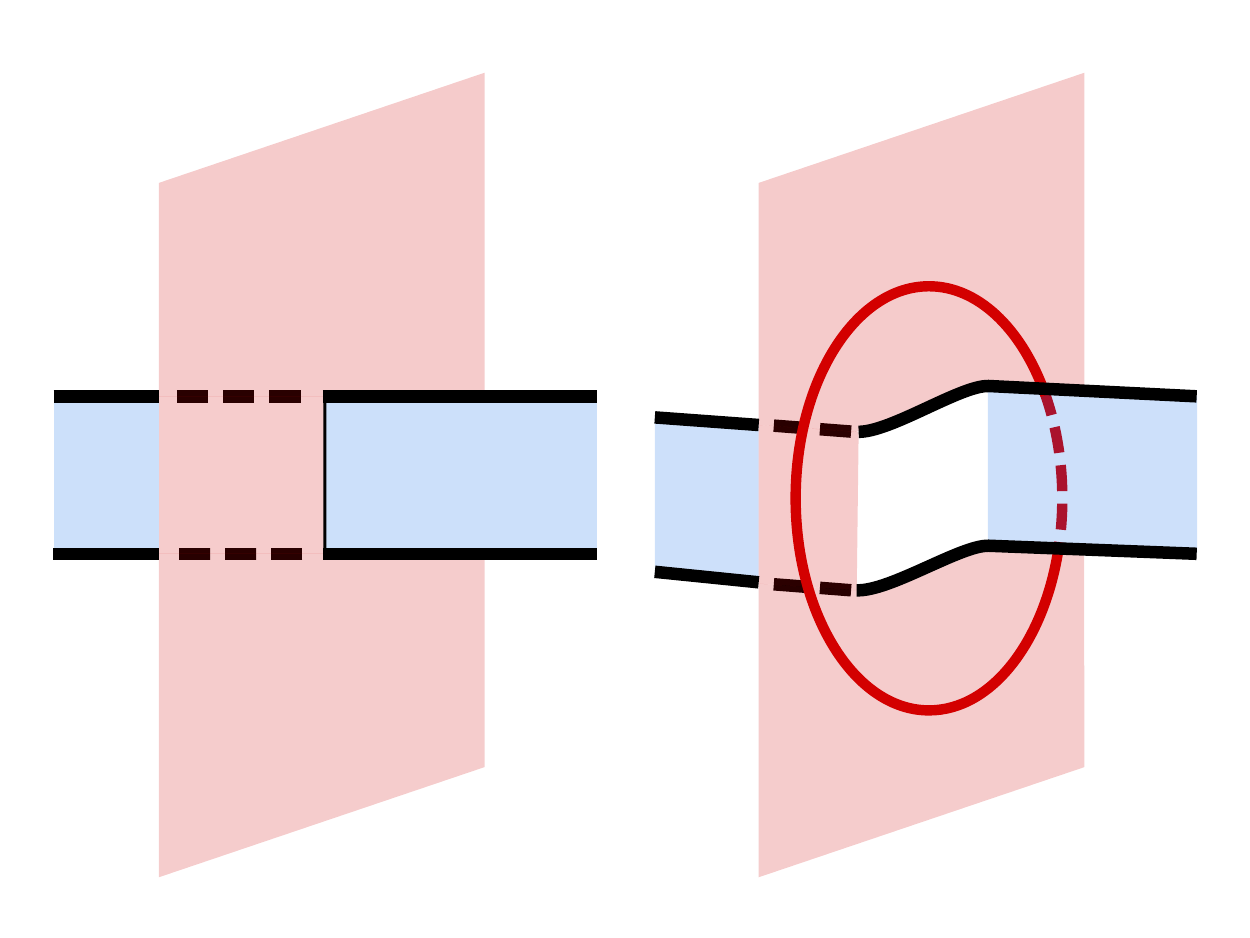}
       
	\end{overpic}
	\caption{Resolving a ribbon disk, and a component of an unlink derivative.}
	\label{fig:cellar}
\end{figure}

Given a ribbon disk, it is straightforward to identify a topological unlink derivative by considering unknots surrounding each ribbon intersection; see \cref{fig:cellar}. The proofs of \cref{cor:slice-benn-sharp} and \cref{cor:CNSrib2} proceed by manipulating the dividing curves of a convex ribbon disk to upgrade this topological unlink derivative to a max-sl transverse unlink, and \cref{thm:BZ1review} then applies. We execute this strategy via a combinatorial analysis of the dividing curves in \cref{sec:ribbon_disk_analysis}.

\subsection{Data on Lagrangian sliceness}
Finally, we provide additional data on Lagrangian sliceness. Regarding \cref{conj:CNS}, Cornwell, Ng, and Sivek recorded a census of Lagrangian slice knots through $14$ crossings in \cite{cornwell2016obstructions}, verifying the conjecture (which we hereafter refer to as \emph{the CNS conjecture} when convenient) through 14 crossings.  Moreover, they narrow down the collection of 15-crossing knots to 48 potentially Lagrangian slice knots.  We explore a subset of these 48 knots that bound star-shaped ribbon disks, excluding two and explicitly verifying the CNS conjecture for the other nine via additional obstruction analysis, cross-referencing databases from \cite{arc11} and~\cite{arc12}, and exhibiting explicit pinch moves on grid diagrams.

\subsection{Organization} In \cref{sec:background} we collect background information primarily on convex surface theory. \cref{sec:ribbon} defines our notion of convexity for ribbon disks and proves \cref{prop:ribbon_convex_exist}, the existence part of \cref{thm:main-ribbon-disk-convex}, while \cref{sec:manipulation} describes the set of manipulations referenced in the same theorem. In \cref{sec:BZ1review} we review the necessary result from \cite{BZ0forthcoming} and sketch the proof in order to make this article nearly self-contained. These ingredients coalesce in \cref{sec:ribbon_disk_analysis} to prove \cref{cor:slice-benn-sharp} and \cref{cor:CNSrib2}. Finally, in \cref{sec:data}, we present the data in support of \cref{conj:CNS}. 

\subsection{Acknowledgments}
The authors would like to thank John Etnyre, Peter Lambert-Cole, Joan Licata, and Gordana Mati\'{c} for discussions and suggestions which helped inspire and progress the framework for convex ribbon intersections. We additionally thank James Hughes, Josh Sabloff, B\"ulent Tosun, and Vera V\'ertesi for interest in our work, Kyle Hayden for a useful correspondence, and Austin Christian for countless conversations, and suggestions regarding grid diagrams. 

\subsection{Statement on AI use} Large language models were not used in the preparation of this article: the ideas and text were conceived and written by the authors.

\section{Background}\label{sec:background}

We begin with some brief background on Legendrian knots, transverse knots, and convex surface theory. Much of this material is pulled in brevity from the background section of \cite{BZ0forthcoming}, with some additional inclusions relevant for the present work. 

\subsection{Legendrian and transverse knots}

We let $\xi_{\mathrm{st}}$ denote the unique tight (positive) contact structure on $S^3$, the (maximally non-integrable) plane field of complex tangencies to $S^3\subset \C^2$, and we study \emph{Legendrian} knots and \emph{transverse} knots in $(S^3, \xi_{\st})$: those that are tangent and (positively) transverse to $\xi_{\st}$, respectively.

Legendrian knots have two classical invariants. The \emph{Thurston-Bennequin number}, denoted $\mathrm{tb}(\cdot)$, is the integral framing induced by the contact structure, and the \emph{rotation number}, denoted $\mathrm{rot}(\cdot)$, is the relative Euler class of the contact structure on a Seifert surface, changing sign according to a choice of orientation of the knot. Both invariants may be computed combinatorially in the \emph{front projection} $\pi_{xz}:\R^3_{(x,y,z)} \to \R^2_{(x,z)}$ of $(\R^3, \xi_{\mathrm{st}}= \ker (dz - y\, dz))$, which uniquely determines the $y$-coordinate via $y =\tfrac{dz}{dx}$ and generically has semi-cubical cusps in place of vertical tangencies:
\begin{align*}
    \mathrm{tb}(\Lambda) &= \mathrm{writhe}(\pi_{xz}(\Lambda))\, - \,\tfrac{1}{2}(\# \text{ of cusps of }\pi_{xz}(\Lambda)), \\
    \mathrm{rot}(\Lambda) &= \tfrac{1}{2}(\# \text{ of down cusps} - \# \text{ of up cusps}).
\end{align*}
Transverse knots $\mathcal{T}\subset (S^3, \xi_{\st})$ only have one classical invariant, the \emph{self-linking number}, defined as follows. Choose a Seifert surface $\Sigma$ and a nonzero section $s: \Sigma \to \xi\mid_{\Sigma}$, which is possible as the vector bundle $\xi\mid_{\Sigma}$ is trivial. Then $\mathrm{sl}(\mathcal{T}):= \mathrm{lk}(\mathcal{T}, s(\mathcal{T}))$. If $\mathcal{T}$ is the positive transverse pushoff of a Legendrian knot $\Lambda$, then $\mathrm{sl}(\mathcal{T}) = \mathrm{tb}(\Lambda) - \mathrm{rot}(\Lambda)$.

\begin{figure}[ht]
	\centering
	\begin{overpic}[scale=.25]{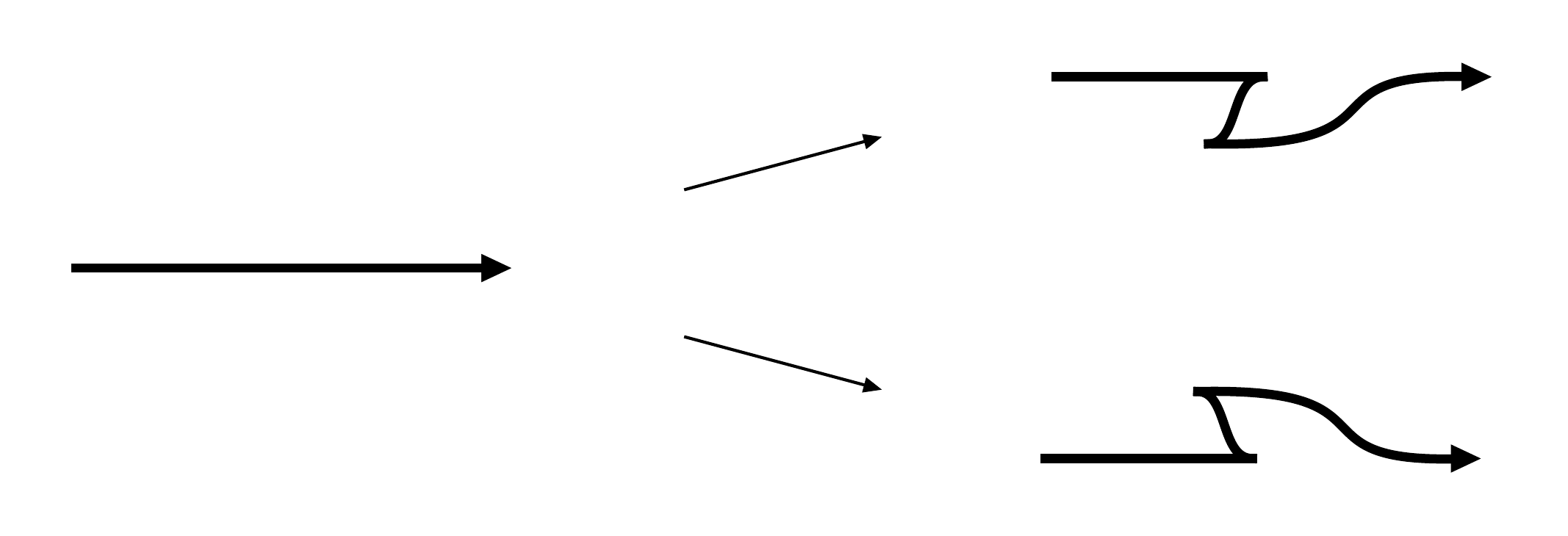}
       \put(48,26){\small $S_+$}
       \put(48,7){\small $S_-$}
	\end{overpic}
	\caption{Positive and negative Legendrian stabilization in the front projection.}
	\label{fig:stabs}
\end{figure}

Given an oriented Legendrian knot $\Lambda \subset (S^3,\xi_{\mathrm{st}})$, we denote by $S_{\pm}(\Lambda)$ the \emph{positive} (resp.\ \emph{negative}) stabilization as in \cref{fig:stabs}. Note that $\mathrm{tb}(S_{\pm}(\Lambda)) = \mathrm{tb}(\Lambda) - 1$ and $\mathrm{rot}(S_{\pm}(\Lambda)) = \mathrm{rot}(\Lambda) \pm 1$. In particular, if $\mathcal{T}, \mathcal{T}_-$ denote the positive transverse pushoffs of $\Lambda, S_-(\Lambda)$, respectively, then $\mathrm{sl}(\mathcal{T}) = \mathrm{sl}(\mathcal{T}_-)$. In fact, $\mathcal{T}$ and $\mathcal{T}'$ are always transversely isotopic \cite[\S 2]{etnyre2005knots}.

Finally, in \cref{sec:data} we consider grid diagrams as presentations of Legendrian knots, using the convention that we rotate a given grid counter-clockwise by $45$ degrees, turning NW and SE corners into cusps, and smoothing NE and SW corners. 

\subsection{Convex surface theory}

Convex surface theory was reviewed more thoroughly in \cite{BZ0forthcoming}, so we point the reader either there or to any of the standard references (e.g.\ \cite{etnyre2004convexlectures}) for more details. For now we record the essential definitions and theorems which are used throughout this work.

\begin{definition}
Let $(Y, \xi)$ be a co-oriented contact $3$-manifold. A \emph{contact vector field} $X$ on $Y$ is one satisfying $\mathcal{L}_X\xi = \xi$, and an oriented surface $\Sigma\subset (Y, \xi)$, compact and possibly with Legendrian boundary, is \emph{convex} if there is a contact vector field $X$ everywhere transverse to $\Sigma$. The \emph{dividing set} is $\Gamma_X := \set{p\in \Sigma}{X_p \in \xi_{p}}$. The \emph{positive} (resp.\ \emph{negative}) region, denoted $R_{\pm}(\Sigma)$, is the set of points $p\in \Sigma$ where $X_p$ is positively (resp.\ negatively) transverse to $\xi_p$.   
\end{definition}

The isotopy class of the multicurve $\Gamma_X$ in $\Sigma$ is independent of the choice of $X$. For this reason, we refer to the dividing set as $\Gamma$ without reference to a specific contact vector field. 

\begin{lemma}[\cite{giroux1991convexite}]\label{lemma:nbd_size}
Let $\Sigma \subset (Y, \xi)$ be a convex surface. Let $\mathcal{U} \subset Y$ be any open neighborhood of $\Sigma$. There is a contact embedding $(\Sigma \times \R_t, \ker \alpha_{\mathrm{inv}})\hookrightarrow (\mathcal{U}, \xi)$ mapping $\Sigma\times \{0\}$ to $\Sigma \subset \mathcal{U}$, such that 
\begin{enumerate}
    \item $\partial_t$ is a strict contact vector field, so that $\Sigma\times \{0\} \subset \Sigma \times \R$ is convex and $\alpha_{\mathrm{inv}}$ is $\R$-invariant, and 
    \item the dividing set of $\Sigma \times \{0\}$ is mapped to the dividing set of $\Sigma\subset \mathcal{U}$. 
\end{enumerate}
Specifically, we may take $\alpha_{\mathrm{inv}}$ to be of the form $\alpha_{\mathrm{inv}} = f\, dt + \beta$ where $f, \beta$ are $t$-invariant, the dividing set is $\Gamma =\{f=0\}$, and $R_{\pm}(\Sigma) = \{\pm f > 0\}$.
\end{lemma}

\begin{theorem}[\cite{giroux1991convexite,honda2000classification}]\label{thm:legendrian_bdry_convex}
Let $\Sigma\subset (Y, \xi)$ be an orientable surface with Legendrian boundary $\partial \Sigma = \Lambda_1 \cup \cdots \cup \Lambda_n$. Assume that $\mathrm{tb}(\Lambda_i)\leq f_i$, where $f_i$ is the smooth framing induced by $\Sigma$. Then there is a perturbation of $\Sigma$ relative to $\partial \Sigma$, $C^0$-small near $\partial \Sigma$ and $C^{\infty}$-small in the interior, such that the resulting surface is convex. Moreover, $\mathrm{tb}(\Lambda_i) = f_i -\frac{1}{2}|\Gamma \cap \Lambda_i|$ and $\mathrm{rot}_{[\Sigma]}(\partial \Sigma) = \chi(R_+(\Sigma)) - \chi(R_-(\Sigma))$. 
\end{theorem}

\begin{figure}[ht]
	\centering
	\begin{overpic}[scale=.3]{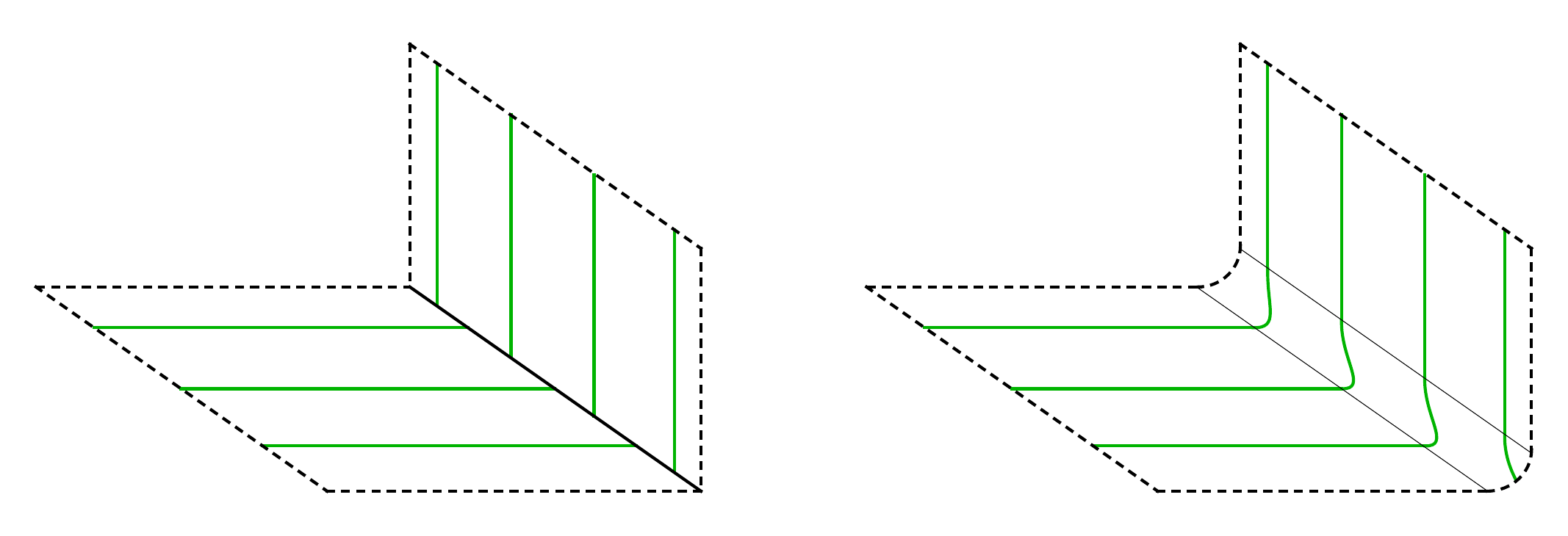}
       
	\end{overpic}
	\caption{Edge-rounding convex surfaces along a Legendrian corner.} 
	\label{fig:edge}
\end{figure}

The edge-rounding principle of \cite[Lemma 3.11]{honda2000classification} describes how to smooth the corners of two convex surfaces meeting along a Legendrian edge. We will state and analyze a local model of this Legendrian edge more explicitly in \cref{sec:ribbon}; for now, we remark that the dividing curves on each sheet are interlacing, and the edge-rounding principle is "up and to the left" as in \cref{fig:edge}.

If $\Sigma$ is convex with dividing set $\Gamma$, we say that a properly embedded graph $C\subset \Sigma$ is \textit{non-isolating with respect to $\Gamma$} if each component of $\Sigma - C$ has nonempty intersection with $\Gamma$.

\begin{theorem}[Legendrian realization principle \cite{kanda1998legendrian,honda2000classification}]\label{thm:LeRP}
Let $\Sigma \subset (Y, \xi)$ be a convex surface with dividing set $\Gamma$. Let $C\subset \Sigma$ be a properly embedded graph transverse to and non-isolating with respect to $\Gamma$. There is an isotopy from the identity $\phi_s:\Sigma \to \mathcal{U}\subset Y$, $s\in [0,1]$, supported in the invariant neighborhood of $\Sigma$, such that
\begin{enumerate}
    \item $\phi_s(\Sigma)$ is graphical (in particular convex) in the invariant neighborhood, and 
    \item $\phi_1(C)$ is Legendrian.
\end{enumerate}
Moreover, wherever $C$ is already Legendrian, the isotopy may be taken to be constant. 
\end{theorem}

Recall that a contact structure is \emph{overtwisted} if there is a Legendrian unknot with $\mathrm{tb} = 0$, and is \emph{tight} otherwise.

\begin{theorem}[Giroux criterion \cite{honda2000classification}]\label{thm:criterion}
If $\Sigma \neq S^2$ is a convex surface in a contact manifold $(Y, \xi)$, then $\Sigma$ has
a tight neighborhood if and only if the dividing set $\Gamma$ has no homotopically trivial curves. When $\Sigma = S^2$, $\Sigma$ has a tight neighborhood if and only if $\Gamma$ has a single connected component.
\end{theorem}

\begin{figure}[ht]
	\centering
	\begin{overpic}[scale=.4]{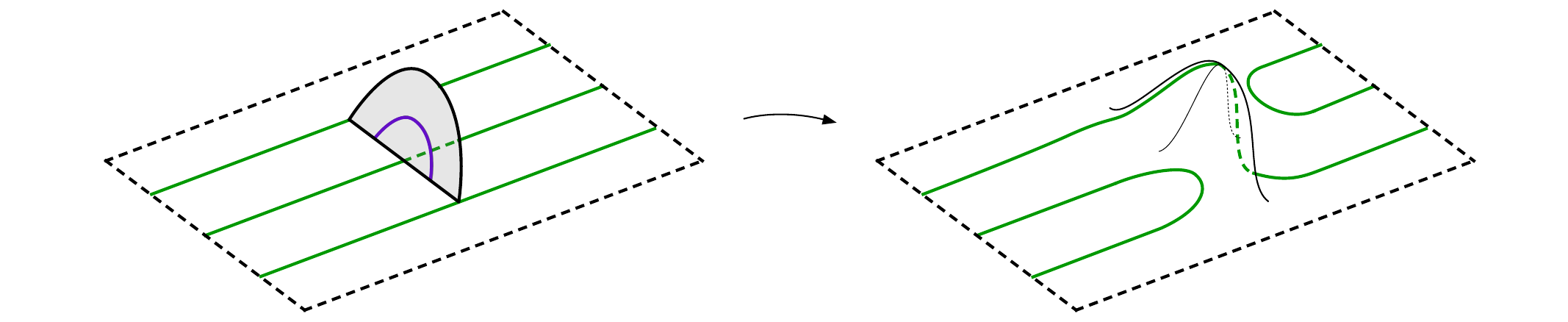}
       \put(10,,1){$\Sigma$}
       \put(43,13){\small \textcolor{darkgreen}{$\Gamma$}}
       \put(60,1){$\Sigma'$}
       \put(93,13){\small \textcolor{darkgreen}{$\Gamma'$}}
	\end{overpic}
	\caption{Bypass attachment.}
	\label{fig:bypassfigure2}
\end{figure}

A \emph{bypass half-disk} for a convex surface $\Sigma$ is topological disk $D$ with piecewise-smooth Legendrian boundary $\partial D = a \cup b$, where $a\subset \Sigma$ is an embedded Legendrian arc transversally intersecting the dividing set thrice, twice at its endpoints, $b$ is properly embedded in the complement of $\Sigma$, the natural smoothing of $a\cup b$ is a $\mathrm{tb}=-1$ unknot, and $D$ is convex with dividing curve given by the purple arc in the middle of \cref{fig:bypassfigure2}. One may isotope $\Sigma$ past $D$ to obtain a new convex surface $\Sigma'$ with dividing set $\Gamma'$ as indicated in \cref{fig:bypassfigure2}; see \cite[Lemma 3.12]{honda2000classification}. The process of isotoping $\Sigma$ past $D$ to obtain $\Sigma'$ is called \emph{attaching a bypass to $\Sigma$}.

\begin{figure}[ht]
	\centering
	\begin{overpic}[scale=.38]{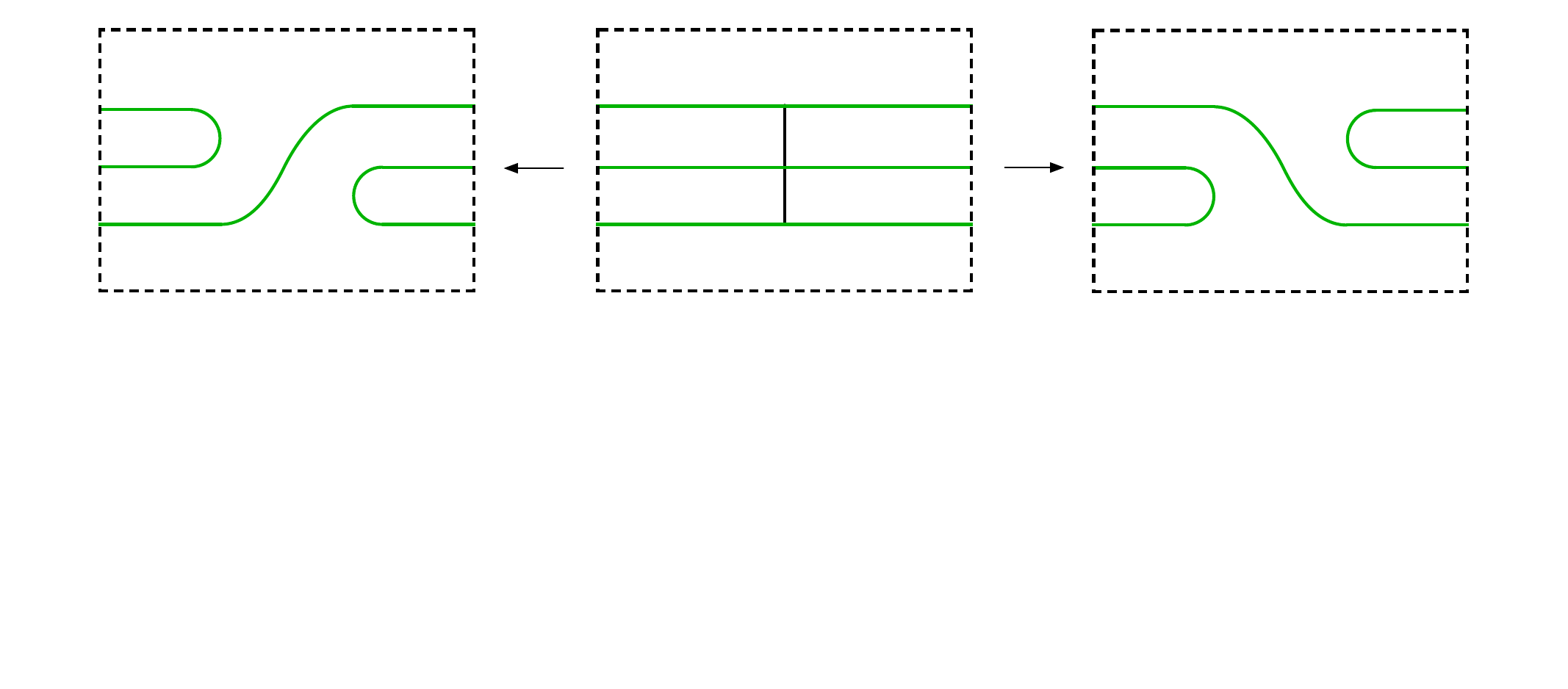}
       \put(10.5,39){\footnotesize from the back}
       \put(73.75,39){\footnotesize from the front}
	\end{overpic}
    \vskip-3cm
	\caption{A bypass attachment from the front or back.}
	\label{fig:frontback}
\end{figure}

\begin{remark}[From the front or back]\label{remark:frontback}
The combinatorial effect of a bypass attachment described above, when viewed from the side of the bypass half disk, is what is known as a bypass \emph{from the front}. If we instead find a bypass half-disk with the same attaching arc below the surface, the effect of isotoping $\Sigma$ past the half-disk results in a combinatorial modification which, when viewed from above, is the reflection of a bypass from the front; see \cref{fig:frontback}. We call this a bypass \emph{from the back}.
\end{remark}

We record one last lemma that will be used frequently in future sections, describing negative stabilization from the convex point of view.

\begin{figure}[ht]
	\centering
	\begin{overpic}[scale=.4]{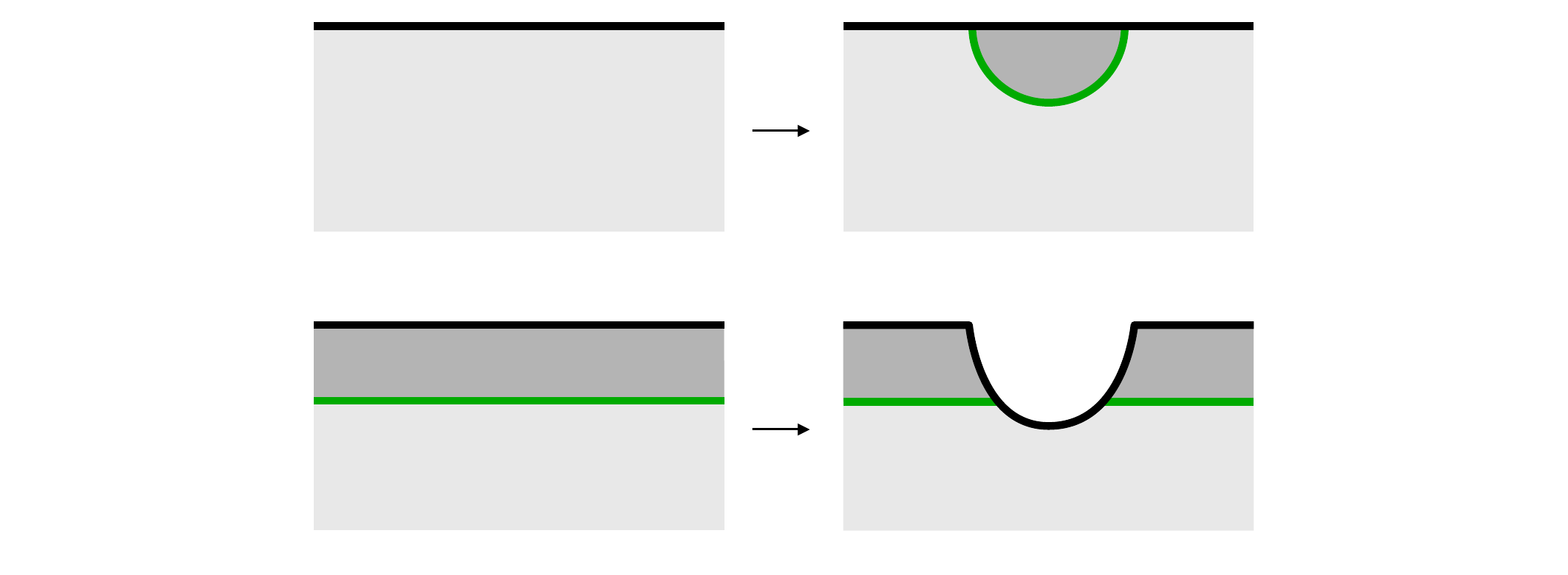}
    
	\end{overpic}
	\caption{Negatively stabilizing a Legendrian via convex surface theory. Negative regions are shaded in dark gray.}
	\label{fig:destab}
\end{figure}

\begin{lemma}\label{lemma:destab}
Let $\Sigma$ be an oriented convex surface with oriented Legendrian boundary $\Lambda$.
\begin{enumerate}
    \item Let $p\in \Lambda$ be a point such that $p\in R_+(\Sigma)$. After performing a negative stabilization of $\Lambda$ near $p$, the resulting surface $\Sigma'$ is convex with an extra boundary parallel dividing curve encasing a negative disk; see the top row of \cref{fig:destab}.
    \item Let $p\in \Lambda$ be a point such that $p\in R_-(\Sigma)$. Let $U\subset\Sigma$ be an open neighborhood of $p$ homeomorphic to a disk and such that $U \cap \Gamma_{\Sigma}$ contains a single arc; assume by the Legendrian realization principle that $\partial U$ is Legendrian. The surface $\Sigma'$ obtained by excising $U$ from $\Sigma$ has boundary Legendrian isotopic to a negative stabilization of $\Lambda$; see the bottom row of \cref{fig:destab}.  
\end{enumerate}
\end{lemma}

\begin{proof}
The first statement (1) is standard; see \cite{etnyre2001knotscontact}. For (2), by Giroux flexibility we may assume that the characteristic foliation $U \cap R_-(\Sigma)$ contains a single (negative) critical point of index $1$; indeed, one can birth a canceling pair and position the index $1$ point accordingly. Removing $U$ thus increases the Euler characteristic of the negative region, hence decreases the rotation number of the boundary (and decreases $\mathrm{tb}$ by $1$). It follows that $\partial \Sigma'$ is a negative stabilization of $\Lambda$. 
\end{proof}

\section{Existence of convex ribbon disks}\label{sec:ribbon}

This section contains the general framework of convex surface theory for ribbon disks. Specifically, we define \emph{normalized ribbon-convexity}, the relevant notion of convexity in this setting.  We show (\cref{prop:ribbon_convex_exist}) that one can perturb certain ribbon disks to satisfy this condition. This gives the first part of \cref{thm:main-ribbon-disk-convex}.

\begin{figure}[ht]
	\centering
	\begin{overpic}[scale=.44]{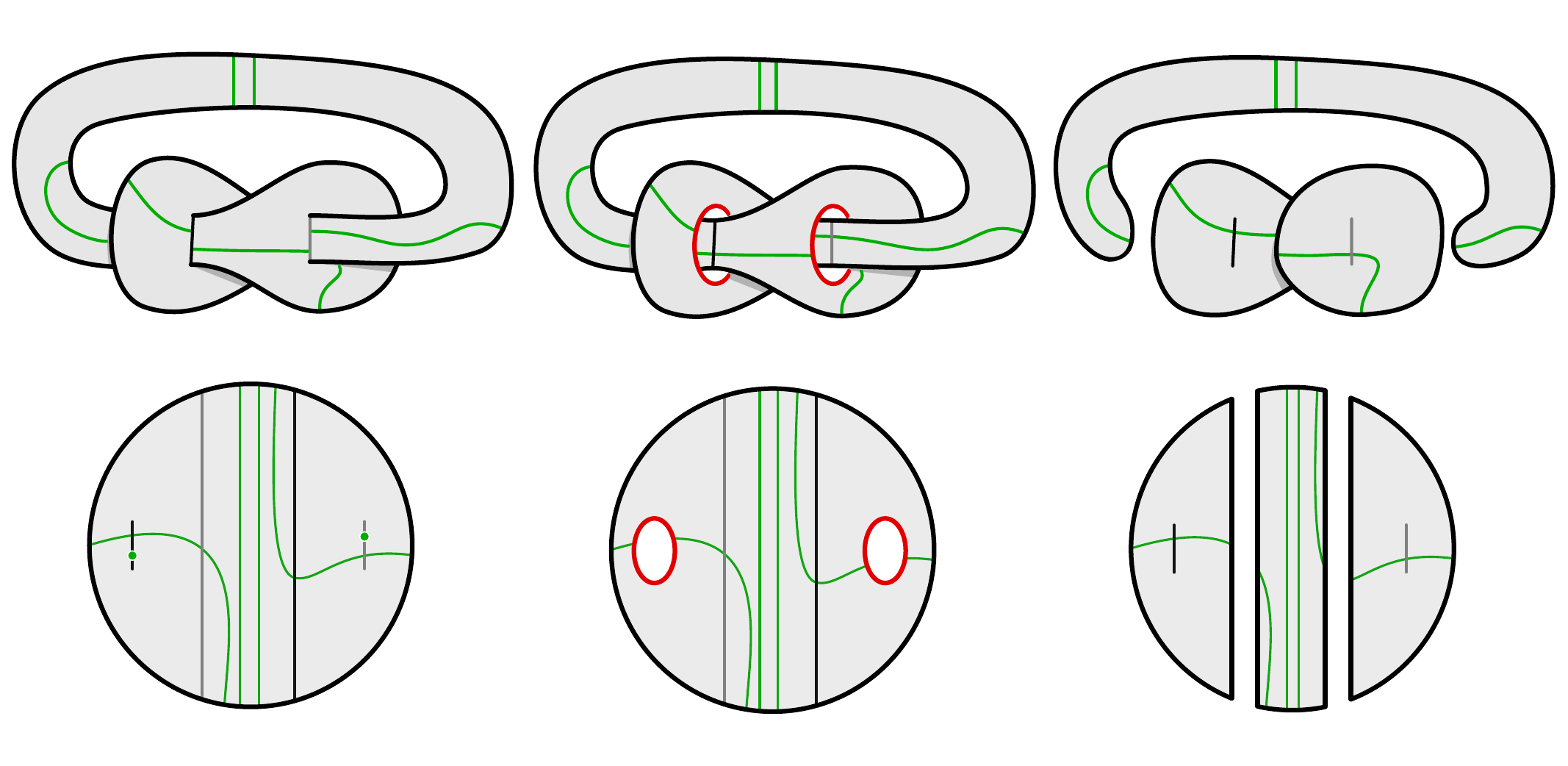}
       \put(1,45){\small$D$}
       \put(33.5,45){\small$D_b$}
       \put(66,45){\small$D_d$}
       \put(5,22){\small $\tilde{D}$}
       \put(7.75,10.5){\footnotesize $\gamma_1^-$}
       \put(17.75,1.75){\footnotesize $\gamma_1^+$}
       \put(12,1.75){\footnotesize \textcolor{gray}{$\gamma_2^+$}}
       \put(14.85,1.75){\footnotesize \textcolor{darkgreen}{$\Gamma_{\tilde{D}}$}}
       \put(22,10.5){\footnotesize \textcolor{gray}{$\gamma_2^-$}}
	\end{overpic}
	\caption{A schematic picture of a (non-normalized) ribbon-convex disk.}
	\label{fig:ribbonconvex}
\end{figure}

Let $D\subset S^3$ be a ribbon disk with $n$ ribbon singularities $\gamma = \{\gamma_1,\dots,\gamma_n\}$. Let $\mathcal{U}(\gamma) = \bigcup_{i=1}^{n} \mathcal{U}(\gamma_i)$ be a sufficiently small open tubular neighborhood of $\gamma$ so that $D \cap \mathcal{U}(\gamma_i) = d_i \cup b_i$, where $\bar{d}_i\subset \overline{\mathcal{U}(\gamma_i)}$ is a properly embedded disk and $\bar{b}_i\subset \overline{\mathcal{U}(\gamma_i)}$ is a properly embedded band with $b_i \cap d_i = \gamma_i$. Define 
\[
D_b := D - \left(\bigcup_{i=1}^{n} d_i\right),\quad  D_d := D - \left(\bigcup_{i=1}^{n} b_i\right), \quad \text{and} \quad \hat{D} := D - \mathcal{U}(\gamma).
\]
Note that $D_b = \hat{D}\,\cup\, \bigcup_{i=1}^n b_i$ is a connected planar surface with $n+1$ boundary components and $D_d = \hat{D} \cup \bigcup_{i=1}^n d_i$ is a collection of $n+1$ disks. See \cref{fig:ribbonconvex}.

\begin{definition}\label{def:normalized}
Let $D \subset (S^3, \xi_{\mathrm{st}})$ be a ribbon disk with Legendrian boundary and ribbon singularities $\gamma = \{\gamma_1,\dots,\gamma_n\}$. We call $D$ \textit{ribbon-convex} if there is a neighborhood $\mathcal{U}(\gamma)$ as above such that 
\begin{enumerate}
\item each ribbon singularity $\gamma_i$ is Legendrian, 
    \item $D_b$ is convex with respect to a contact vector field $X_b$, and 
    \item $D_d$ is convex with respect to a contact vector field $X_d$ such that $X_d = X_b$ away from $\mathcal{U}(\gamma)$.
\end{enumerate}
We call the ribbon-convex structure \textit{normalized} if for all $i\in \{1,\dots,n\}$ and every $x\in \gamma_i$, we have $(X_b)_x \in T_xd_i$ and $(X_d)_x \in T_xb_i$. 
\end{definition}

Let $\tilde{D}$ denote an abstract disk which is the domain of the immersion $\tilde{D} \looparrowright D \subset (S^3, \xi_{\mathrm{st}})$, and for each ribbon singularity $\gamma_i$ let $\gamma_i^{\pm}$ denote the preimages of $\gamma_i$ in $\tilde{D}$, where $\gamma_{i}^-$ is the arc embedded in the interior of $\tilde{D}$. As $D_b$ and $D_d$ are convex with respect to vector fields that agree outside of $\mathcal{U}(\gamma)$, the dividing sets $\Gamma_{b_i}$ and $\Gamma_{d_i}$ glue to the dividing set $\Gamma_{\hat{D}}$ along $\partial \hat{D}$. We let 
\[
\Gamma_{D} := \Gamma_{\hat{D}} \cup\,  \bigcup_{i=1}^n \, (\Gamma_{b_i} \cup \Gamma_{d_i}) \subset D \subset S^3.
\]
Likewise, we let $\Gamma_{\tilde{D}}$ denote the preimage of $\Gamma_{D}$ under the ribbon immersion. Generically, $\Gamma_{\tilde{D}}\subset \tilde{D}$ is an embedded multicurve and a collection of points recording where $\Gamma_{D}\pitchfork D$. If the ribbon-convex structure is normalized, we can say more.

\begin{figure}[ht]
	\centering
	\begin{overpic}[scale=.4]{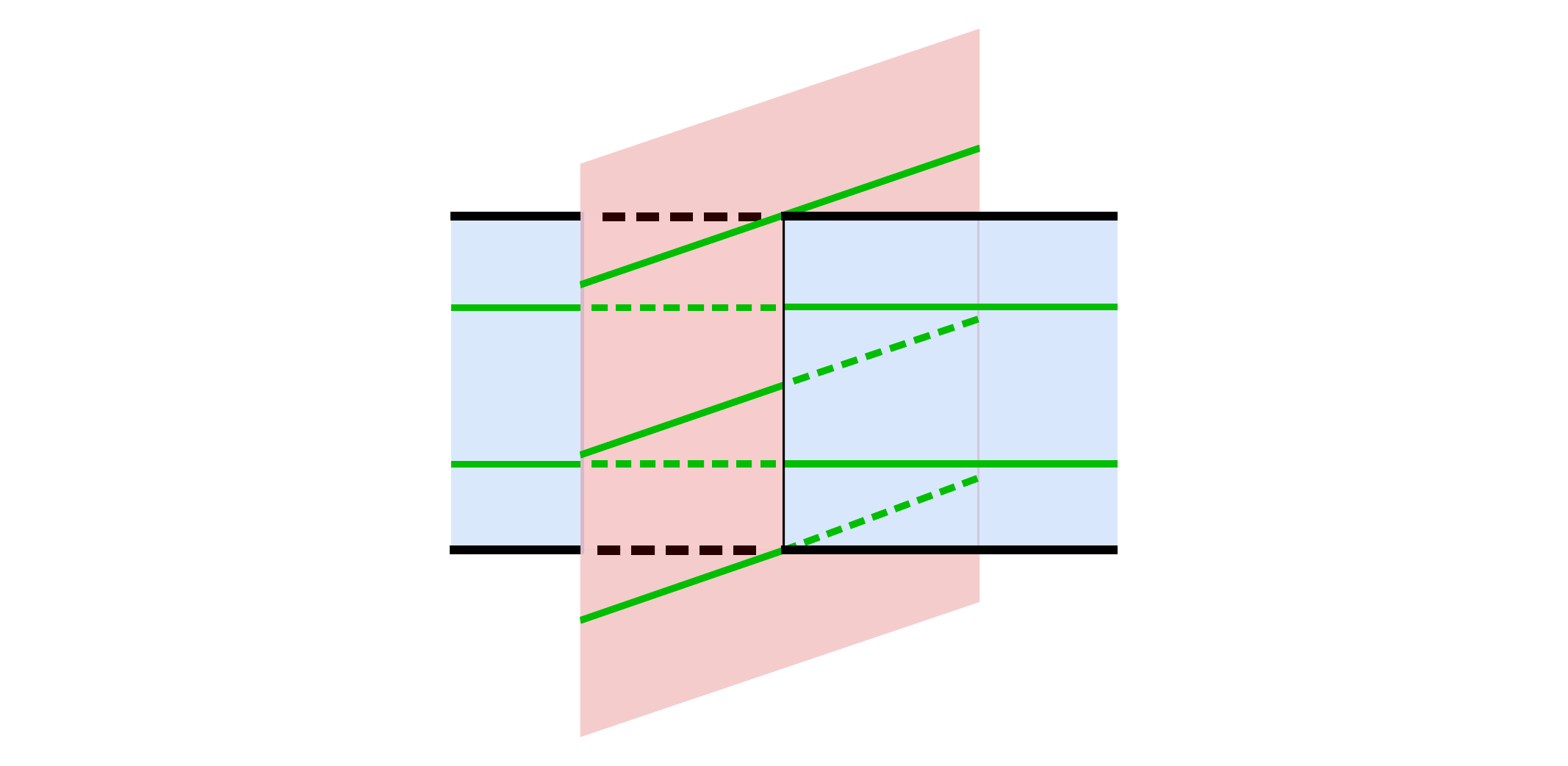}
        \put(49,11){$\gamma_i$}
        \put(72.5,12){\textcolor{black}{$\Lambda$}}

        \put(25,32){$b_i$}
        \put(42,43){$d_i$}
	\end{overpic}
	\caption{The normalized dividing sets $\Gamma_{b_i}$ and $\Gamma_{d_i}$ are interleaved near $\gamma_i$.  The case pictured is $k_i = 2$.}
	\label{fig:normalized}
\end{figure}

\begin{lemma}\label{lemma:normalized_disk}
If $D$ is normalized ribbon-convex, then $\Gamma_{\tilde{D}} = \Gamma_0 \cup \{p_1, \dots, p_m\}$ where $\Gamma_0$ is a properly embedded multicurve disjoint from the points $p_1, \dots, p_m$. Moreover, for each $i\in \{1,\dots, n\}$,
\begin{enumerate}
    \item $\Gamma_0$ intersects both endpoints of $\gamma_i^-$ transversally, 
    \item $\Gamma_0$ intersects the interior of $\gamma_i^+$ transversally $k_i \geq 1$ times, and 
    \item along $\gamma_i^{\pm}$, transverse intersections with $\Gamma_0$ alternate with points in $\{p_1, \dots, p_m\}$; consequently $\Gamma_0$ intersects the interior of $\gamma_i^-$ transversally $k_i - 1$ times. 
\end{enumerate}
In particular, the dividing curves are 
"interleaved" near a ribbon singularity, as indicated in \cref{fig:normalized}.
\end{lemma}

\begin{remark}
Though we specialize to disks in the present work, we point out that both \cref{def:normalized} and \cref{lemma:normalized_disk} hold in the more general setting of ribbon surfaces of any genus.     
\end{remark}

\begin{proof}
First we observe (1). Since $D$ is normalized convex, the contact vector field $X_d$ witnessing convexity of each $d_i$ satisfies $(X_d)_{\partial \gamma_i}\in T_{\partial \gamma_i}b_i$. Since both $\Lambda$ and $\gamma_i$ are Legendrian, $\xi$ is tangent to $b_i$ at $\partial \gamma_i$. This implies that $(X_d)_{\partial \gamma_i} \in \xi$, hence the dividing set $\Gamma_{d_i}$ passes through $\partial \gamma_i$, giving (1).

Next we prove (2). Suppose for the sake of contradiction that $\Gamma_0$ does not intersect the interior of $\gamma_i^+$. Reinterpreted, this means $\gamma_i \cap \Gamma_{b_i} = \emptyset$, and in particular there is no twisting of the contact structure from one endpoint of $\gamma_i$ to the other relative to $b_i$. Since $\Lambda$ and $\gamma_i$ are Legendrian, we may thus identify a tubular neighborhood of $\gamma_i$ with 
\[
\left([-\delta, 1+\delta]_q \times [-\ve,\ve]_p \times \R_t, \, \ker(dt - p\, dq)\right)
\]
where $\delta, \ve > 0$ are small, $\gamma_i = [0,1]_q \times \{p=0, t=0\}$, $b_i = [0,1]_q \times [-\ve,\ve]_p \times \{t=0\}$, $\Lambda = \{q=0,1\}\times [-\ve,\ve]_p \times \{t=0\}$, and $X_b = \partial_t$. (We have not yet identified $d_i$ or $X_d$ in the model. In fact, we will reach a contradiction by arguing that no such pair can exist.) 

Since $D$ is normalized ribbon-convex, there is a contact vector field $X_d$ transverse to $d_i$ such that for all $0\leq q \leq 1$, $(X_{d})_{(q,0,0)}$ is tangent to $b_i$. In particular, if we write $X_d = f\, \partial_q + g\, \partial_p + h\, \partial_t$, then up to orientation we have 
\begin{align}
    g(q,0,0) &> 0, \label{eq:cond1} \\
    h(q,0,0) &= 0.\label{eq:cond2}
\end{align}
The former condition ensures transversality to $d_i$ along $\gamma_i$, while the latter ensures tangency to $b_i$ along $\gamma_i$. Note that 
\begin{align*}
    \mathcal{L}_{X_d}(dt - p\, dq) &= \iota_{X_d}(dq \wedge dp) + d\iota_{X_d}(dt - p\, dq) \\
    &= dh - p\, df - g\, dq.
\end{align*}
Since $X_d$ is a contact vector field, $\ker \mathcal{L}_{X_d}(dt - p\, dq) \supset \ker(dt - p\, dq)$. In particular, along $\{p=0\}$ we have $\ker(dh - g\, dq)\mid_{p=0} \supset \ker dt$. Writing $dh = \frac{\partial h}{\partial q}\, dq + \frac{\partial h}{\partial p}\, dp + \frac{\partial h}{\partial t}\, dt$, 
\[
dh - g\, dq = \left(\frac{\partial h}{\partial q} - g\right)\, dq + \frac{\partial h}{\partial p}\, dp + \frac{\partial h}{\partial t}\, dt,
\]
so along $\{p=0\}$ we require $\frac{\partial h}{\partial q} = g$. Recalling from \eqref{eq:cond1} that $g(q,0,0) > 0$, it follows that $\frac{\partial h}{\partial q}(q,0,0) > 0$. This contradicts the condition \eqref{eq:cond2} that $h(q,0,0) = 0$. Therefore, we conclude $\gamma_i \cap \Gamma_{b_i} \neq \emptyset$ and thus $k_i := |\Gamma_0 \cap \gamma_i^+| > 0$, giving (2).

Because $\gamma_i \cap \Gamma_{b_i} \neq \emptyset$ and $\gamma_i$ and $\Lambda$ are Legendrian, up to a small contact isotopy we may take a sufficiently small standard tubular neighborhood of $\gamma_i$ to be contactomorphic to 
\[
\overline{\mathcal{U}}_i' := \left([-\delta,1+\delta]_z \times \{x^2 + y^2 \leq \ve^2\},\, \ker\left(\sin(\pi k_i z)\, dx + \cos(\pi k_i z)\, dy\right)\right) 
\]
where $\delta, \ve > 0$ are small, $\gamma_i = [0,1] \times \{(0,0)\}$, $b_i = [0,1] \times \{y=0\}$, and the contact vector field witnessing convexity of $D_b$ restricts to the neighborhood as $X_b = \partial_y$. The normalized assumption allows us to arrange $d_i = [-\delta,1+\delta] \times \{x = 0\}$, which is transverse to the contact vector field $X_d = \pm \partial_x$. This is essentially the local model for edge-rounding \cite[Lemma 3.11]{honda2000classification} for two convex surfaces meeting transversally along a closed Legendrian, modified to account for the fact that $\gamma_i$ has boundary with endpoints on the Legendrian $\Lambda$. By computing when $\partial_y \in \xi$ and $\pm \partial_x\in \xi$ respectively, we have 
\begin{align*}
    \Gamma_{b_i} &=\left\{z = \frac{2\ell - 1}{2k_i}\right\}_{\ell=1,\dots,k_i} \times [-\ve,\ve]_x \times \{y=0\},\\
    \Gamma_{d_i} &= \left\{z = \frac{2\ell}{2k_i}\right\}_{\ell=0,\dots,k_i} \times \{x=0\} \times [-\ve,\ve]_y.  
\end{align*}
See \cref{fig:normalized}. In particular, along $\gamma_i$ the dividing set intersections of $b_i$ and $d_i$ alternate. This gives (3), completing the proof of the lemma. 
\end{proof}

\begin{proposition}\label{prop:ribbon_convex_exist}
Let $D\subset (S^3,\xi_{\mathrm{st}})$ be an oriented ribbon disk bounding a Legendrian knot $\Lambda$, with $\Lambda$ oriented as $\partial D$. There is an isotopy relative to $\Lambda$ so that $D$ is normalized ribbon-convex. Moreover, $\mathrm{tb}(\Lambda) = -\tfrac{1}{2}|\Gamma_0 \cap \Lambda|$ and $\mathrm{rot}(\Lambda) = \chi(R_+(\tilde{D})) - \chi(R_-(\tilde{D}))$.
\end{proposition}

\begin{proof}
Perform a smooth isotopy of $D$ relative to $\Lambda$ to eliminate all boundary parallel ribbon singularities. Next, let $\mathcal{U}(\gamma) = \bigcup_{i=1}^{n} \mathcal{U}(\gamma_i)$ be a sufficiently small open tubular neighborhood of the ribbon singularities as earlier so that $D \cap \mathcal{U}(\gamma_i) = d_i \cup b_i$. By $C^{\infty}$-density of convex surfaces, we may assume that each $S^2_i := \partial \overline{\mathcal{U}(\gamma_i)}$ is convex. The dividing set $\Gamma_{S^2_i}$ is an embedded circle and the intersection $b_i \cap S^2_i$ consists of two arcs. After a smooth isotopy of $D_d$ supported in a neighborhood of $S^2_i - (b_i\cap S^2_i)$, we may assume that $\partial d_i \cap \Gamma_{S^2_i}\neq \emptyset$. By the Legendrian realization principle, a further isotopy ensures that $L_i := \partial d_i$ is Legendrian.    

\vspace{2mm}
\textsc{Step 1:} \textit{Perturb $D_b$ to be convex.}
\vspace{2mm}

Next, we argue that we can perturb $D_b$ to be convex. At this stage, $D_b$ is a connected planar surface with Legendrian boundary $\partial D_b = \Lambda \cup L_1\cup \cdots \cup L_n$. Observe that each (unknotted) $L_i$ is $0$-framed by $D_b$, as the $D_b$-framing agrees with the $\bar{d}_i$-framing. Additionally, we claim that $\Lambda$ is $0$-framed by $D_b$. To see this, perform a "cellar-door" resolution of the ribbon intersections of $D$ to construct a Seifert surface for $\Lambda$. Compared to $D_b$ (the middle of \cref{fig:ribbonA}), the cellar-door resolution induces a full twist of either sign along the two strands of $\Lambda$ for each ribbon singularity, netting no change to the induced framing; see the right side of \cref{fig:ribbonA}. Consequently, the $D_b$-framing agrees with the Seifert framing. Since each $L_i$ is an unknot in a tight contact manifold, we have $\mathrm{tb}(L_i) \leq -1$, and since $\Lambda$ is ribbon, the slice-Bennequin inequality implies $\mathrm{tb}(\Lambda) \leq -1$; it follows that the contact framing of each component of $\partial D_b$ relative to $D_b$ is negative. By \cref{thm:legendrian_bdry_convex}, there is a perturbation of $D_b$ relative to $\partial D_b$ which is convex with respect to some contact vector field $X_b$. 

\begin{figure}[ht]
	\centering
    \vskip-1cm
	\begin{overpic}[scale=.4]{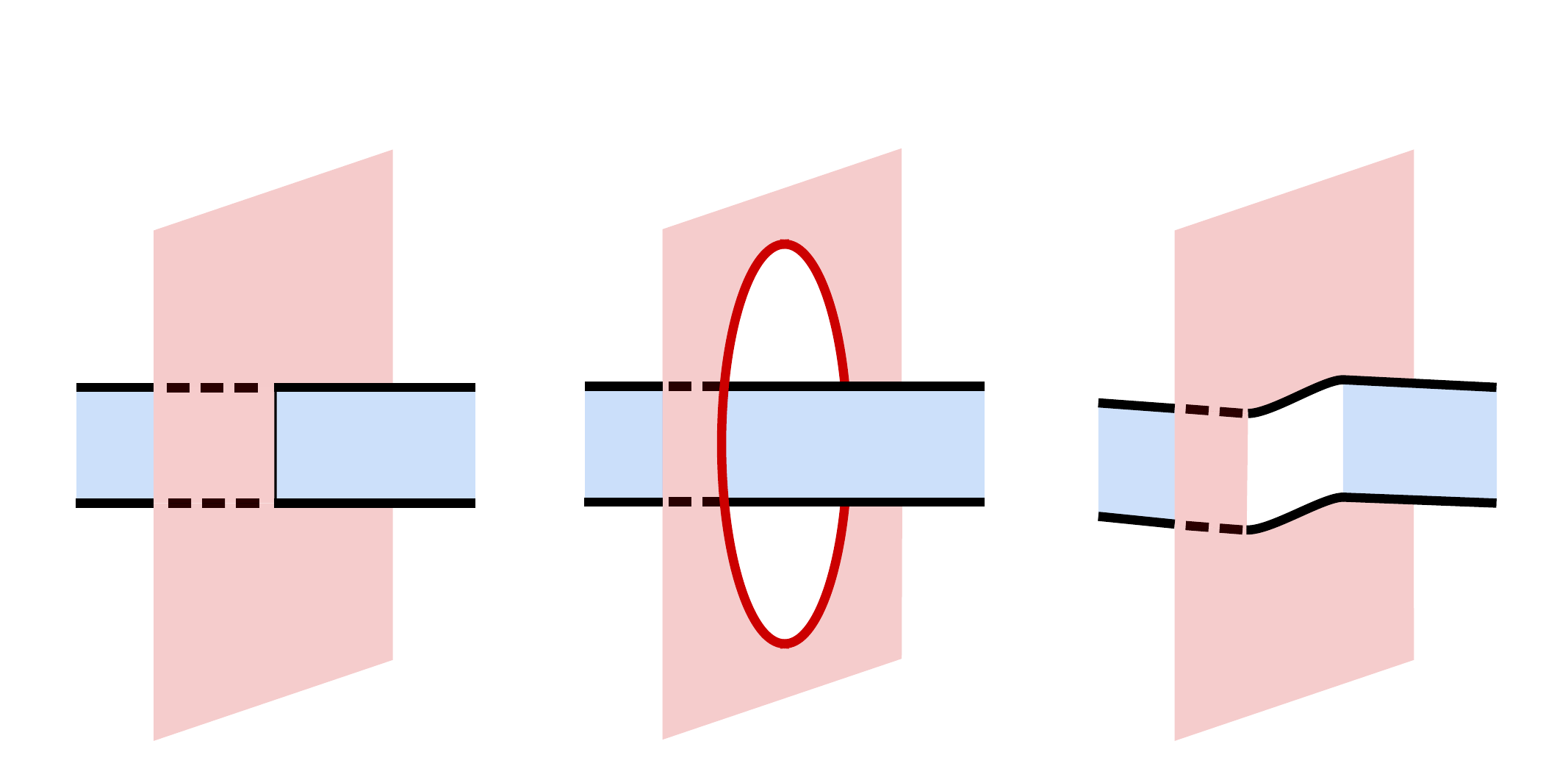}
        \put(16,13){$\gamma_i$}
        \put(48.5,11){\textcolor{darkred}{$L_i$}}

        \put(6,35){$D$}
        \put(38,35){$D_b$}
	\end{overpic}
	\caption{A "cellar-door" resolution.}
	\label{fig:ribbonA}
\end{figure}

\vspace{2mm}
\textsc{Step 2:} \textit{Legendrian realize all ribbon singularities in $D_b$.}
\vspace{2mm}

Next, we argue that each ribbon singularity $\gamma_i \subset b_i \subset D_b$ may be Legendrian realized. As $\gamma_i$ is not a boundary parallel ribbon singularity, it is not $\Lambda$-parallel in $D_b$. In particular, $\gamma_i$ is non-isolating (see the discussion before \cref{thm:LeRP}) with respect to the boundary components $L_1, \dots, L_n$. For each $1\leq j\leq n$, $L_j$ is $0$-framed by $D_b$ and $\mathrm{tb}(L_j)\leq -1$, so $\Gamma_{D_b}\cap L_j \neq \emptyset$, and consequently $\gamma_i$ is non-isolating with respect to $\Gamma_{D_b}$. By the Legendrian realization principle, an isotopy of $D_b$ supported in $\mathcal{U}(\gamma_i)$ witnesses $\gamma_i$ as a Legendrian arc.

If $\mathrm{tb}(\Lambda) = -1$, we claim that each $\gamma_i$ contains an essential intersection with $\Gamma_{D_b}$. Suppose for the sake of contradiction that $\gamma_i\cap \Gamma_{D_b} = \emptyset$. As $\gamma_i$ is separating and not $\Lambda$-parallel, cutting $D_b$ along $\gamma_i$ produces two connected planar surfaces $D_b^{\pm}$, each admitting a (smoothed) Legendrian boundary component $\Lambda^{\pm}$; see \cref{fig:ribbonB}. Since $|\Gamma_{D_b}\cap \Lambda|=2$ and $\Gamma_{D_b}\cap \gamma_i = \emptyset$, it is either the case that $|\Gamma_{D_b^{\pm}}\cap \Lambda^{\pm}| = 1$, or (up to relabeling) $|\Gamma_{D_b^{+}}\cap \Lambda^{+}| = 2$ and $|\Gamma_{D_b^{-}}\cap \Lambda^{-}| = 0$. The former case is not possible, as orientability of $D_b^{\pm}$ implies evenness of $|\Gamma_{D_b^{\pm}}\cap \Lambda^{\pm}|$. In the latter case, we claim that a contradiction arises from $|\Gamma_{D_b^{-}}\cap \Lambda^{-}| = 0$. Note that $\Lambda^-$ is given by $0$-framed topological handleslides of the unknotted components $L_i$ which are bounded by $D_b^-$, and $0$-framed handleslides preserve ribbonness. As $\Lambda^-$ is a ribbon knot, the slice-Bennequin inequality implies $\mathrm{tb}(\Lambda^-) \leq -1$, and in particular the dividing set of any convex Seifert surface bounding $\Lambda^-$ must intersect $\Lambda^-$. However, $D_b^-$ can be extended to a convex Seifert surface by filling in the unknotted (and unlinked) boundary components with disks and performing cellar-door resolutions when necessary. In particular, $|\Gamma_{D_b^{-}}\cap \Lambda^{-}| \neq 0$, reaching a contradiction. We conclude that $\gamma_i\cap \Gamma_{D_b} \neq \emptyset$ as desired.

\begin{figure}[ht]
	\centering
	\begin{overpic}[scale=.36]{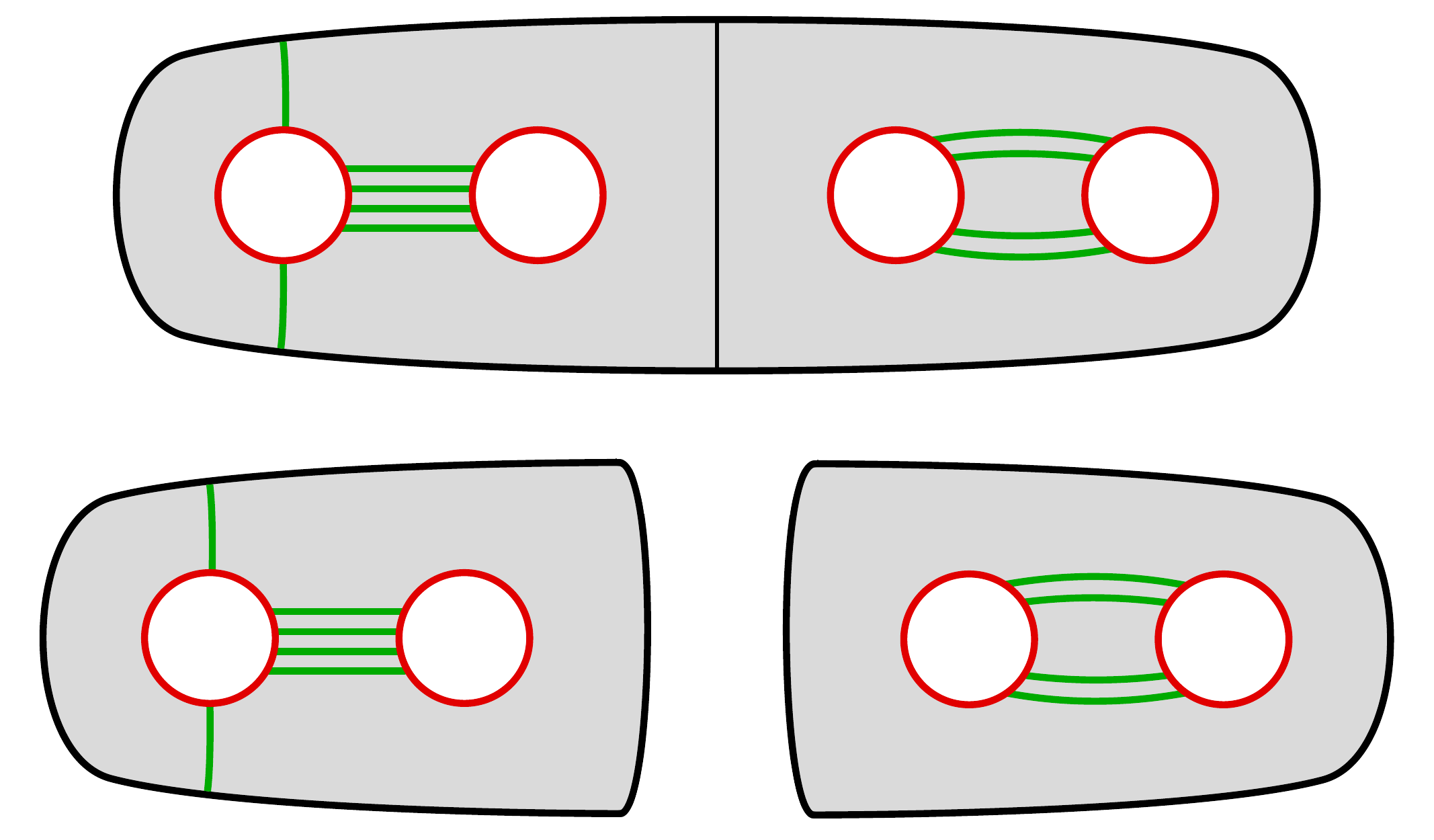}
        \put(5,53){$D_b$}
        \put(3,26){$D_b^+$}
        \put(94,25){$D_b^-$}
        \put(51,36){$\gamma_i$}
	\end{overpic}
	\caption{A dividing set configuration on $D_b$ which is not possible.}
	\label{fig:ribbonB}
\end{figure}

If $\mathrm{tb}(\Lambda) < -1$, it is possible that $\gamma_i\cap \Gamma_{D_b} = \emptyset$. In this case, we Legendrian realize an isotopic copy of $\gamma_i$ obtained by performing a finger move on any nearby dividing curve (which is possible, since $\gamma_i$ is non-isolating). 

\vspace{2mm}
\textsc{Step 3:} \textit{Perturb each $d_i$ so that $D$ is normalized ribbon-convex.}
\vspace{2mm}

Summarizing so far, we have perturbed $D_b$ to be convex, so that for each $i$ we have $|\gamma_i\cap \Gamma_{D_b}| = k_i > 0$. We may identify a local neighborhood $\mathcal{U}_i' \subset \mathcal{U}(\gamma_i)$ of $\gamma_i\subset D_b$ from \cref{lemma:normalized_disk} and assume after a further isotopy of $D$ that $D_d$ also agrees with the model description. Since $\partial \mathcal{U}_i' \cap d_i\mid_{\mathcal{U}'_i}$ intersects $\Gamma_{d_i\mid_{\mathcal{U}'_i}}$, we may further perturb $D$ away from $\gamma_i$ to ensure that $d_i\mid_{\mathcal{U}'_i}$ is Legendrian. To finally perturb all of $d_i$ to be convex, it remains to convexify the (unknotted) annulus $d_i \setminus \mathcal{U}'_i$, which now has unknotted and unlinked Legendrian boundary. This is possible from \cref{thm:legendrian_bdry_convex}.

\vspace{2mm}
\textsc{Step 4:} \textit{Compute $\mathrm{tb}$ and $\mathrm{rot}$.}
\vspace{2mm}

It remains to prove the "moreover" statement of the proposition. The fact that $\mathrm{tb}(\Lambda) = -\tfrac{1}{2}|\Gamma_0 \cap \Lambda|$ follows from the earlier observation that $D_b$, and hence $D$, determines the Seifert framing of $\Lambda$. 

For the rotation number, assume an orientation of $D$. We consider the oriented link $\Lambda \cup L_1 \cup \cdots \cup L_n = \partial D_b$, where $D_b$ inherits the orientation of $D$. Viewing $[D_b] \in H_2(S^3, \Lambda \cup L_1 \cup \cdots \cup L_n)$, we have 
\[
\mathrm{rot}_{[D_b]}(\Lambda \cup L_1 \cup \cdots \cup L_n) = \langle e(\xi), [D_b] \rangle = \chi(R_+(D_b)) - \chi(R_-(D_b)).
\]
Thus 
\[
\mathrm{rot}(\Lambda) = \chi(R_+(D_b)) - \chi(R_-(D_b)) - \sum_{i=1}^n \mathrm{rot}(L_i).
\]
Note that $-L_i = \partial d_i$, where again $d_i$ inherits the orientation of $D$. Thus, by additivity of Euler characteristic, we have  
\begin{align*}
\mathrm{rot}(\Lambda) &= \chi(R_+(D_b)) - \chi(R_-(D_b)) + \sum_{i=1}^n \mathrm{rot}(-L_i) \\
&= \chi(R_+(D_b)) - \chi(R_-(D_b)) + \sum_{i=1}^n \chi(R_+(d_i)) - \chi(R_-(d_i)) \\
&= \chi(R_+(\tilde{D})) - \chi(R_-(\tilde{D}))
\end{align*}
as desired.
\end{proof}

\section{Manipulation of convex ribbon disks}\label{sec:manipulation}

\begin{figure}[ht]
	\centering
    \vskip-0.75cm
	\begin{overpic}[scale=.4]{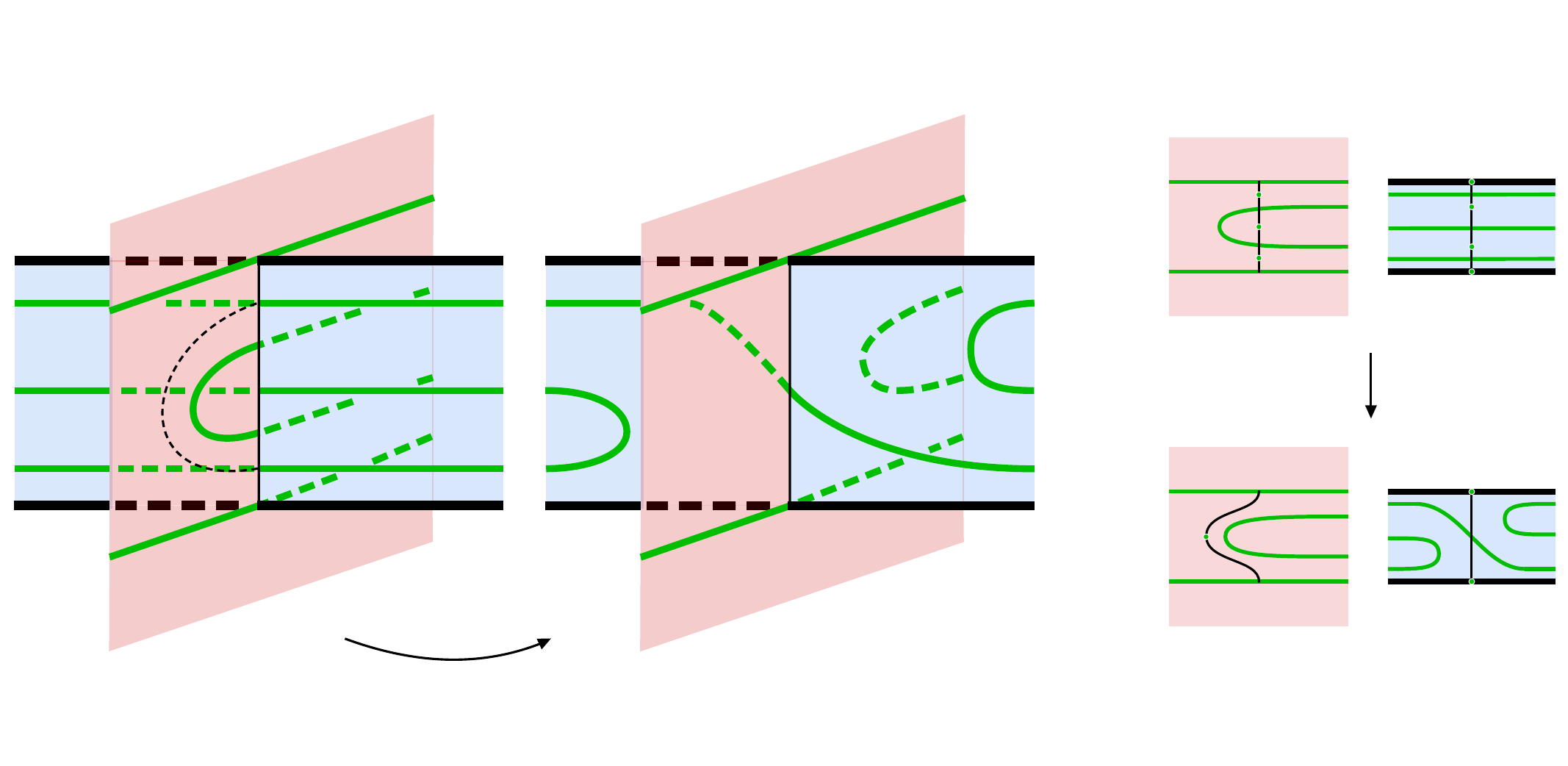}
       
	\end{overpic}
    \vskip-0.75cm
	\caption{The statement of \cref{lemma:bigon_removal} in the $-$ case, removing a $\Gamma_0$-$\gamma_i^-$ bigon.}
	\label{fig:normalized_bypass2}
\end{figure}

\begin{figure}[ht]
	\centering
    \vskip-0.75cm
	\begin{overpic}[scale=.4]{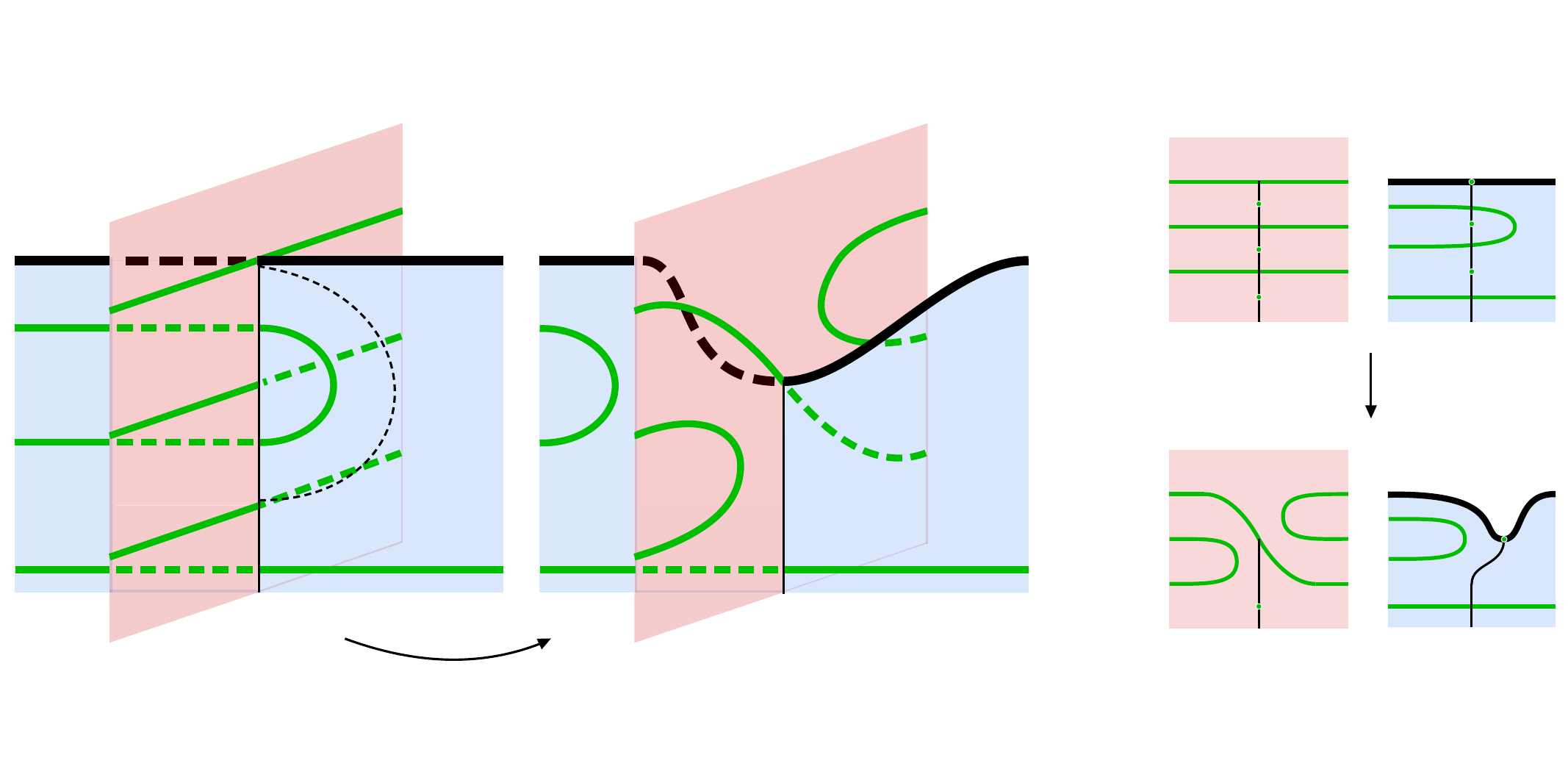}
       
	\end{overpic}
    \vskip-0.75cm
	\caption{The statement of \cref{lemma:bigon_removal} in the more complicated $+$ case, when one of $p_{j_1},p_{j_3}$ is an endpoint of $\gamma_i^+$. The left part of the figure depicts $D$, while the right part depicts $\tilde{D}$.}
	\label{fig:normalized_bypass}
\end{figure}

Here we establish several lemmas which describe how to manipulate the dividing curves of normalized ribbon-convex disks: bigon and half-bigon creation, removal, exchange, and lollipop reduction. Many of these preserve the Legendrian isotopy type of the boundary, while others do not. We begin by removing bigons adjacent to $\gamma_i^-$ (\cref{fig:normalized_bypass2}) and $\gamma_i^+$ (\cref{fig:normalized_bypass}).

\begin{lemma}[Bigon removal]\label{lemma:bigon_removal}
Let $D$ be a normalized ribbon-convex disk with boundary $\Lambda$ and dividing set $\Gamma_{0} \cup \{p_1, \dots, p_m\} \subset \tilde{D}$. Suppose $B\subset \tilde{D}$ is a bigon such that $\partial B = \partial_{\gamma_i^{\pm}}B \cup \partial_{\Gamma_0}B$, where $\partial_{\gamma_i^{\pm}}B$ is subarc in the interior of $\gamma_i^{\pm}$, $\partial_{\Gamma_0}B$ is a subarc of $\Gamma_0$, and $B$ intersects no other ribbon singularities or components of $\Gamma_0$. In the $+$ case, we ask for at least one other intersection of $\Gamma_0$ with $\gamma_i^+$ beyond those adjacent to $B$. Then there is an isotopy of $D$ to a normalized ribbon-convex disk with boundary Legendrian isotopic to $\Lambda$ which modifies the dividing set $\Gamma_{\tilde{D}}$ as follows: 
\begin{enumerate}
    \item The arc $\gamma_i^{\pm}$ is isotoped across $B$, reducing the number of intersections with $\Gamma_0$ by $2$ and consolidating the three interlacing intersection points $p_{j_1}, p_{j_2}, p_{j_3}$ to one point. 
    \item The multicurve $\Gamma_0$ changes by a bypass attachment whose attaching arc is the arc contained in $\gamma_i^{\mp}$ with endpoints corresponding to the points $p_{j_1}, p_{j_3}$. In particular, the number of intersections of $\gamma_i^{\mp}$ with $\Gamma_0$ reduces by $2$. 
\end{enumerate}
In the $+$ case, if either $p_{j_1}$ or $p_{j_3}$ is an endpoint of $\gamma_i^+$, then the remaining intersection point is an endpoint of $\gamma_i^+$ and the ribbon singularity $\gamma_i^-$ terminates on the diagonal of the bypass-affected dividing curve; see \cref{fig:normalized_bypass}.
\end{lemma}

\begin{remark}
Depending on the relative orientations of $b_i$ and $d_i$, the effect of the bypass attachment may be a bypass from the front or from the back when viewed in $\tilde{D}$; see \cref{remark:frontback}.
\end{remark}

\begin{proof}
Given a bigon $B$ as in the statement of the lemma with $\partial B = \partial_{\gamma_i^{\pm}}B \cup \partial_{\Gamma_0}B$, (3) of \cref{lemma:normalized_disk} implies that there are three points $p_{j_1}, p_{j_2}, p_{j_3}\subset \gamma_{i}^{\pm}$, alternating along $\gamma_i^{\pm}$ with the endpoints of $\partial_{\Gamma_0}B$, such that $p_{j_2}\in \partial_{\gamma_i^{\pm}}B$. If $B$ is incident to $\gamma_i^-$, normalization forces both $p_{j_1},p_{j_3}$ to also be in the interior of $\gamma_i^-$, whereas if $B$ is incident to $\gamma_i^+$ it is possible that one is an endpoint of $\gamma_i^+$. 

In either case, let $c\subset \tilde{D}$ be an arc with endpoints $p_{j_1}, p_{j_3}$ that forms a bigon $\tilde{B}$ with $\gamma_{i}^{\pm}$ containing $B$ and intersecting no other ribbon singularities or components of $\Gamma_0$. We abuse notation and let $c\subset D$ denote the corresponding arc in the ribbon disk. In the $-$ case (resp.\ $+$ case), $c$ is an arc in $D_d$ (resp.\ $D_b$); see the dashed arc in \cref{fig:normalized_bypass2} or in \cref{fig:normalized_bypass}. We wish to apply the Legendrian realization principle (\cref{thm:LeRP}) to $c$ in order to witness $\tilde{B}$ as a bypass half-disk, but some care must be taken in the presence of intersecting surfaces. Legendrian realization appeals to a perturbation of the underlying surface which is graphical with respect to its transverse contact vector field. Because there are two different contact vector fields ($X_b, X_d$) describing ribbon-convexity of the immersed surface $D$, such an isotopy may a priori introduce unwanted self-intersections. We need to ensure that there is enough room in the invariant neighborhood described by the relevant contact vector field to perform the isotopy which Legendrian realizes $c$. This is again where normalization is useful. 

The standard "edge-rounding" neighborhood $\overline{\mathcal{U}}_i' \subset \mathcal{U}(\gamma_i)$ of $\gamma_i$ from \cref{lemma:normalized_disk} allows us to choose $c$ to be Legendrian near its endpoints on $\gamma_i$. Therefore, the isotopy which Legendrian realizes $c$ may be assumed constant in $\overline{\mathcal{U}}_i'$. The Legendrian rulings in the "edge-rounding" neighborhood also allow us to take $\hat{D}':= D - \bigcup_{j=1}^n\mathcal{U}_j'$, a slight extension of $\hat{D}$, to have Legendrian boundary near which $c\mid_{\hat{D}'}$ is already Legendrian. 

As $\hat{D}'$ is an extension of $\hat{D}$ into the neighborhood $\mathcal{U}(\vec\gamma)$ from \cref{def:normalized}, neither $X_b$ nor $X_d$ is necessarily everywhere transverse to $\hat{D}'$ (whereas convexity of $\hat{D}$ is witnessed by $X_b = X_d$). Instead, we will describe convexity of $\hat{D}'$ by a different contact vector field, engineered to allow for the Legendrian realization of $c$. If $c\subset D_{\ast}$, where $\ast\in\{b,d\}$, let $Y$ be a contact vector field transverse to $\hat{D}'$ which agrees with $X_{\ast}$ near $\tilde{B}\cap \hat{D}'$. Note that the $Y$-dividing curves may a priori disagree with the original dividing curves (i.e.\ those with respect to $X_b$ and $X_d$) of $\hat{D}'$ elsewhere, but they will agree near $c\mid_{\hat{D}'}$. By multiplying by a bump contact Hamiltonian function supported in a sufficiently small tubular neighborhood of $\hat{D}'$, we decay the support of $Y$. Then $\hat{D}'$ is convex with respect to $Y$, has Legendrian boundary, and $c\mid_{\hat{D}'}$ is a non-isolating arc which is Legendrian near $\partial \hat{D}'$. Applying Legendrian realization with respect to $Y$, constant near $\partial \hat{D}'$, then allows us to Legendrian realize all of $c$ without introducing unwanted intersections with $D$.

With $c$ Legendrian realized, $\tilde{B}$ is a bypass half-disk for $D_b$ (resp. $D_d$) with attaching arc a subarc of $\gamma_{i}^-$ (resp. $\gamma_i^+$). Isotoping $D_b$ (resp. $D_d$) past the bypass half-disk gives the desired modification described in the lemma. Remaining is the more complicated $+$ case, where one of $p_{j_1}, p_{j_3}$ is an endpoint of $\gamma_i^+$. We may assume without loss of generality that the situation is as drawn on the left side of \cref{fig:normalized_bypass}. Isotoping $D_d$ across the bypass half-disk may a priori break normalization, resulting in a surviving intersection point in the interior of $\gamma_i^+$ without any intersection point at the endpoint. From the point of view of $D_d$, this means that $\Lambda$ intersects $D_d$ off of the dividing set after the bypass attachment. If this is the case, we can simply use the characteristic foliation of $D_d$, which is transverse to its dividing set, to induce a Legendrian isotopy of $\Lambda$ supported near $D_d$ which pushes the endpoint of $\gamma_i^-$ onto the "diagonal" dividing set component after the bypass attachment. The result of such an isotopy is indicated in \cref{fig:normalized_bypass}.

In general, given a Legendrian strand intersecting a convex surface transversally near its dividing set, there is enough flexibility --- whether via an isotopy of the Legendrian, or a perturbation of the contact vector field and hence dividing set --- to push the Legendrian onto the dividing set, or past it at most one time. This is a consequence of the local model near a dividing curve described by \cref{lemma:nbd_size}.
\end{proof}

By a \emph{half-bigon} we mean a trigon $\Delta\subset \tilde{D}$ whose piecewise-smooth boundary is
\[
\partial \Delta = \partial_{\Lambda}\Delta \,\cup \,\partial_{\gamma_i^+}\Delta \,\cup\, \partial_{\Gamma_0}\Delta
\]
where $\partial_{\Lambda}\Delta$ is a subarc of $\Lambda$ with one endpoint on $\Gamma_0$ and the other endpoint on some $\gamma_i^+$ ribbon singularity arc, $\partial_{\gamma_i^+}\Delta$ is a subarc of $\gamma_i^+$, and $\partial_{\Gamma_0}\Delta$ is a subarc of $\Gamma_0$.

\begin{figure}[ht]
	\centering
	\begin{overpic}[scale=.4]{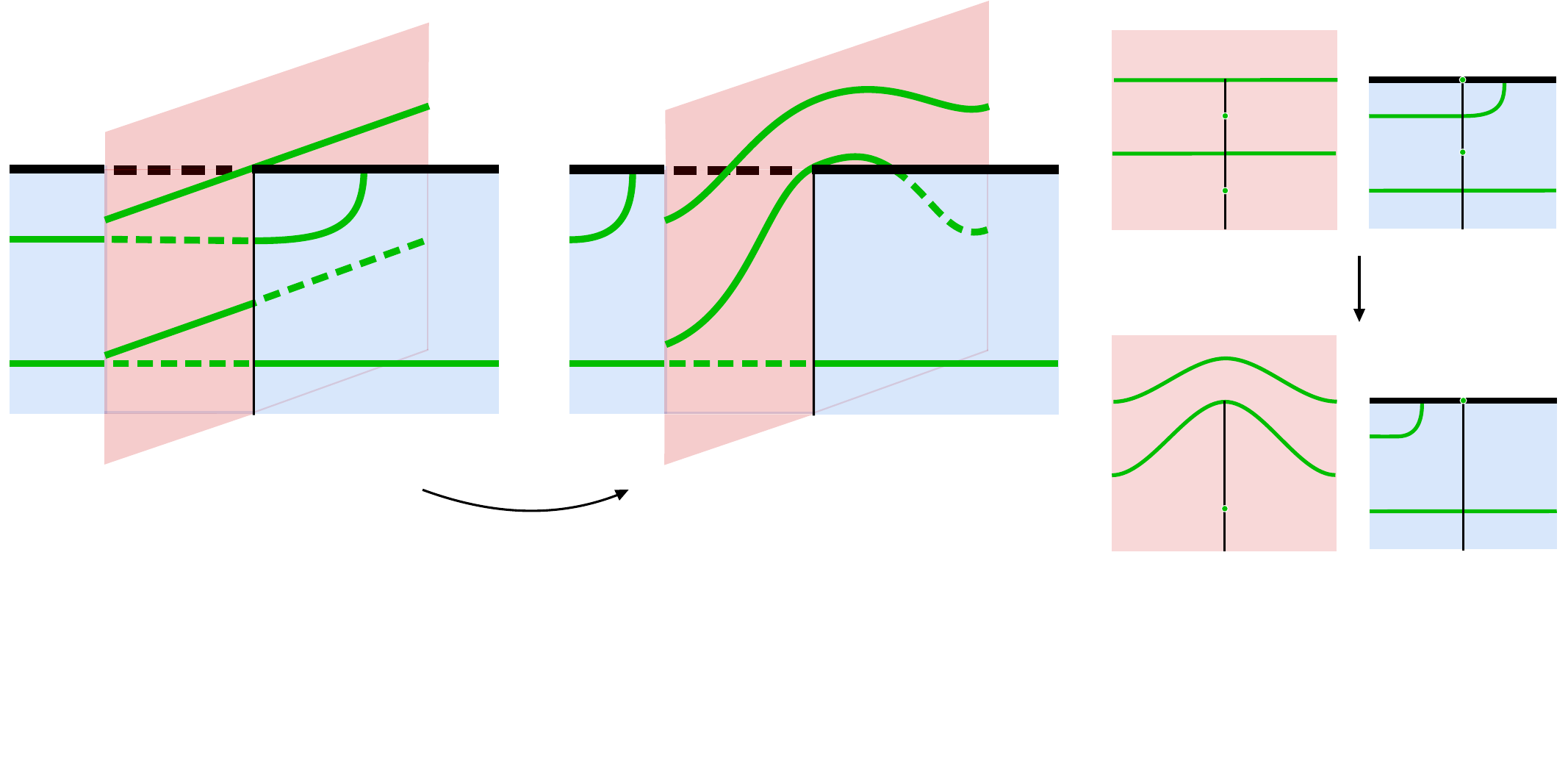}
       \put(2,40){\textcolor{black}{$\Lambda$}}
        \put(38,40){\textcolor{black}{$\Lambda'$}}
	\end{overpic}
    \vskip-1.75cm
	\caption{The statement of the half-bigon removal lemma, \cref{lemma:half-bigon}.}
	\label{fig:half-bigon-reflection}
\end{figure}

\begin{lemma}[Half-bigon removal]\label{lemma:half-bigon}
Let $D$ be a normalized ribbon-convex disk with boundary $\Lambda$ and dividing set $\tilde\Gamma = \Gamma_{0} \cup \{p_1, \dots, p_m\} \subset \tilde{D}$. Suppose $\Delta\subset \tilde{D}$ is a half-bigon such that $\gamma_i^+ - \partial_{\gamma_i^+}\Delta$ contains an intersection with $\Gamma_0$. Then there is an isotopy of $D$ to a normalized ribbon-convex disk $D'$ with Legendrian boundary $\Lambda'$, Legendrian isotopic to $\Lambda$, which modifies the dividing set $\Gamma_{\tilde{D}}$ as indicated in \cref{fig:half-bigon-reflection}.
\end{lemma}

\begin{figure}[ht]
	\centering
	\begin{overpic}[scale=.4]{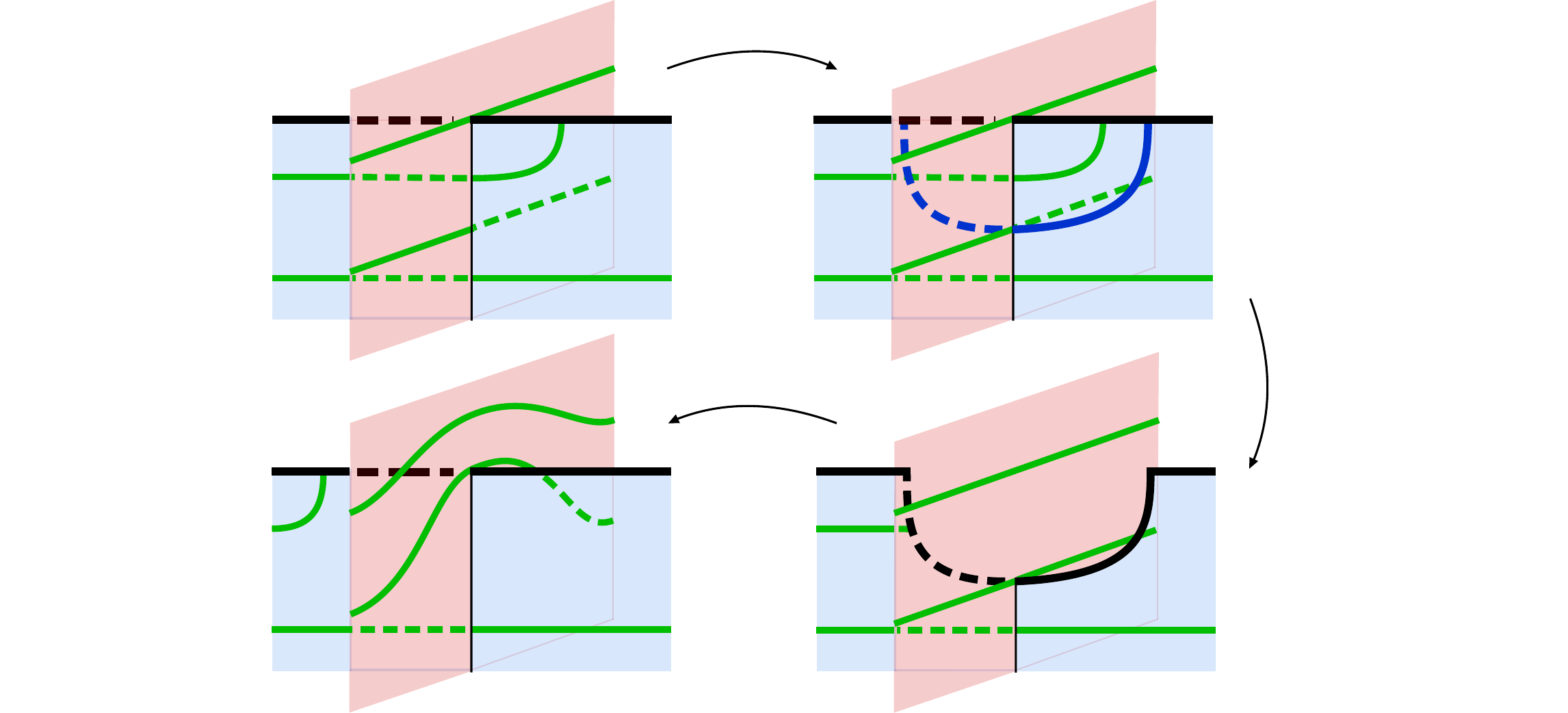}
    
	\end{overpic}
    	\caption{The proof of the half-bigon removal lemma.}
	\label{fig:half-bigon-reflection-proof}
\end{figure}

\begin{proof}
We may assume that our local model is as drawn in the top left of \cref{fig:half-bigon-reflection-proof}. Let $a\subset b_i$ be the blue arc drawn in the top right of the figure. Near $\gamma_i$, we can take $a$ to be Legendrian by the standard model near a normalized ribbon singularity. Because $a$ is non-isolating with respect to $\Gamma$ on either side of the (already Legendrian) arc $\gamma_i$, we may apply the Legendrian realization principle to $a$ and take it to be Legendrian. We obtain $D'$ by excising the part of $b_i$ cobounded by $a$ and $\Lambda$ and smoothing the corners. This is depicted in the bottom row of the figure. Finally, to observe that $\Lambda':= \partial D'$ is Legendrian isotopic to $\Lambda$, note that the part of $b_i$ that was carved out to produce $D'$ is a convex disk with a single dividing arc, with one endpoint on $\Lambda$ and one endpoint on $\Lambda'$. The standard characteristic foliation on such a disk guides a Legendrian isotopy from $\Lambda'$ to $\Lambda$.
\end{proof}

It will be useful to undo the bigon and half-bigon removal moves. This is possible, with some care regarding the relative orientation of the intersection between the two sheets.

\begin{figure}[ht]
	\centering
    \vskip-0.75cm
	\begin{overpic}[scale=.4]{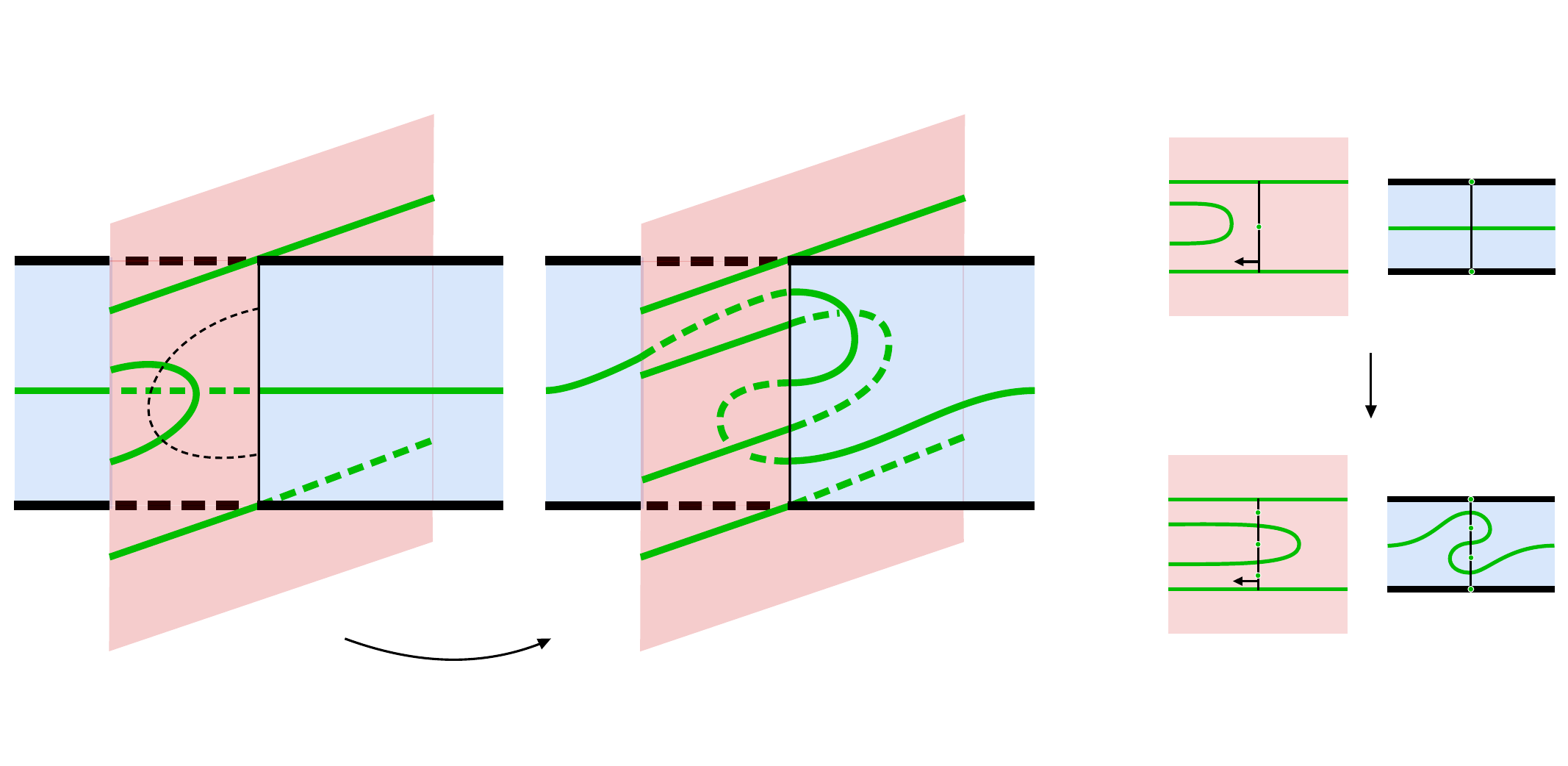}
       
	\end{overpic}
    \vskip-0.75cm
	\caption{The statement of \cref{lemma:bigon_creation}. Note that the direction of the wiggled dividing set on the blue sheet is fixed; the arrows on the right side of the figure indicate the outward co-orientation of the blue sheet. The figure also only describes a finger move through $\gamma_i^-$, but finger moves through $\gamma_i^+$ work the same way.}
	\label{fig:finger}
\end{figure}

\begin{lemma}[Bigon creation]\label{lemma:bigon_creation}
Let $D$ be a normalized ribbon-convex disk with boundary $\Lambda$ and dividing set $\Gamma_{0} \cup \{p_1, \dots, p_m\} \subset \tilde{D}$. Let $\gamma_0\subset D$ be an arc with endpoints on some $\gamma_i$ defining a bigon $B \subset \tilde{D}$ such that $\partial B = \gamma_0 \cup \partial_{\gamma_i^{\pm}} B$ and $\Gamma_0 \mid_B$ is a single arc with endpoints on $\gamma_0$ and a single point on $\gamma_{i}^{\pm}$, as depicted on the left side of \cref{fig:finger}. Then there is an isotopy of $D$ relative to $\Lambda$ which performs a finger move on the arc in $\Gamma_0\mid_B$ and creates a bigon with as indicated in the same figure. 
\end{lemma}

\begin{proof}
We begin with a sanity check. Note that the newly produced dividing set in \cref{fig:finger} has three bigons to which \cref{lemma:bigon_removal} may be applied: two in $b_i$ and one in $d_i$. The reader may verify that removing any of these three bigons results in the original dividing set configuration on the left. Moreover, this clarifies that the direction of the "wiggling" induced by the finger move is important, as "wiggling" the other way and then removing the created bigon would result in a homotopically trivial dividing curve component in $d_i$. 

The sanity check is useful, but not a proof. Begin by Legendrian realizing  $\gamma_0$, which is possible to do since the graph $\gamma_0 \cup \gamma_i$ is non-isolating in $D_d$. Let $N(B) \cong [-\ve, \ve]\times B$ be a small standard neighborhood of the convex bigon $B$ bounded by the Legendrians $\gamma_0$ and $\gamma_i$. Let $\partial_- N(B):= \partial N(B) \cap b_i$. Replacing $b_i$ with 
\[
b_i - (\partial_- N(B)) \cup (\partial N(B) - \partial_- N(B))
\]
gives a new piecewise smooth convex surface. Performing the usual corner rounding to smooth the perturbed band produces the desired dividing set configuration.
\end{proof}

\begin{figure}[ht]
	\centering
	\begin{overpic}[scale=.4]{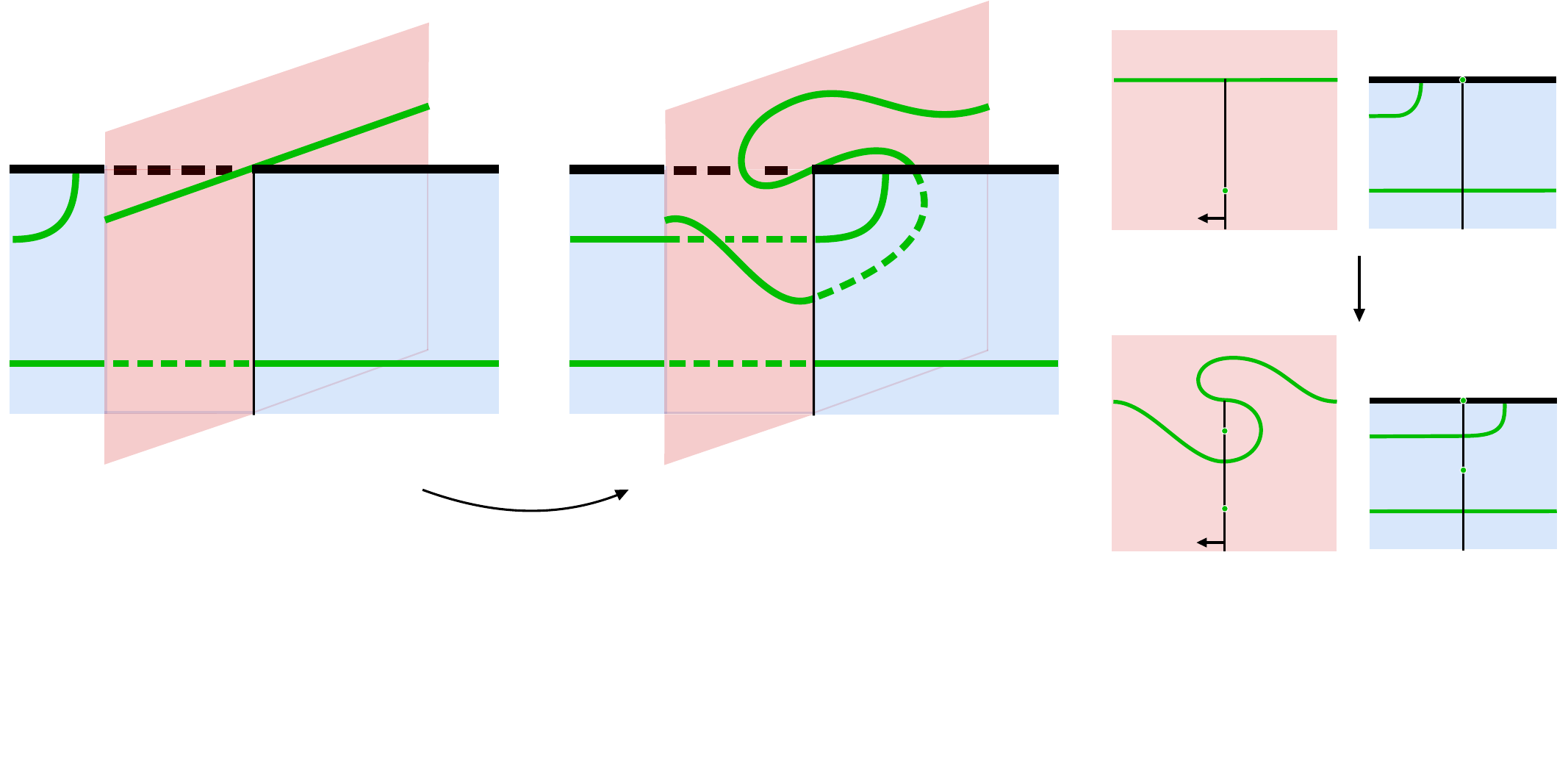}
       
	\end{overpic}
    \vskip-1.75cm
	\caption{The statement of the half-bigon creation lemma, \cref{lemma:half-bigon-creation}.}
	\label{fig:half-bigon-creation}
\end{figure}

\begin{lemma}[Half-bigon creation]\label{lemma:half-bigon-creation}
Let $D$ be a normalized ribbon-convex disk with boundary $\Lambda$ and dividing set $\Gamma_{0} \cup \{p_1, \dots, p_m\} \subset \tilde{D}$. Let $\Lambda_0 \subset \Lambda$ be a subarc of $\Lambda$ with one endpoint on $\Gamma_0$, another endpoint on $\gamma_i$, and no other such intersections in its interior. Then, in a neighborhood of $\Lambda$, there is an isotopy of $D$ to a normalized ribbon-convex disk $D'$ with Legendrian boundary $\Lambda'$, Legendrian isotopic to $\Lambda$, which modifies the dividing set $\Gamma_{\tilde{D}}$ as indicated in \cref{fig:half-bigon-creation}.
\end{lemma}

\begin{proof}
This follows from a combination of \cref{lemma:bigon_creation} and \cref{lemma:half-bigon}. Specifically, perform a finger move on $\Gamma_0$ through $\gamma_i$ using \cref{lemma:bigon_creation} along a path in $b_i$ parallel to $\Lambda_0$. Then remove the created half-bigon by \cref{lemma:half-bigon}. The dividing set is then modified as indicated in the figure. 
\end{proof}

\begin{figure}[ht]
	\centering
	\begin{overpic}[scale=.8]{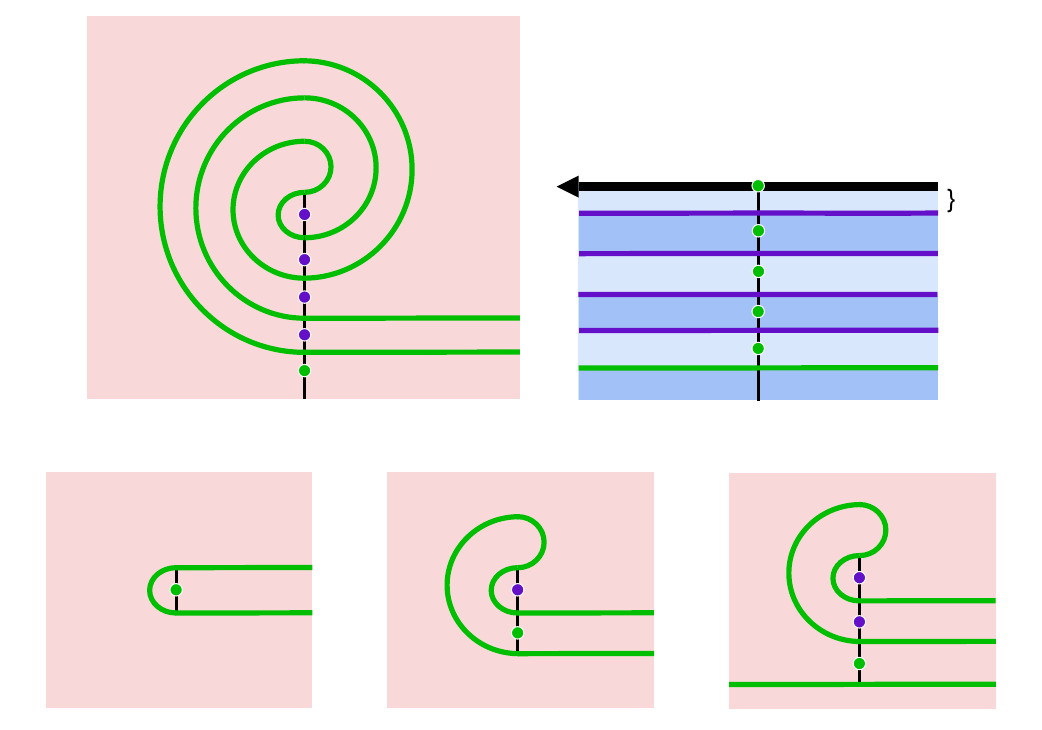}
       \put(10,3){\small not a lollipop}
       \put(11.5,0){\small i.e.\ index $0$}
       \put(46,0){\small index $1$}
       \put(79,0){\small index $2$}
       \put(92.25,50.25){\small positive}
	\end{overpic}
    \caption{A positive lollipop of index $4$ (top row). Given the orientation of $D$, the darker regions on the right represent the negative region $R_-$. Several additional examples of lollipops (bottom row).}
	\label{fig:lollipop}
\end{figure}

\begin{definition}
Fix $k\geq 1$. A \emph{lollipop of index $k$} is a neighborhood $U\subset \tilde{D}$ of an endpoint of $\gamma_i^-$ with dividing curves isotopic (up to rotation and reflection) to those on the left side of the top row of \cref{fig:lollipop}, where $k$ is the number of intersections of the component of $\Gamma_0$ with $\gamma_i^+$ above the lowermost dividing arc in the band on the right; for clarity, the arcs which contribute to the index are colored in purple, and the referenced lowermost dividing arc is green. Given an orientation of $D$, we call the lollipop \emph{positive} (resp.\ \emph{negative}) if the region adjacent to $\Lambda$ at the endpoint of $\gamma_i^+$ is positive (resp.\ negative).
\end{definition}

The next lemma shows that, by negatively stabilizing the boundary Legendrian, we may remove all lollipops. 

\begin{lemma}\label{lemma:lollipop-reduce}
Let $D$ be an oriented normalized ribbon-convex disk with Legendrian boundary $\Lambda$. There is an isotopy of $D$ to a normalized ribbon-convex disk $D_-$ such that    
\begin{enumerate}
    \item $\Lambda_- := \partial D_-$ is a negative multistabilization of $\Lambda$, and 
    \item $D_-$ has no lollipops.
\end{enumerate}
In particular, $\Lambda$ and $\Lambda_-$ have a common transverse pushoff.
\end{lemma}

\begin{figure}[ht]
	\centering
	\begin{overpic}[width=\textwidth]{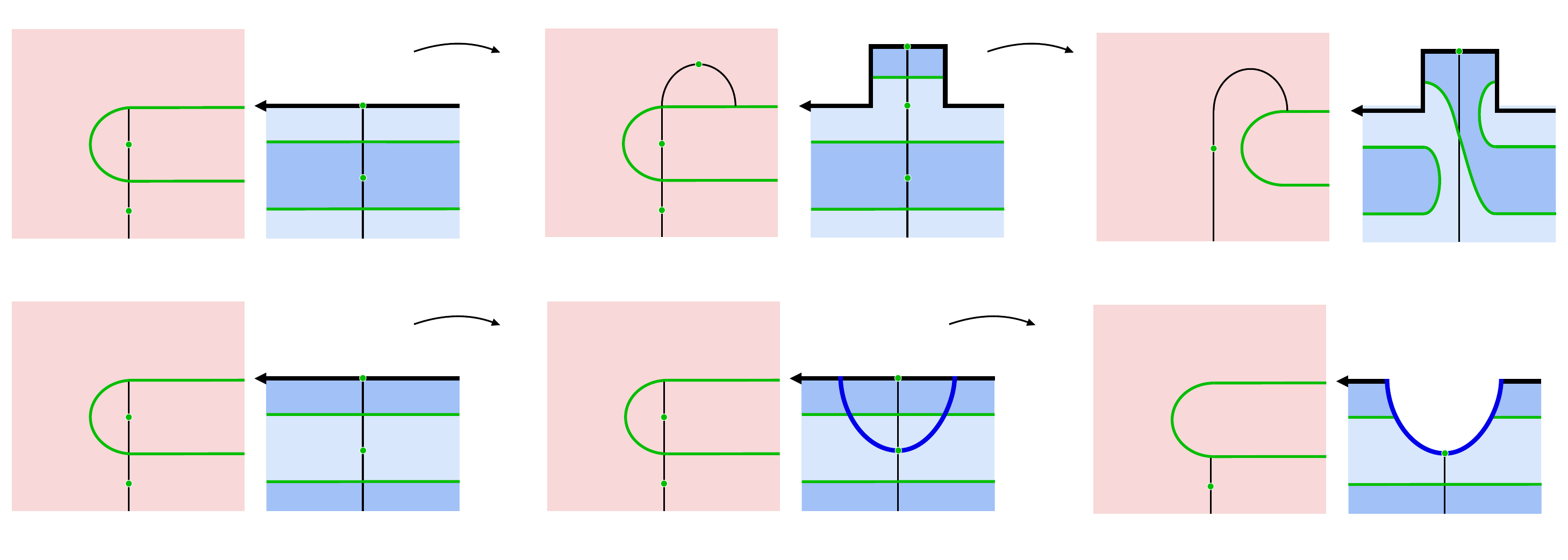}
       
	\end{overpic}
	\caption{The proof of \cref{lemma:lollipop-reduce}.}
	\label{fig:lollipop-reduce}
\end{figure}

\begin{proof}
Suppose $D$ has a lollipop of index $k>0$. We argue that by negatively stabilizing $\Lambda$, we can reduce the index of the lollipop by $1$. We consider separately the cases the lollipop is positive or negative.

First, suppose the lollipop is positive. Since the lollipop is index $k>1$, we may localize further to a "core" lollipop of index $1$. As in the proof of \cref{lemma:half-bigon}, we can Legendrian realize an extension of of $\gamma_i^-$ (top row of \cref{fig:lollipop-reduce}), and append a standard neighborhood to thicken $b_i$, noting that positivity of the lollipop means the Legendrian stabilization involved will be negative by \cref{lemma:destab}. After this extension, the core of the former lollipop is a bigon to which \cref{lemma:bigon_removal} can be applied. The index of the resulting lollipop is $\leq k-1$.

Second, suppose the lollipop is negative. Again, we localize to a "core" lollipop of index $1$. Here we simply Legendrian realize the blue curve in the bottom row of \cref{fig:lollipop-reduce} (c.f.\ the proof of \cref{lemma:half-bigon}) and excise the boundary parallel region it bounds. The index of the resulting lollipop is $\leq k-1$, and by \cref{lemma:destab} the stabilization involved is negative.
\end{proof}

We present one last modification which exchanges non-removable bigons for non-removable half-bigons, again at the cost of negatively stabilizing the boundary. 

\begin{figure}[ht]
	\centering
	\begin{overpic}[scale=.35]{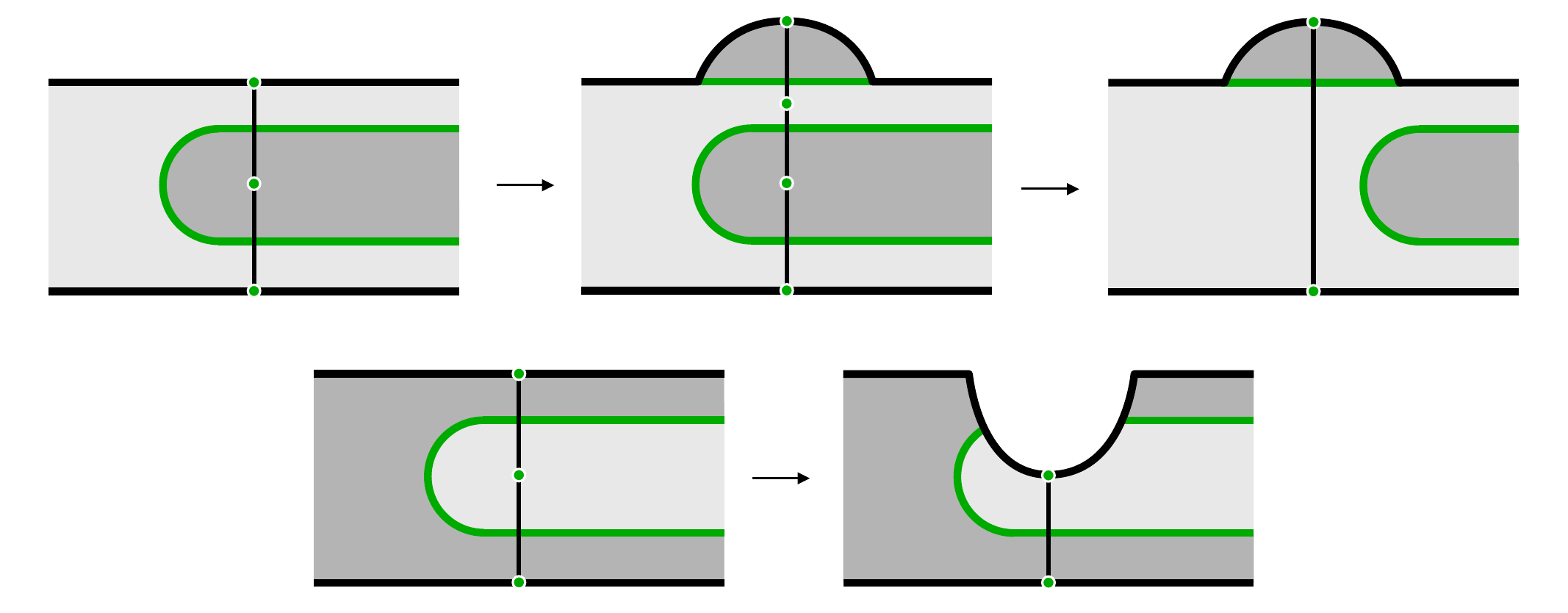}
       
	\end{overpic}
	\caption{Exchanging a non-removable bigon for a half-bigon as in \cref{lemma:non-removable-exchange}.}
	\label{fig:non-removable-exchange}
\end{figure}

\begin{lemma}[Non-removable bigon exchange]\label{lemma:non-removable-exchange}
Let $D$ be an oriented normalized ribbon-convex disk with Legendrian boundary $\Lambda$. Suppose $B\subset \tilde{D}$ is a non-removable bigon; i.e.\ a bigon such that $\partial B = \partial_{\gamma_i^+}B \cup \partial_{\Gamma_0}B$, where $\Gamma_0 \cap \gamma_i^+ = \partial(\partial_{\Gamma_0}B)$ Then there is an isotopy of $D$ to a normalized ribbon-convex disk $D_-$ such that 
\begin{enumerate}
    \item $\Lambda_- := \partial D_-$ is a negative stabilization of $\Lambda$, and
    \item $B$ is replaced by a half-bigon as in \cref{fig:non-removable-exchange}, according to whether $B$ is part of the positive region (bottom row) or the negative region (top row). 
\end{enumerate}
In particular, $\Lambda$ and $\Lambda_-$ have a common transverse pushoff.
\end{lemma}

\begin{proof}
If $B$ is part of the negative region of $\tilde{D}$, then as in the proof of \cref{lemma:lollipop-reduce} (see \cref{fig:lollipop-reduce}) we begin by stabilizing the boundary Legendrian by Legendrian realizing an extension of $\gamma_i$. By \cref{lemma:destab}, we have performed a negative stabilization of the boundary Legendrian. The original bigon is now removable, and we apply \cref{lemma:bigon_removal}. This is the top row of \cref{fig:non-removable-exchange}.

On the other hand, if $B$ is part of the positive region --- see the bottom row of \cref{fig:non-removable-exchange} --- we perform a negative stabilization of the boundary by instead "eating in" into the band. As in \cref{fig:lollipop-reduce}, we Legendrian realize a curve in the band region, which, upon excision of the created boundary parallel region, leaves the lower right configuration in \cref{fig:non-removable-exchange}. By \cref{lemma:destab}, this negatively stabilizes the boundary Legendrian and what remains from $B$ is a half-bigon. 
\end{proof}

\section{Transverse derivative links}\label{sec:BZ1review}

Here we briefly review the relevant result from \cite{BZ0forthcoming}, which is the implication (3) $\Rightarrow$ (1) of Theorem 1.12 in that article. To make this work as self-contained as possible, we sketch the proof, deferring to the above reference for a thorough treatment. For convenience, we restate the result: 

\bz*

\begin{proof}[Proof sketch.]
Let $\Sigma$ be a genus-$g$ Seifert surface and $T:= T_1 \cup \cdots \cup T_g$ the assumed max-sl unlink derivative. Because $T$ is a max-sl transverse unlink, inadmissible transverse $0$ surgery along $T$ yields $S^3_T(0) = (\#^g (S^1 \times S^2), \, \xi_{\st})$. Moreover, the image of $\mathcal{T}$ in $S^3_T(0)$ is a topological unknot $U$. 

We first claim that $\mathrm{sl}(U) = -1$, which can be verified by a characteristic foliation computation. Assume the characteristic foliation of $\Sigma$ is Morse-Smale, so that $-1 = \mathrm{sl}(\mathcal{T}) = -(e_+ - h_+) + (e_- - h_-)$, where $e_{\pm}, h_{\pm}$ are the counts of positive (resp.\ negative) elliptic and hyperbolic singular points, respectively. Each component of the derivative $T$ contributes to two boundary components of the cut-open Seifert surface $\Sigma':=\Sigma - \nu(T)$ in $S^3_T(0)$, and when we cap off these two boundary components in $S^3_T(0)$ with standard meridional disks in the transverse surgery torus, the orientation forces one capping disk to contribute a positive elliptic point and the other to contribute a negative elliptic point. Thus, 
\[
\mathrm{sl}(U) = -(e_+ + g - h_+) + (e_- + g - h_-) = -(e_+ - h_+) + (e_- - h_-) = -1.
\]  
Since $U$ is a max-sl unknot in $(\#^g (S^1\times S^2), \xi_{\mathrm{st}})$, performing a further inadmissible transverse $0$ surgery along $U$ produces $(\#^{g+1} (S^1\times S^2), \xi_{\mathrm{st}})$. 

Summarizing so far, we have shown that inadmissible transverse $0$ surgery along the link $\mathcal{T} \cup T\subset (S^3, \xi_{\st})$ yields $(\#^{g+1} (S^1\times S^2), \xi_{\mathrm{st}})$. This contact manifold is Stein filled by $(\natural^{g+1} (S^1 \times B^3), \omega_{\mathrm{st}})$. Using \cite{gay2002symplectichandles,wendl2013symplectichandles}, we may attach symplectic $2$-handles with symplectic co-cores to $\natural^g(S^1 \times B^3)$ along the dual link of $\mathcal{T} \cup T$ to induce a weak symplectic filling of the boundary. Smoothly, the boundary is $S^3$, and since it weakly filled, must be $(S^3, \xi_{\mathrm{st}})$. Moreover, $\mathcal{T} \cup T$ bounds the symplectic co-core disks of the $2$-handles. Each component (in particular, $\mathcal{T}$) is therefore symplectically slice in the weak filling, which by Gromov \cite{gromov1985pseudo} must be $(B^4, \omega_{\mathrm{st}})$. 
\end{proof}

\section{Star-shaped ribbon disks}\label{sec:ribbon_disk_analysis}

In this section, we consider \emph{star-shaped} ribbon disks:  A ribbon disk $D$ is called \emph{star-shaped} if for each singular arc $\gamma_i^+$, we have that $D \setminus \gamma_i^+$ has a component $D_i$ containing exactly one singular arc of the form $\gamma_j^-$.  We call $D_i$ a \emph{subdisk}, and for ease of notation, we reindex the arcs of the form $\gamma_j^-$ so that $\gamma_i^- \subset D_i$.  Note that this reindexing destroys the previous convention, in which $\gamma_i^+$ and $\gamma_i^-$ were identified.  It follows that the $D \setminus \{\gamma_1^+,\dots,\gamma_g^+\}$ contains $g+1$ components, the subdisks $D_1,\dots,D_g$ and one additional region, $S$, containing no arcs of the form $\gamma_j^-$.  We call $S$ the \emph{star-shaped region}.  A depiction of a star-shaped disk appears in \cref{fig:starshaped1}.  The term "star-shaped" is inspired by the fact that the corresponding ribbon code is a star graph in which each edge contains a single marking.  See~\cite{friedl-misev-zupan} for further details about ribbon codes.

Next, given an oriented star-shaped ribbon disk which is normalized ribbon-convex, we consider the possible intersections of the dividing set $\Gamma$ (in this section, we will not distinguish between $\Gamma$ or its lift to $\Gamma_{\tilde{D}}$) with the components $D_1,\dots,D_g,S$, where $\Gamma$ meets the singular set of $D_i$ minimally.  First, define $I(D_i) := |\Gamma \cap \gamma_i^-| - |\Gamma \cap \gamma_i^+| - 1$.  We will use the following lemma.

\begin{lemma}\label{lem:I}
A normalized ribbon-convex ribbon disk satisfies $\sum I(D_i) = 0$.
\end{lemma}

\begin{proof}
Suppose $\gamma_i^+$ is identified with $\gamma_j^-$.  By the interleaving condition (see \cref{lemma:normalized_disk}), we have $|\Gamma \cap \gamma_i^+|-1 = |\Gamma \cap \gamma_j^-|$, or equivalently, $|\Gamma \cap \gamma_j^-| - |\Gamma \cap \gamma_i^+| - 1 = 0$.  Taking the sum over all $i$ yields the statement of the lemma.
\end{proof}

Now, suppose $D$ is a normalized ribbon-convex disk with boundary $\Lambda$, chosen to minimize $(|\Gamma \cap \gamma|, -\tb(\Lambda))$ with the dictionary ordering among all possible choices up to Legendrian isotopy and negative stabilizations and destabilizations of $\Lambda$. Consider an arc component of $\Gamma \cap D_i$ or $\Gamma \cap S$.  We call the arc a \emph{$\Lambda$-arc} if both boundary points are contained in $\Lambda := \partial D$ and a \emph{$\gamma_i^+$-arc} if both boundary points of are contained in $\gamma_i^+$.  Otherwise, the arc is a \emph{seam}.  We operate under the following principles:

\begin{itemize}

\item Since the ambient contact structure is tight, and each subdisk $D_i$ is an embedded convex disk, Giroux's criterion implies that $\Gamma \cap D_i$ contains no closed components. (However, we cannot a priori obstruct closed components of $\Gamma$.)

\item We may assume that every component of $\Gamma \cap D_i$ meets $\gamma_i^-$ in at most one point, with the possible exception of an arc that meets $\gamma_i^-$ twice along $\partial \gamma_i^-$; otherwise $\Gamma$ would cut out a bigon with $\gamma_i^-$, and this bigon is either an interior bigon or part of a lollipop, both of which can be removed using \cref{lemma:bigon_removal} or \cref{lemma:lollipop-reduce}. 

\item By \cref{lemma:bigon_removal} and \cref{lemma:non-removable-exchange}, no $\Lambda$-arc of $\Gamma \cap D_i$ cobounds a bigon with an arc in $\gamma_i^+$ unless the bigon contains $\gamma_i^-$.

\item No $\Lambda$-arc of $\Gamma \cap D_i$ cobounds a bigon unless the bigon contains $\gamma_i^-$, and no $\Lambda$-arc of $\Gamma \cap S$ cobounds a bigon; otherwise, we can perform a negative destabilization to increase $\tb(\Lambda)$.

\item If a seam avoids $\gamma_i^-$, then $|\Gamma \cap \gamma_i^+| = 1$; otherwise, we could apply half-bigon removal (\cref{lemma:half-bigon}). 

\item Finally, by exchanging bigons with half-bigons if necessary, we may assume $S \setminus \Gamma$ contains no bigons.
\end{itemize}

We let $\rot(D)$ denote the signed sum of the Euler characteristics of $D \setminus \Gamma$, so that $\tb(\Lambda) - \rot(D) = -1$.  We can also define these quantities locally:  First, for a subdisk $D_i$, define $\tb(D_i) := -| \Gamma \cap D_i \cap \Lambda|/2$, and for $S$, define $\tb(S) := -| \Gamma \cap S \cap \Lambda|/2$.  It follows immediately that $\tb(\Lambda) = \tb(S) + \sum \tb(D_i)$.  Next, let $R$ be a region of $D_i \setminus \Gamma$.  Since $\Gamma$ meets $D_i$ in arcs, it follows that $R$ must be a disk.  Define $\rot(R) := \pm (1 - |R \cap \gamma_i^+|/2)$, where the sign $\pm$ is determined by the sign of the region $R$, and $|R \cap \gamma_i^+|$ is the number of arc components of $\partial R$ contained in $\gamma_i^+$.  Similarly, for a region $R$ of $S \setminus \Gamma$, define $\rot(R) := \pm (1 - |R \cap (\bigcup \gamma_i^+)|/2)$.  Finally, let $\rot(D_i)$ denote the sum of $\rot(R)$ over all regions $R$ of $D_i \setminus \Gamma$, and let $\rot(S)$ denote the same for $S \setminus \Gamma$.  We have

\begin{lemma}\label{lem:chi}
A normalized ribbon-convex star-shaped ribbon disk following the above principles satisfies
\[ \rot(S) + \sum \rot(D_i) = \rot(D).\]
\end{lemma}

\begin{proof}
Let $\chi^+(D)$ denote the sum of the Euler characteristics of the positive regions, and let $\chi^-(D)$ denote the sum of the Euler characteristics of the negative regions, so that $\rot(D) = \chi^+(D) - \chi^-(D)$.  Define $\chi^{\pm}(D_i)$ and $\chi^{\pm}(S)$ similarly, so that $\rot(S) = \chi^+(S) - \chi^-(S)$ and $\rot(D_i) = \chi^+(D_i) - \chi^-(D_i)$.  Consider a positive region $R$ of $D \setminus \Gamma$, where $R$ is the union of regions $R_1,\dots,R_n$, each contained in some $D_i \setminus \Gamma$ or $S \setminus \Gamma$.  As noted above, each region $R_i$ is a disk.  If we add up the Euler characteristics of $R_i$ in order to obtain $\chi(R)$, we have counted edge of $R_i \cap R_j$ twice (instead of once), and thus we have
\[ \chi^+(R) = n - \left|R \cap \left(\bigcup \gamma_i^+\right)\right|.\]
On the other hand,
\[ \sum \chi^+(R_i) = \sum \left(1 - \left|R_i \cap \left(\bigcup \gamma_i^+\right)\right|/2\right) = n - \left|R \cap \left(\bigcup \gamma_i^+\right)\right|,\]
since each intersection of $R$ with $\bigcup \gamma_i^+$ appears in two regions $R_i$ and $R_j$.  A similar calculation shows that if $R$ is negative, then $\chi^-(R) = \sum \chi^-(R_i)$, and the statement of the lemma follows.
\end{proof}

Next, we consider the effect of parallel arcs of $\Gamma \cap D_i$ on the computation of $\rot(D_i)$.  Two arcs of $\Gamma \cap D_i$ or $\Gamma \cap S$ are \emph{parallel} if they cobound a rectangular region with arcs contained in $\partial D_i$ that do not contain the endpoints of $\gamma_i^+$, and a collection of arcs is called \emph{parallel} if any pair is parallel.  The cardinality of a collection of parallel arcs is called the \emph{weight} of the class of arcs.  Note that all arcs in a collection of parallel arcs must be of the same type ($\Lambda$-arcs, $\gamma_i^+$-arcs, or seams).

\begin{lemma}\label{lem:parallel}
Adding two parallel arcs to $\Gamma \cap D_i$ does not change $\rot(D_i)$, and adding two parallel arcs to $\Gamma \cap S$ does not change $\rot(S)$.
\end{lemma}

\begin{figure}[ht]
	\centering
	\begin{overpic}[scale=0.3]{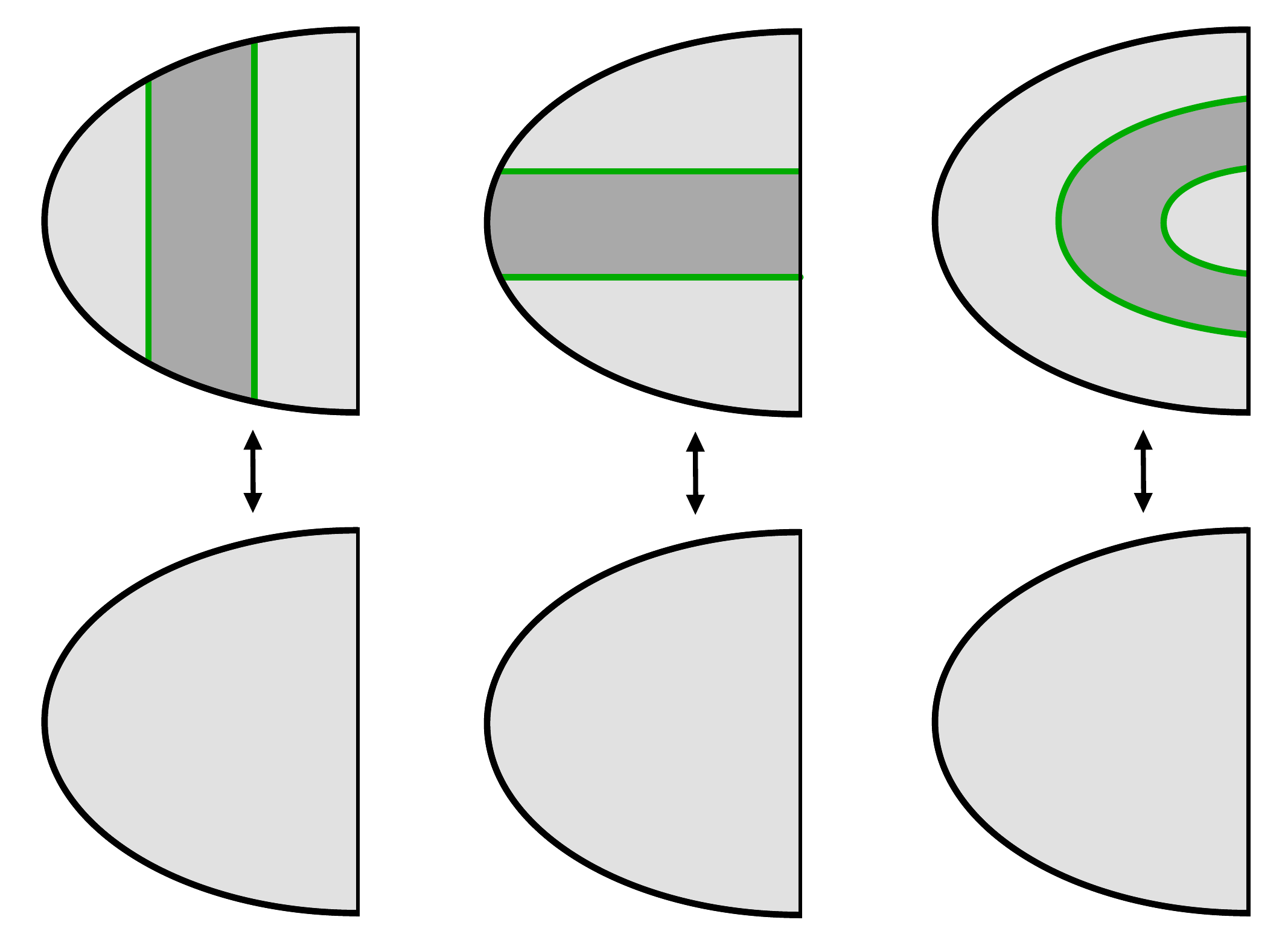}
     \put(51,20){\small $R^*$}
     \put(43,13){\small $+1 - \tfrac{k' + k'' -1}{2}$}
     \put(15,20){\small $R^*$}
     \put(11,13){\small $+1 - \tfrac{k}{2}$}
     \put(86,20){\small $R^*$}
     \put(78,13){\small $+1 - \tfrac{k' + k'' -2}{2}$}

     \put(6.75,60.5){\small $R'$}
     \put(6.75,55){\tiny $+1$}
     \put(14,65){\small $R$}
     \put(14,55){\tiny $-1$}
     \put(21.5,67){\small $R''$}
     \put(20.25,55){\tiny $+1 - \tfrac{k}{2}$}

     \put(57.5,64){\small $R'$}
     \put(47,64){\tiny $+1 - \tfrac{k'}{2}$}
     \put(58,55){\small $R$}
     \put(49,55){\tiny $-\tfrac{1}{2}$}
     \put(57.5,46){\small $R''$}
     \put(47,46){\tiny $+1 - \tfrac{k''}{2}$}

     \put(92.5,67.5){\small $R'$}
     \put(98,67.5){\tiny $+1 - \tfrac{k'}{2}$}
     \put(93,61.5){\small $R$}
     \put(98,61.5){\tiny $-0$}
     \put(92.5,55){\small $R''$}
     \put(98,55){\tiny $+1 - \tfrac{k''}{2}$}
	\end{overpic}
	\caption{Adding two parallel arcs does not change $\rot(D_i)$ or $\rot(S)$, as shown in \cref{lem:parallel}.}
	\label{fig:parallel}
\end{figure}

\begin{proof}
There are three cases to consider, depending on whether the two new arcs are $\Lambda$-arcs, $\gamma_i^+$-arcs, or seams; see \cref{fig:parallel}.  Suppose the two arcs cobound a rectangle $R$, which is adjacent to regions $R'$ and $R''$ on either side.  We suppose without loss of generality that $R$ is negative and $R'$ and $R''$ are positive.  If both arcs are $\Lambda$-arcs, then $\chi^-(R) = 1$, and only one of the regions $R'$ or $R''$ may meet $\gamma_i^+$, so suppose that $R''$ meets $\gamma_i^+$ in $k$ edges.  Then $\chi^+(R') = 1$ and $\chi^+(R'') = 1 - \tfrac{k}{2}$.  On the other hand, removing the two parallel arcs combines $R$, $R'$, and $R''$ into a single positive region $R^*$ which meets $\gamma_i^+$ in $k$ edges, so that $\chi^+(R^*) = 1 - \tfrac{k}{2}2$.  We verify that $\chi^+(R') + \chi^+(R'') - \chi^-(R) = \chi^+(R^*) = 1 - \tfrac{k}{2}$.

Similarly, if both arcs are seams, then $R$ meets $\gamma_i^+$ in a single edge, so that $\chi^-(R) = \tfrac{1}{2}$.  We suppose that $R'$ meets $\gamma_i^+$ in $k'$ arcs and $R'$ meets $\gamma_i^+$ in $k''$ arcs; hence $\chi^+(R') = 1 - \tfrac{k'}{2}$ and $\chi^+(R'') = 1 - \tfrac{k''}{2}$.  Removing the two parallel arcs merges $R$, $R'$, and $R''$ into a single positive region $R^*$ which meets $\gamma_i^+$ in $k' + k'' - 1$ edges (since a pair of edges merges into one).  We verify $\chi^+(R') + \chi^+(R'') - \chi^-(R) = \chi^+(R^*) = \tfrac{3}{2} - \tfrac{k' + k''}{2}$. 

Third, if both arcs are $\gamma_i^+$-arcs, $R$ meets $\gamma_i^+$ in two edges and $\chi^-(R) = 0$, while $R'$ meets $\gamma_i^+$ in $k'$ edges and $R''$ meets $\gamma_i^+$ in $k''$ edges as above.  Removing the two parallel arcs combines $R$, $R'$, and $R''$ into a single positive region $R^*$ which meets $\gamma_i^+$ in $k' + k'' - 2$ edges (here two pairs are merged to two edges), and we have $\chi^+(R') + \chi^+(R'') - \chi^-(R) = \chi^+(R^*) = 2 - \tfrac{k'+k''}{2}$.  Finally, the computations for the star-shaped region $S$ follow similarly.
\end{proof}

\begin{lemma}\label{lem:star}
The star-shaped region $S$ satisfies $\tb(S) - \rot(S) \leq \tfrac{g}{2}-1.$
\end{lemma}

\begin{proof}
We induct on the number of arcs of $\Gamma \cap S$ that intersect $\Lambda$.  For the base case, suppose that no arc of $\Gamma \cap S$ meets $\Lambda$.  Then $S \setminus \Gamma$ consists of some number of polygonal disks, where each has an even number of sides.  Since $\Gamma \cap S$ contains no bigons, each polygonal disk has at least four sides, and by \cref{lem:parallel}, we can add pairs of parallel arcs to $\Gamma$ until every non-rectangle has six sides without altering $\rot(S)$.  An example is shown in \cref{fig:star1}.  After this addition, $S \setminus \Gamma$ contains some number of rectangles and some number of hexagons, where each rectangle $R$ satisfies $\rot(R) = 0$ and each hexagon $H$ satisfies $\rot(H) = \pm \tfrac{1}{2}$, and such that there are exactly $g-2$ hexagons.  Since no arcs of $\Gamma \cap S$ meet $\Lambda$, we have $\tb(S) = 0$; hence,
\[\tb(S) - \rot(S) \leq (g-2) \cdot \tfrac{1}{2}  = \tfrac{g}{2} - 1.\]

\begin{figure}[ht]
	\centering
	\begin{overpic}[scale=.35]{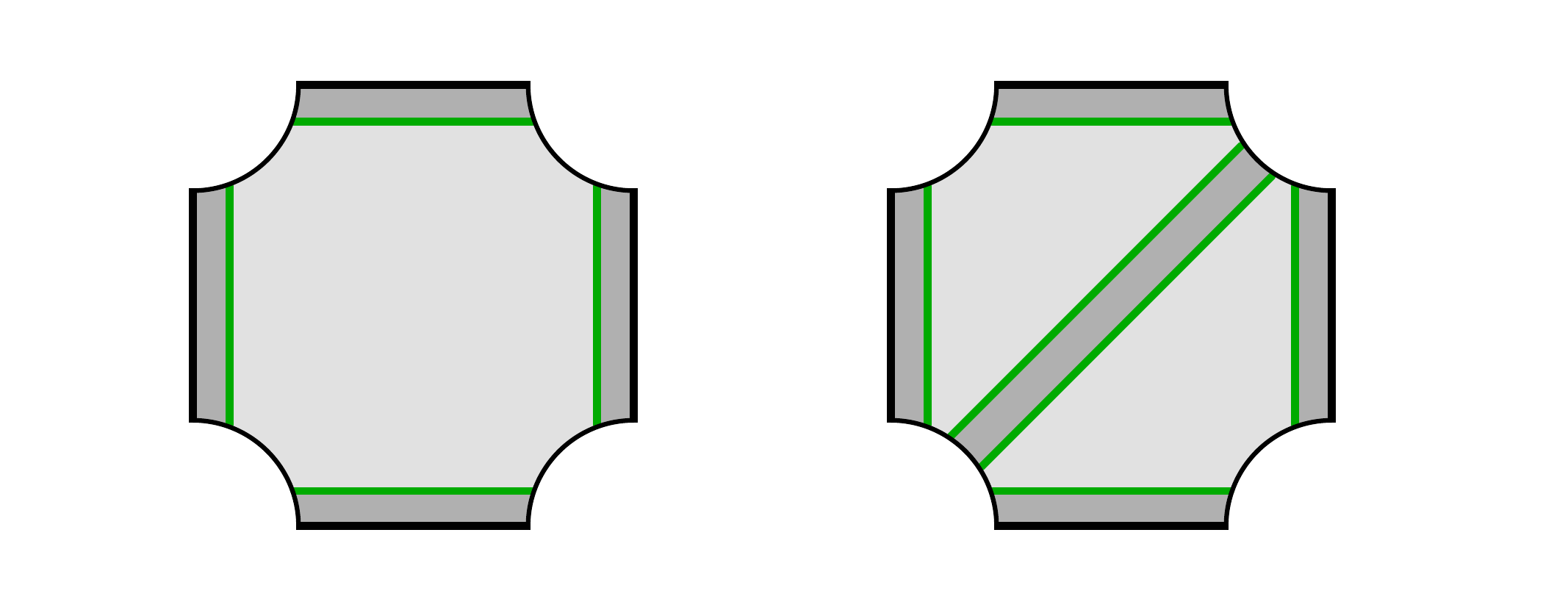}
        \put(24,19){\small $-1$}
        \put(13,19){\tiny $0$}
        \put(38.75,19){\tiny $0$}
        \put(26,6){\tiny $0$}
        \put(26,31.85){\tiny $0$}

        \put(70.5,19){\tiny $0$}
        \put(57.5,19){\tiny $0$}
        \put(83.25,19){\tiny $0$}
        \put(70.5,6){\tiny $0$}
        \put(70.5,31.85){\tiny $0$}
        \put(63,23){\small $-\tfrac{1}{2}$}
        \put(74,15){\small $-\tfrac{1}{2}$}
	\end{overpic}
	\caption{In the base case, we can add pairs of parallel arcs to the star-shaped region of $S$ shown at left, so that each region of $S \setminus \Gamma$ shown at right is a rectangle or a hexagon, for a total of $g-2$ hexagons. Here, and in all similar figures, $\Lambda$ is thick and $\gamma_i^{\pm}$ is thin.}
	\label{fig:star1}
\end{figure}

Now, suppose that $\Gamma \cap S$ contains an arc $a$ of weight one with at least one endpoint in $\Lambda$.  If $a$ does not cut off a half-bigon from $S$, we let $S \setminus a = S' \cup S''$, where $S'$ contains $g'$ boundary arcs of the form $\gamma_i^+$, $S''$ contains $g''$ boundary arcs of $\gamma_i^+$, and $g' + g'' = g$ (if $a$ meets $\Lambda$ twice) or $g'+g'' = g+1$ (if $a$ meets $\Lambda$ once).  We can treat $S'$ and $S''$ as star-shaped regions and apply the inductive hypothesis, since $S'$ and $S''$ contain fewer arcs that meet $\Lambda$ than $S$ does.  By induction, we have $\tb(S') - \rot(S') \leq \frac{g'}{2} - 1$ and $\tb(S'') - \rot(S'') \leq \tfrac{g''}{2} - 1$.  Ignoring the contribution to 
$\tb(S)$ coming from $a$, we have $\tb(S) \leq \tb(S') + \tb(S'')$.  Thus,
\begin{align*}
    \tb(S) - \rot(S) &\leq (\tb(S') - \rot(S')) + (\tb(S'') - \rot(S'')) \\
    &\leq (\tfrac{g'}{2} - 1) + (\tfrac{g''}{2} - 1) \\
    &\leq \tfrac{g+1}{2} - 2\\ 
    &< \tfrac{g}{2}-1.
\end{align*}

See \cref{fig:star2} and \cref{fig:star3} for examples of possible $S$, $S'$, and $S''$.

\begin{figure}[ht]
	\centering
	\begin{overpic}[scale=.35]{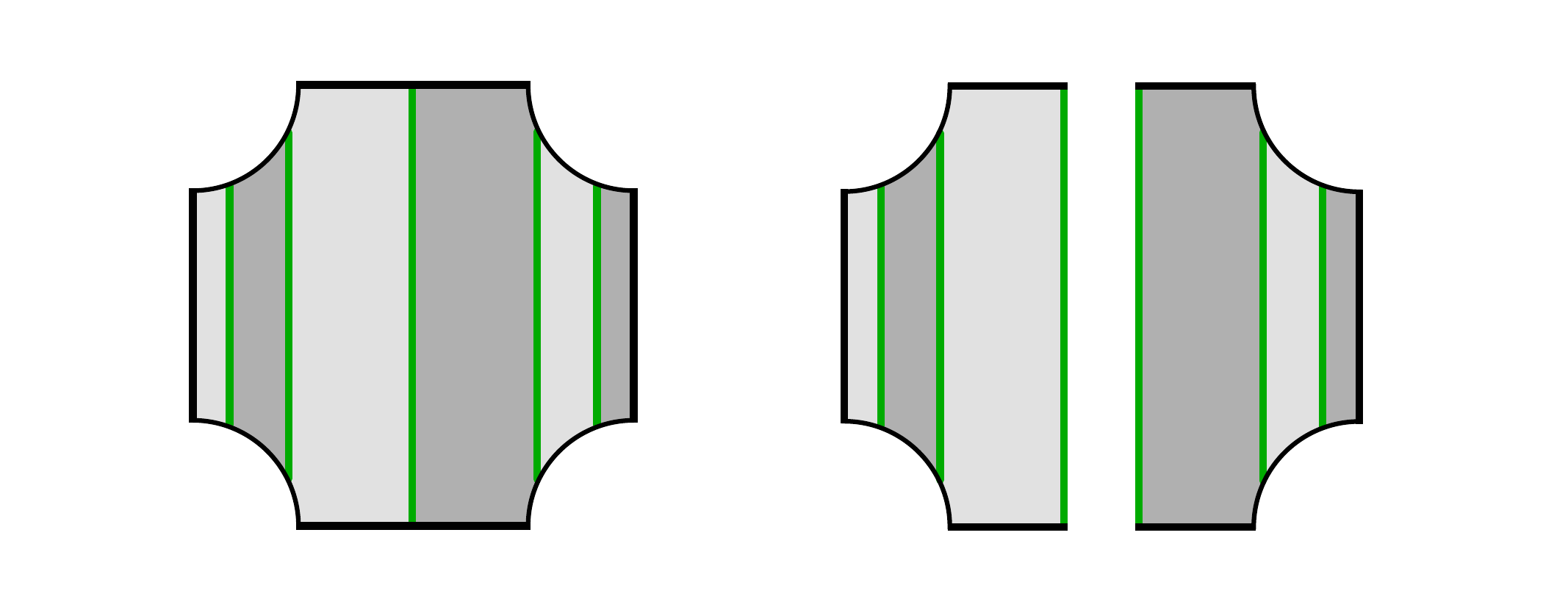}
        \put(12.75,19){\tiny $0$}
        \put(16,19){\tiny $0$}
        \put(22,19){\tiny $0$}
        \put(30,19){\tiny $0$}
        \put(35.75,19){\tiny $0$}
        \put(38.75,19){\tiny $0$}
        \put(25.25,3){\small \textcolor{darkgreen}{$a$}}

        \put(54.25,19){\tiny $0$}
        \put(57.5,19){\tiny $0$}
        \put(63.5,19){\tiny $0$}
        \put(76,19){\tiny $0$}
        \put(81.75,19){\tiny $0$}
        \put(85.25,19){\tiny $0$}
	\end{overpic}
	\caption{If $\Gamma \cap S$ contains an arc $a$ with both endpoints in $\Lambda$ (left), we cut $S$ along $a$ to get $S'$ and $S''$ (right).}
	\label{fig:star2}
\end{figure}

\begin{figure}[ht]
	\centering
	\begin{overpic}[scale=.35]{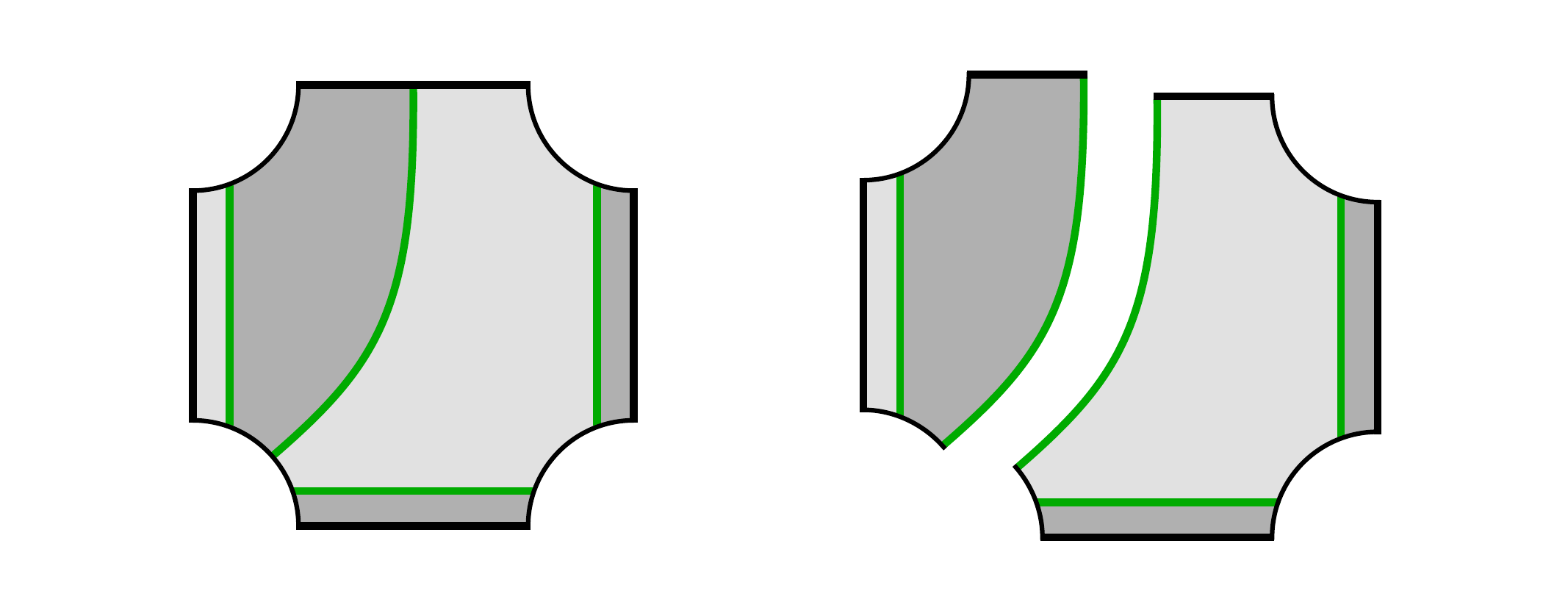}

        \put(28,18){\small $-\tfrac{1}{2}$}
        \put(13,19){\tiny $0$}
        \put(38.75,19){\tiny $0$}
        \put(26,6){\tiny $0$}
        \put(18.5,21){\small $0$}
        \put(25.5,34.75){\small \textcolor{darkgreen}{$b$}}

        \put(76,18){\small $-\tfrac{1}{2}$}
        \put(86,18.25){\tiny $0$}
        \put(55.75,19.25){\tiny $0$}
        \put(74,5.25){\tiny $0$}
        \put(61.5,21){\small $0$}
	\end{overpic}
	\caption{If $\Gamma \cap S$ contains an arc $b$ with one endpoint in $K$ and one endpoint in $\gamma_i^+$ (left), we cut $S$ along $b$ to get $S'$ and $S''$ (right).}
	\label{fig:star3}
\end{figure}

On the other hand, if the arc $a$ cuts off a half-bigon $R$ from $S$, we let $S' = S \setminus R$, so that we can treat $S'$ as its own star-shaped region containing one fewer arc that meets $\Lambda$.  By the inductive hypothesis, $\tb(S') - \rot(S') \leq \tfrac{g}{2} - 1$.  Since $R$ is a half-bigon, we have $\rot(R) = \pm \tfrac{1}{2}$, which implies that $\rot(S) = \rot(S') \pm \tfrac{1}{2}$ and $\rot(S) \geq \rot(S') - \tfrac{1}{2}$.  Since $a$ meets $\Lambda$ once, $\tb(S) = \tb(S') - \tfrac{1}{2}$.  Thus,
\[ \tb(S) - \rot(S) \leq (\tb(S') - \tfrac{1}{2}) + (\rot(S') - \tfrac{1}{2}) = \tb(S') - \rot(S') \leq \tfrac{g}{2}-1.\]

Finally, if $\Gamma \cap S$ contains arcs meetings $\Lambda$ with weight more than one, then by \cref{lem:parallel} we can remove a parallel pair, which preserves $\rot(S)$ and increases $\tb(S)$ by one; once again, the inductive hypothesis applied after removal yields $\tb(S) - \rot(S) \leq \frac{g}{2} - 1$, completing the proof.
\end{proof}

\begin{proposition}\label{prop:star}
Suppose $D$ is a normalized ribbon-convex star-shaped disk satisfying the principles above.  Then $D$ admits a max-sl transverse unlink ribbon derivative; c.f.\ \cref{fig:cellar}.
\end{proposition}

\begin{proof}
For each subdisk $D_i$, we examine possible intersection patterns of $\Gamma \cap D_i$ subject to the constraints above.  Our strategy is to show that
\begin{equation}\label{eq:index}
I(D_i) + \tb(D_i) -\rot(D_i) \leq -\frac{1}{2}
\end{equation}
for all possible subdisks.  We invoke \cref{lem:I} and \cref{lem:chi} to obtain 
\begin{align*}
    -1 &= \tb(K) - \rot(D) = (\tb(S) - \rot(S)) + \sum (\tb(D_i) - \rot(D_i))\\ 
    &= (\tb(S) - \rot(S)) + \sum( I(D_i) + \tb(D_i) - \rot(D_i)).
\end{align*}
Then, if we assume \eqref{eq:index} and use \cref{lem:star}, we have
\[(\tb(S) - \rot(S)) + \sum( I(D_i) + \tb(D_i) - \rot(D_i)) \leq (\tfrac{g}{2}-1) - \tfrac{g}{2} = -1.\]
If some $D_i$ satisfies $I(D_i) + \tb(D_i) -\rot(D_i) < -\tfrac{1}{2}$, the inequality above becomes strict, yielding the contradictory $-1 < -1$, from which it follows that $I(D_i) + \tb(D_i) -\rot(D_i) = -\tfrac{1}{2}$ for all subdisks $D_i$.

To this end, we call a subdisk $D_i$ with $I(D_i) + \tb(D_i) -\rot(D_i) = -\tfrac{1}{2}$ \emph{admissible}.  To prove \eqref{eq:index}, we split possible subdisks $D_i$ into cases based on the types of arcs in $\Gamma \cap D_i$ that meet $\partial \gamma_i^-$.  In the first six cases, we suppose that the two endpoints of $\gamma_i^-$ are contained in distinct arcs, and in the final three cases, we suppose these endpoints are contained in the same arc.  By \cref{lem:parallel}, in order to compute $\rot(D_i)$, we need only consider the parities of the weights of each collection of parallel arcs.  Additionally, because our goal is to compute an upper bound on $I(D_i) + \tb(D_i) -\rot(D_i)$, in each case we find a lower bound for $\rot(D_i)$, choosing negative and positive regions so that $\rot(D_i) \leq 0$ (if not, we could transpose negative and positive regions).

\vspace{2mm}
\textsc{Case \RNum{1}:} \textit{Both boundary points of $\gamma_i^-$ are contained in distinct $\Lambda$-arcs of $\Gamma \cap D_i$.}
\vspace{2mm}

There are three sub-types in this case.  First, suppose the boundary points of $\gamma_i^-$ are contained in non-parallel arcs, with $m$ total arcs of $\Gamma \cap D_i$ parallel to one of these, and $n$ total arcs parallel to the other.  It is possible that $\Gamma \cap D_i$ contains an additional collection of parallel arcs, so suppose these additional arcs are a collection of $p$ parallel seams, each of which meets $\gamma_i^-$ once (Type $\text{\RNum{1}}_{m,n,p}$, $m,n,p > 0$).  To compute $\rot(D_i)$ for Type $\text{\RNum{1}}_{m,n,p}$, we consider the eight possibilities corresponding to the parities of the parameters $m$, $n$, and $p$, shown in \cref{fig:type1a}.  By inspection, we see that in each of the eight cases, $\rot(D_i) \geq -\tfrac{3}{2}$.  We also compute $I(D_i) = m+n - 1$ and $\tb(D_i) = -m - n -\tfrac{p}{2}$.  It follows that $I(D_i) + \tb(D_i) -\rot(D_i) \leq -\tfrac{p}{2} - \tfrac{1}{2} < -\tfrac{1}{2}$, since $p > 0$, so $D_i$ is not admissible.

\begin{figure}[ht]
	\centering
	\begin{overpic}[scale=.35]{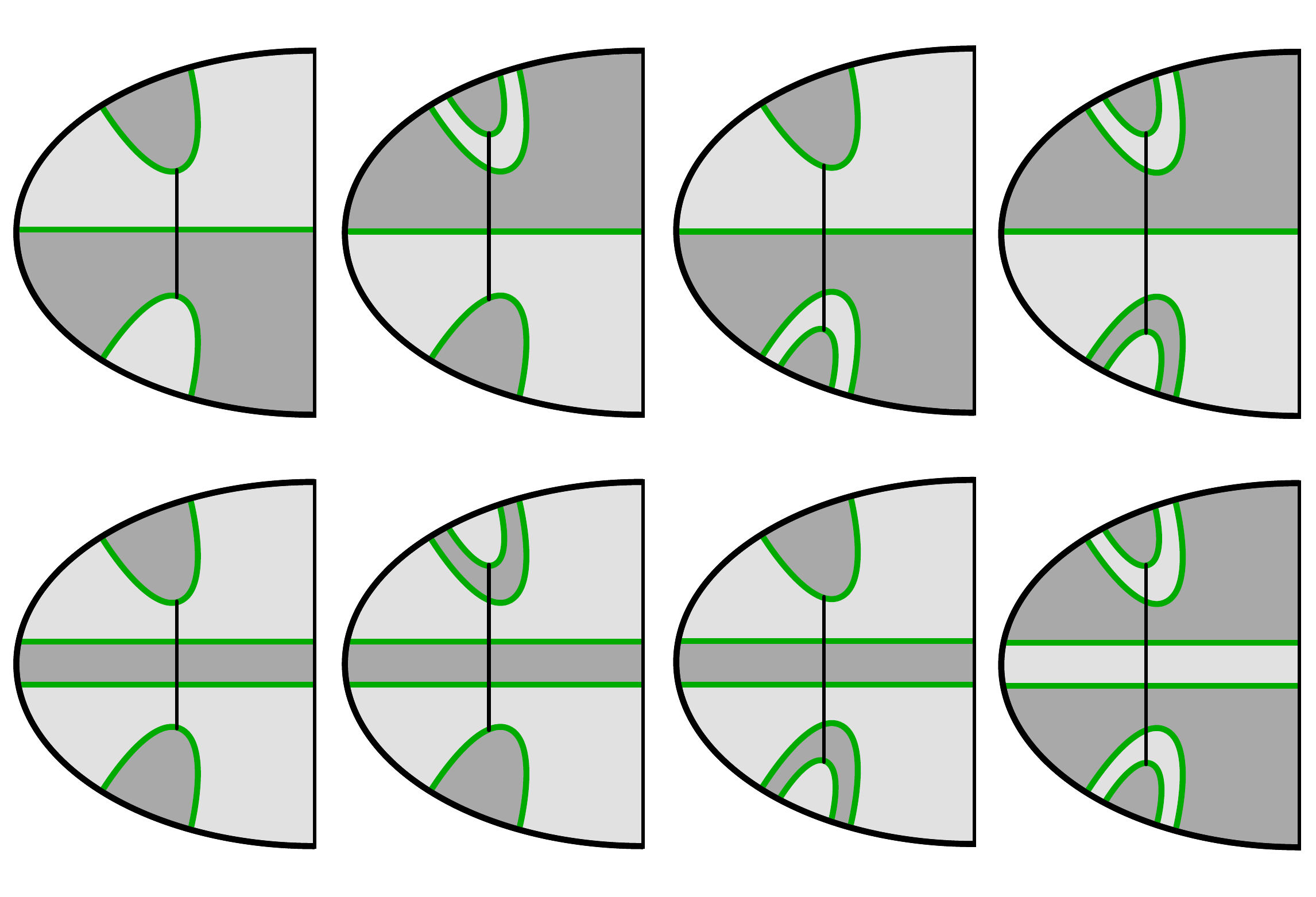}
        \put(-1,50){\small \textcolor{darkgreen}{$p$}}
        \put(6,39){\small \textcolor{darkgreen}{$n$}}
        \put(6,61){\small \textcolor{darkgreen}{$m$}}

        \put(6,46){\tiny $-\tfrac{1}{2}$}
        \put(6,54){\tiny $+\tfrac{1}{2}$}
        \put(11,41){\tiny $+1$}
        \put(11,59){\tiny $-1$}

        \put(30,46){\tiny $+\tfrac{1}{2}$}
        \put(30,54){\tiny $-\tfrac{1}{2}$}
        \put(35.25,60){\tiny $-1$}
        \put(37.25,56.25){\tiny $+1$}
        \put(35.75,41){\tiny $-1$}

        \put(56,46){\tiny $-\tfrac{1}{2}$}
        \put(56,54){\tiny $+\tfrac{1}{2}$}
        \put(60.75,40.25){\tiny $-1$}
        \put(63,43.25){\tiny $+1$}
        \put(61.5,58){\tiny $-1$}

        \put(80.5,46){\tiny $+\tfrac{1}{2}$}
        \put(80.5,54){\tiny $-\tfrac{1}{2}$}
        \put(85.25,40){\tiny $+1$}
        \put(87.25,43.25){\tiny $-1$}
        \put(85.25,60){\tiny $-1$}
        \put(87.25,56.25){\tiny $+1$}

        \put(6,13){\tiny $+\tfrac{1}{2}$}
        \put(6,17.25){\tiny $-\tfrac{1}{2}$}
        \put(6,22){\tiny $+\tfrac{1}{2}$}
        \put(11,9){\tiny $-1$}
        \put(11,26){\tiny $-1$}

        \put(31,13){\tiny $+\tfrac{1}{2}$}
        \put(31,17.25){\tiny $-\tfrac{1}{2}$}
        \put(31,22){\tiny $+\tfrac{1}{2}$}
        \put(36,9){\tiny $-1$}
        \put(35.5,26.5){\tiny $+1$}
        \put(37.5,24){\tiny $-1$}

        \put(56,13){\tiny $+\tfrac{1}{2}$}
        \put(56,17.25){\tiny $-\tfrac{1}{2}$}
        \put(56,22){\tiny $+\tfrac{1}{2}$}
        \put(60.75,8){\tiny $+1$}
        \put(62.75,10.75){\tiny $-1$}
        \put(60.5,26){\tiny $-1$}

        \put(81,13){\tiny $-\tfrac{1}{2}$}
        \put(81,17.25){\tiny $+\tfrac{1}{2}$}
        \put(81,22){\tiny $-\tfrac{1}{2}$}
        \put(85.5,8){\tiny $-1$}
        \put(87.5,10.25){\tiny $+1$}
        \put(85.5,26.5){\tiny $-1$}
        \put(87.5,24){\tiny $+1$}
       
	\end{overpic}
	\caption{Subdisks of Type $\text{\RNum{1}}_{m,n,p}$ split into cases based on parities of $m$, $n$, and $p$.}
	\label{fig:type1a}
\end{figure}

For the second sub-type, Type $\text{\RNum{1}}'_{m,n,p}$ with $m,n > 0$, suppose that the boundary points of $\gamma_i^-$ are contained in non-parallel arcs and $\Gamma \cap D_i$ contains a collection of $p$ additional parallel $\Lambda$-arcs that avoid $\gamma_i^-$.  Since $\Gamma$ necessarily meets $\gamma_i^+$, there must be one additional seam that avoids $\gamma_i^-$.  Cases by parities of $m$, $n$, and $p$ are shown in \cref{fig:type1b}.  Observe that $\rot(D_i) \geq -2$ (occurring when $m$ and $n$ are odd and $p$ is even); otherwise $\rot(D_i) \geq -1$ in all other cases.  As above, we compute $I(D_i) = m+n-2$ and $\tb(D_i) = -m - n - p - \tfrac{1}{2}$.  Thus, if $m$ and $n$ are odd and $p$ is even, we have $I(D_i) + \tb(D_i) -\rot(D_i) \leq -2 -p -\tfrac{1}{2} + 2 = -p -\tfrac{1}{2}$, so that $D_i$ is admissible if and only if $p=0$.  Otherwise, we have $I(D_i) + \tb(D_i) -\rot(D_i) \leq -2 - p -\tfrac{1}{2} +1 = - p-\tfrac{3}{2} < -\tfrac{1}{2}$, and $D_i$ is not admissible.

\begin{figure}[ht]
	\centering
	\begin{overpic}[scale=.35]{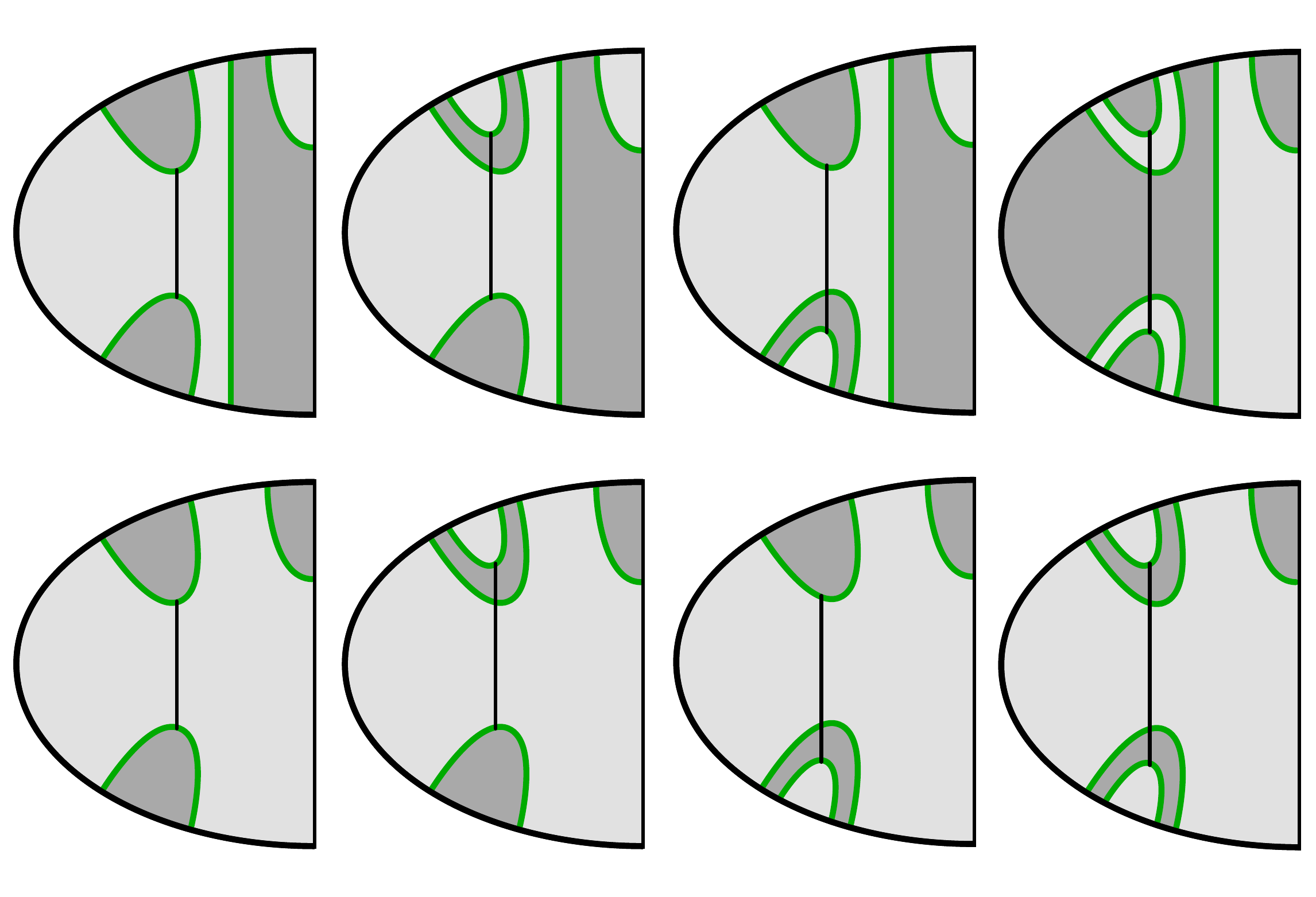}
        \put(16.5,65.5){\small \textcolor{darkgreen}{$p$}}
        \put(6,39){\small \textcolor{darkgreen}{$n$}}
        \put(6,61){\small \textcolor{darkgreen}{$m$}}

        \put(19.5,50){\tiny $-\tfrac{1}{2}$}
        \put(6,50){\tiny $+1$}
        \put(11,41){\tiny $-1$}
        \put(11,59){\tiny $-1$}
        \put(21,61.25){\tiny $+\tfrac{1}{2}$}

        \put(30,50){\tiny $+1$}
        \put(44,50){\tiny $-\tfrac{1}{2}$}
        \put(35.25,60){\tiny $+1$}
        \put(37.75,56.5){\tiny $-1$}
        \put(35.75,41){\tiny $-1$}
        \put(46,61.25){\tiny $+\tfrac{1}{2}$}

        \put(70,50){\tiny $-\tfrac{1}{2}$}
        \put(56,50){\tiny $+1$}
        \put(60.75,40.25){\tiny $+1$}
        \put(63,43.5){\tiny $-1$}
        \put(61.5,58){\tiny $-1$}
        \put(71.25,61.25){\tiny $+\tfrac{1}{2}$}

        \put(94.5,50){\tiny $+\tfrac{1}{2}$}
        \put(80.5,50){\tiny $-1$}
        \put(85.25,40){\tiny $+1$}
        \put(87.5,43.25){\tiny $-1$}
        \put(85.25,60){\tiny $+1$}
        \put(87.5,56.25){\tiny $-1$}
        \put(96,61.25){\tiny $-\tfrac{1}{2}$}

        \put(6,17.25){\tiny $+\tfrac{1}{2}$}
        \put(11,9){\tiny $-1$}
        \put(11,26){\tiny $-1$}
        \put(21,27.5){\tiny $-\tfrac{1}{2}$}

        \put(31,17.25){\tiny $+\tfrac{1}{2}$}
        \put(36,9){\tiny $-1$}
        \put(35.5,26.5){\tiny $+1$}
        \put(37.75,23.75){\tiny $-1$}
        \put(46,27.5){\tiny $-\tfrac{1}{2}$}

        \put(56,17.25){\tiny $+\tfrac{1}{2}$}    
        \put(60.75,8){\tiny $+1$}
        \put(63,10.75){\tiny $-1$}
        \put(60.5,26){\tiny $-1$}
        \put(71.25,27.5){\tiny $-\tfrac{1}{2}$}

        \put(81,17.25){\tiny $+\tfrac{1}{2}$}
        \put(85.5,8){\tiny $-1$}
        \put(87.5,10.25){\tiny $+1$}
        \put(85.5,26.5){\tiny $-1$}
        \put(87.5,23.75){\tiny $+1$}
        \put(96,27.5){\tiny $-\tfrac{1}{2}$}
       
	\end{overpic}
	\caption{Subdisks of Type $\text{\RNum{1}}'_{m,n,p}$ split into cases based on parities of $m$, $n$, and $p$.}
	\label{fig:type1b}
\end{figure}

For the third sub-type, Type $\text{\RNum{1}}''_{n,p}$ with $n > 1$, we suppose the boundary points of $\gamma_i^-$ are contained in distinct parallel arcs, with a total weight of $n > 1$, along with an additional collection of $p$ parallel $\Lambda$-arcs avoiding $\gamma_i^-$.  As in the previous type, there must be one more seam avoiding $\gamma_i^-$.  Cases are shown in \cref{fig:type1c}.  We compute $\rot(D_i) \geq -1$, with $I(D_i) = n-2$ and $\tb(D_i) = -n -p - \tfrac{1}{2}$.  Thus, $I(D_i) + \tb(D_i) -\rot(D_i) \leq  -2 - p -\tfrac{1}{2} + 1 < -\tfrac{1}{2}$, and $D_i$ is not admissible.

\begin{figure}[ht]
	\centering
	\begin{overpic}[scale=.35]{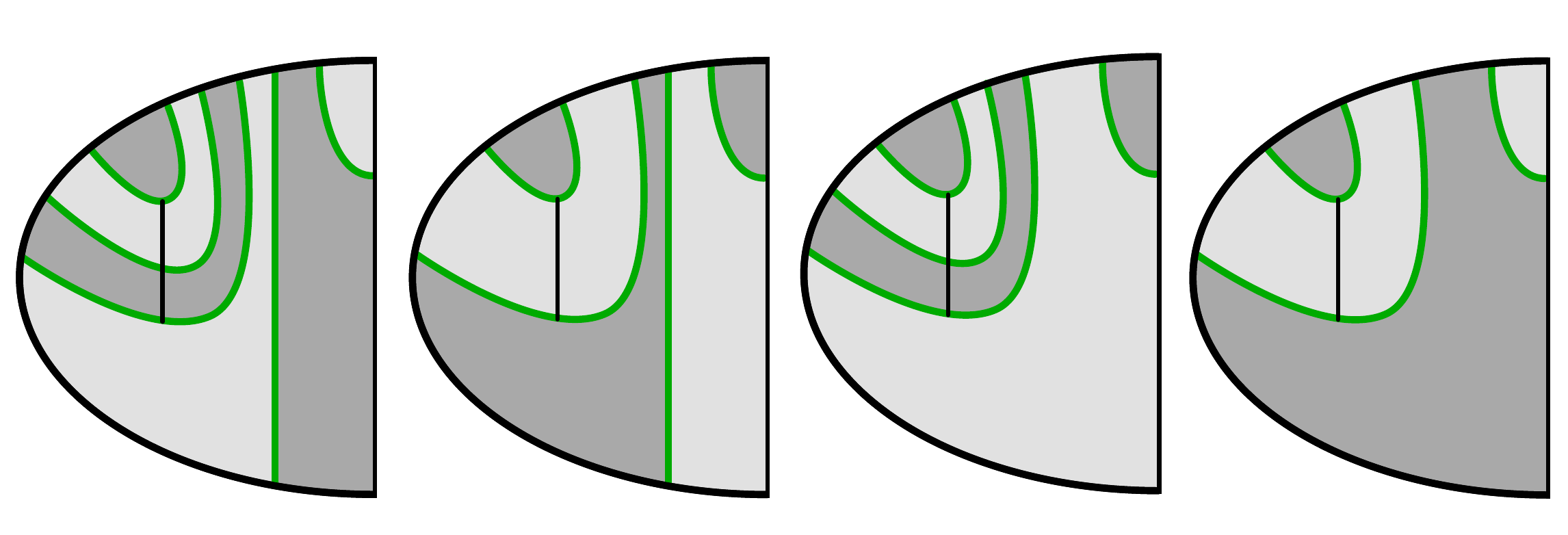}
        \put(16.75,0.25){\small \textcolor{darkgreen}{$p$}}
        \put(5,28){\small \textcolor{darkgreen}{$n$}}

        \put(19,16){\tiny $-\tfrac{1}{2}$}
        \put(11,9){\tiny $+1$}
        \put(21,27.5){\tiny $+\tfrac{1}{2}$}
        \put(8,24){\tiny $-1$}
        \put(5.5,20){\tiny $+1$}
        \put(3,17){\tiny $-1$}

        \put(45,16){\tiny $+\tfrac{1}{2}$}
        \put(36,9){\tiny $-1$}
        \put(33.5,23.5){\tiny $-1$}
        \put(46,27.5){\tiny $-\tfrac{1}{2}$}
        \put(30.5,18){\tiny $+1$}

        \put(69,16){\tiny $+\tfrac{1}{2}$}    
        \put(71.25,27.5){\tiny $-\tfrac{1}{2}$}
        \put(58,24){\tiny $-1$}
        \put(56.5,20){\tiny $+1$}
        \put(55,17){\tiny $-1$}

        \put(94,16){\tiny $-\tfrac{1}{2}$}
        \put(96,27.5){\tiny $-\tfrac{1}{2}$}
        \put(80.5,18){\tiny $+1$}
        \put(83,24){\tiny $-1$}
       
	\end{overpic}
	\caption{Subdisks of Type $\text{\RNum{1}}''_{n,p}$ split into cases based on the parities of $n$ and $p$.}
	\label{fig:type1c}
\end{figure}

\vspace{2mm}
\textsc{Case \RNum{2}:} \textit{Both boundary points of $\gamma_i^-$ are contained in different $\gamma_i^+$-arcs.}
\vspace{2mm}

As with the first case, there are three sub-types to consider.  First, suppose that the boundary points of $\gamma_i^-$ are contained in non-parallel $\gamma_i^+$-arcs, with $m$ total arcs of $\Gamma \cap D_i$ parallel to one of these, and $n$ total arcs parallel to the other.  It is possible that $\Gamma \cap D_i$ contains an additional collection of parallel arcs, so suppose these additional arcs are a collection of $p$ parallel seams, each of which meets $\gamma_i^-$ once (Type $\text{\RNum{2}}_{m,n,p}$, $m,n> 0$).  As above, to compute $\rot(D_i)$ we consider eight possibilities, shown in \cref{fig:type2a}.  We compute $\rot(D_i) \geq -\tfrac{3}{2}$, $I(D_i) = -m - n -1$, and $\tb(D_i) = -\tfrac{p}{2}$.  It follows that $I(D_i) + \tb(D_i) -\rot(D_i) \leq -m - n - \tfrac{p}{2} + \tfrac{1}{2} < -\tfrac{1}{2}$, since $m,n > 0$, and $D_i$ is not admissible.

\begin{figure}[ht]
	\centering
	\begin{overpic}[scale=.35]{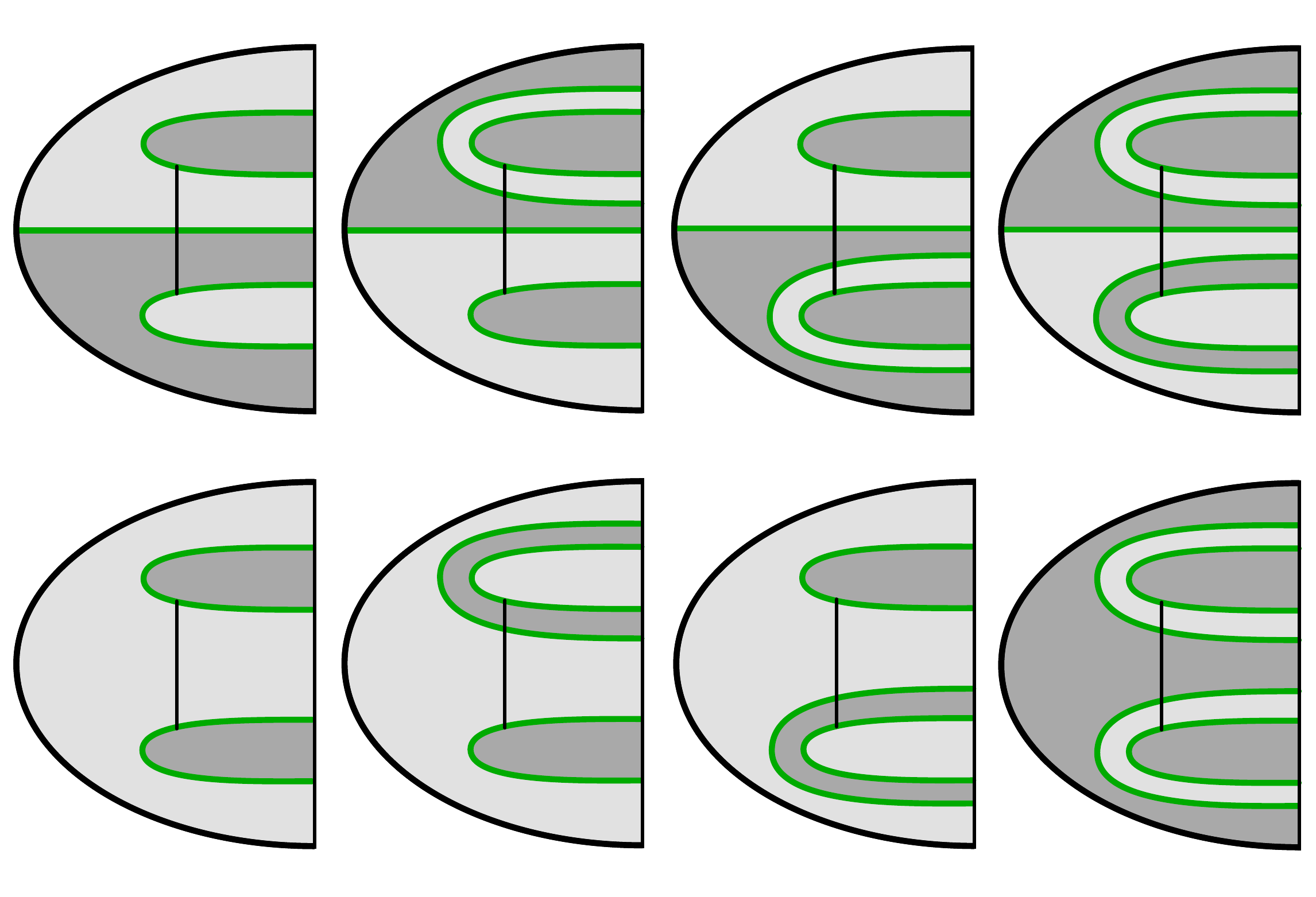}
        \put(-1,50){\small \textcolor{darkgreen}{$p$}}
        \put(25,41){\small \textcolor{darkgreen}{$n$}}
        \put(25,59){\small \textcolor{darkgreen}{$m$}}

        \put(18,43.5){\tiny $+\tfrac{1}{2}$}
        \put(6,53){\tiny $+0$}
        \put(6,47){\tiny $-0$}
        \put(18,57){\tiny $-\tfrac{1}{2}$}

        \put(42,43.5){\tiny $-\tfrac{1}{2}$}
        \put(30,53){\tiny $-0$}
        \put(30,47){\tiny $+0$}
        \put(42,57){\tiny $-\tfrac{1}{2}$}
        \put(42,53.4){\tiny $+0$}

        \put(68,43.5){\tiny $-\tfrac{1}{2}$}
        \put(56,53){\tiny $+0$}
        \put(56,47){\tiny $-0$}
        \put(68,57){\tiny $-\tfrac{1}{2}$}
        \put(68,47){\tiny $+0$}

        \put(92,43.5){\tiny $+\tfrac{1}{2}$}
        \put(80,53){\tiny $-0$}
        \put(80,47){\tiny $+0$}
        \put(92,57){\tiny $-\tfrac{1}{2}$}
        \put(92,53.25){\tiny $+0$}
        \put(92,47){\tiny $-0$}

        \put(6,17.25){\tiny $-\tfrac{1}{2}$}
        \put(18,10.5){\tiny $-\tfrac{1}{2}$}
        \put(18,24){\tiny $-\tfrac{1}{2}$}

        \put(31,17.25){\tiny $-\tfrac{1}{2}$}
        \put(42,10.5){\tiny $-\tfrac{1}{2}$}
        \put(42,24){\tiny $+\tfrac{1}{2}$}
        \put(42,20.2){\tiny $-0$}

        \put(56,17.25){\tiny $-\tfrac{1}{2}$} 
        \put(68,10.5){\tiny $+\tfrac{1}{2}$}
        \put(68,24){\tiny $-\tfrac{1}{2}$}
        \put(68,14){\tiny $-0$}

        \put(81,17.25){\tiny $+\tfrac{1}{2}$}
        \put(92,10.5){\tiny $-\tfrac{1}{2}$}
        \put(92,24){\tiny $-\tfrac{1}{2}$}
        \put(92,20.2){\tiny $+0$}
        \put(92,13.85){\tiny $+0$}
       
	\end{overpic}
	\caption{Subdisks of Type $\text{\RNum{2}}_{m,n,p}$ split into cases based on parities of $m$, $n$, and $p$.}
	\label{fig:type2a}
\end{figure}

For the second sub-type, Type $\text{\RNum{2}}'_{m,n,p}$ with $m,n > 0$, suppose boundary points of $\gamma_i^-$ are contained in non-parallel arcs and $\Gamma \cap D_i$ contains a collection of $p$ additional parallel $\gamma_i^+$-arcs that avoid $\gamma_i^-$.  Cases of this type in which $p$ is even agree with the first sub-type and are shown in the lower row of \cref{fig:type2a}.  Cases by parities of $m$ and $n$ with $p$ odd are shown in \cref{fig:type2b}.  We compute $\rot(D_i) \geq -\tfrac{3}{2}$, $I(D_i) = -m-n-2p-1$, and $\tb(D_i) = 0$.  It follows that $I(D_i) + \tb(D_i) -\rot(D_i) \leq -m-n-2p + \tfrac{1}{2} < -\tfrac{1}{2}$, since $m,n > 0$, and $D_i$ is not admissible.

\begin{figure}[ht]
	\centering
	\begin{overpic}[scale=.35]{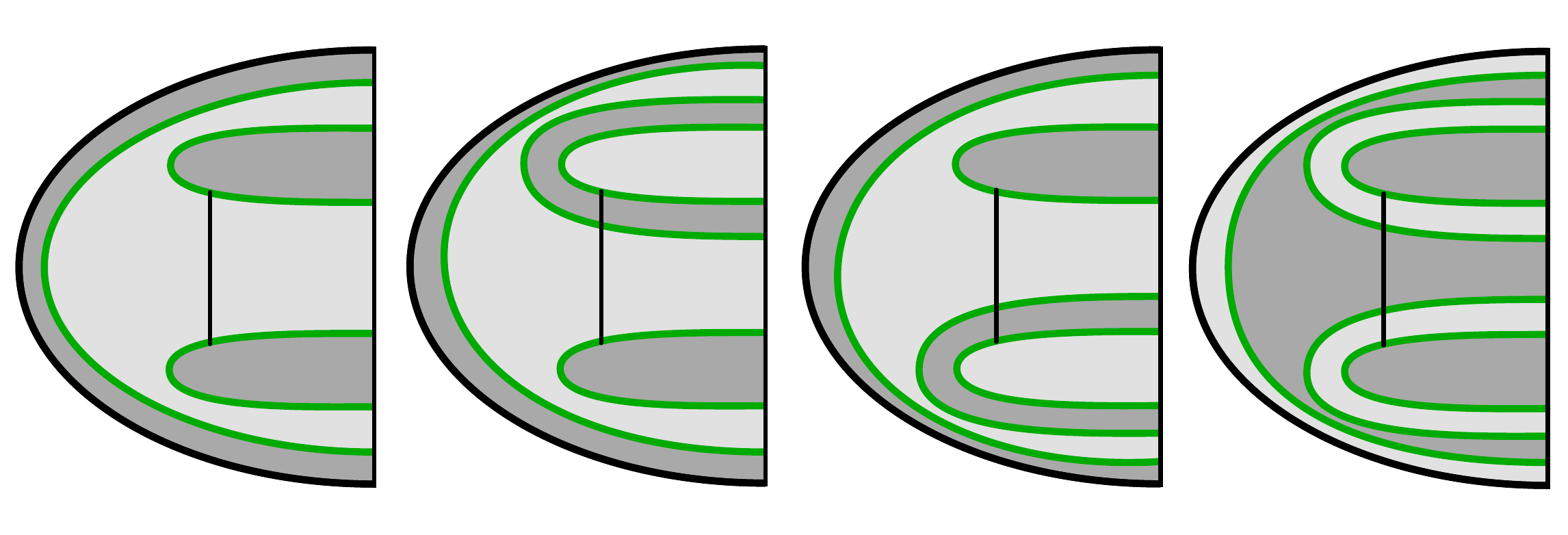}
        \put(25,5){\small \textcolor{darkgreen}{$p$}}
        \put(25,12){\small \textcolor{darkgreen}{$n$}}
        \put(25,25){\small \textcolor{darkgreen}{$m$}}
        
        \put(6,17){\tiny $-\tfrac{1}{2}$}
        \put(18,10){\tiny $-\tfrac{1}{2}$}
        \put(18,23){\tiny $-\tfrac{1}{2}$}
        \put(6,5){\tiny $-0$}
        
        \put(31,17){\tiny $-\tfrac{1}{2}$}
        \put(42,10){\tiny $-\tfrac{1}{2}$}
        \put(42,23){\tiny $+\tfrac{1}{2}$}
        \put(42,19.75){\tiny $-0$}
        \put(31,5){\tiny $-0$}

        \put(56,17){\tiny $-\tfrac{1}{2}$} 
        \put(68,10){\tiny $+\tfrac{1}{2}$}
        \put(68,23){\tiny $-\tfrac{1}{2}$}
        \put(68,13.5){\tiny $-0$}
        \put(56,5){\tiny $-0$}

        \put(81,17){\tiny $+\tfrac{1}{2}$}
        \put(92,10){\tiny $-\tfrac{1}{2}$}
        \put(92,23){\tiny $-\tfrac{1}{2}$}
        \put(92,19.25){\tiny $+0$}
        \put(92,13.25){\tiny $+0$}
        \put(81,5){\tiny $+0$}
       
	\end{overpic}
	\caption{Subdisks of Type $\text{\RNum{2}}'_{m,n,p}$ split into cases based on parities of $m$ and $n$, with $p$ odd.}
	\label{fig:type2b}
\end{figure}

For the third sub-type, Type $\text{\RNum{2}}''_{n,p}$ with $n > 1$, we suppose the boundary points of $\gamma_i^-$ are contained in distinct parallel arcs, with a total weight of $n > 1$, along with an additional collection of $p$ parallel $\gamma_i^+$-arcs avoiding $\gamma_i^-$.  Cases are shown in \cref{fig:type2c}.  We compute $\rot(D_i) \geq -\tfrac{1}{2}$, with $I(D_i) = -n-2p-1$ and $\tb(D_i) = 0$.  Thus, $I(D_i) + \tb(D_i) -\rot(D_i) \leq  -n - 2p -\tfrac{1}{2} < -\tfrac{1}{2}$, since $n>1$, and $D_i$ is not admissible.

\begin{figure}[ht]
	\centering
	\begin{overpic}[scale=.35]{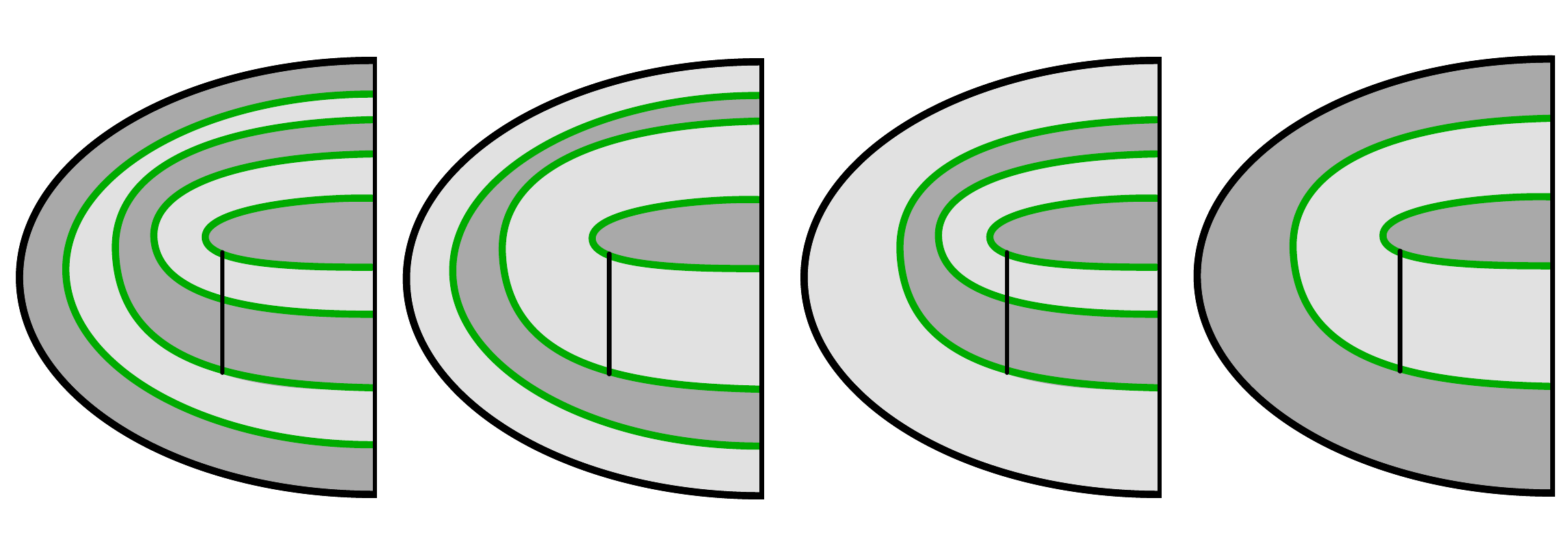}
        \put(25,5){\small \textcolor{darkgreen}{$p$}}
        \put(25,24){\small \textcolor{darkgreen}{$n$}}

        \put(18,7){\tiny $+0$}
        \put(18,11){\tiny $-0$}
        \put(18,15){\tiny $+0$}
        \put(18,19){\tiny $-\tfrac{1}{2}$}
        \put(6,8){\tiny $-0$}

        \put(31,8){\tiny $+0$}
        \put(42,7){\tiny $-0$}
        \put(42,13){\tiny $+0$}
        \put(42,18.5){\tiny $-\tfrac{1}{2}$}

        \put(68,11){\tiny $-0$}
        \put(68,15){\tiny $+0$}
        \put(68,19){\tiny $-\tfrac{1}{2}$}
        \put(60,8){\tiny $+0$}

        \put(93,13){\tiny $+0$}
        \put(93,19){\tiny $-\tfrac{1}{2}$}
        \put(85,8){\tiny $-0$}
       
	\end{overpic}
	\caption{Subdisks of Type $\text{\RNum{2}}''_{n,p}$ split into cases based on parities of $n$ and $p$.}
	\label{fig:type2c}
\end{figure}

\vspace{2mm}
\textsc{Case \RNum{3}:} \textit{One boundary point of $\gamma_i^-$ is contained in a $\Lambda$-arc, and the other is contained in a $\gamma_i^+$-arc.}
\vspace{2mm}

In this case, $\Gamma \cap D_i$ contains $m > 0$ arcs parallel to the $\Lambda$-arc, $n$ arcs parallel to the $\gamma_i^+$-arc, and $p$ parallel seams, each meeting $\gamma_i^-$ once (Type $\text{\RNum{3}}_{m,n,p}$ with $m,n>0$, shown at left in \cref{fig:type3}).  Observe that $\rot(D_i) \geq -\tfrac{3}{2}$ (occurring when $m$ and $n$ are odd and $p$ is even); otherwise $\rot(D_i) \geq -1$ in all other cases.  In addition, $I(D_i) = m-n-1$ and $\tb(D_i) = -m-\tfrac{p}{2}$.  Thus, if $m$ and $n$ are odd and $p$ is even, we have $I(D_i) + \tb(D_i) -\rot(D_i) \leq -n -\tfrac{p}{2} + \tfrac{1}{2} \leq -\tfrac{p}{2} - \tfrac{1}{2}$ (since $n>0$), so that $D_i$ is admissible if and only if $p=0$.  Otherwise, we have $I(D_i) + \tb(D_i) -\rot(D_i) \leq -n - \tfrac{p}{2} < -\tfrac{1}{2}$, and $D_i$ is not admissible.

\begin{figure}[ht]
	\centering
	\begin{overpic}[scale=.35]{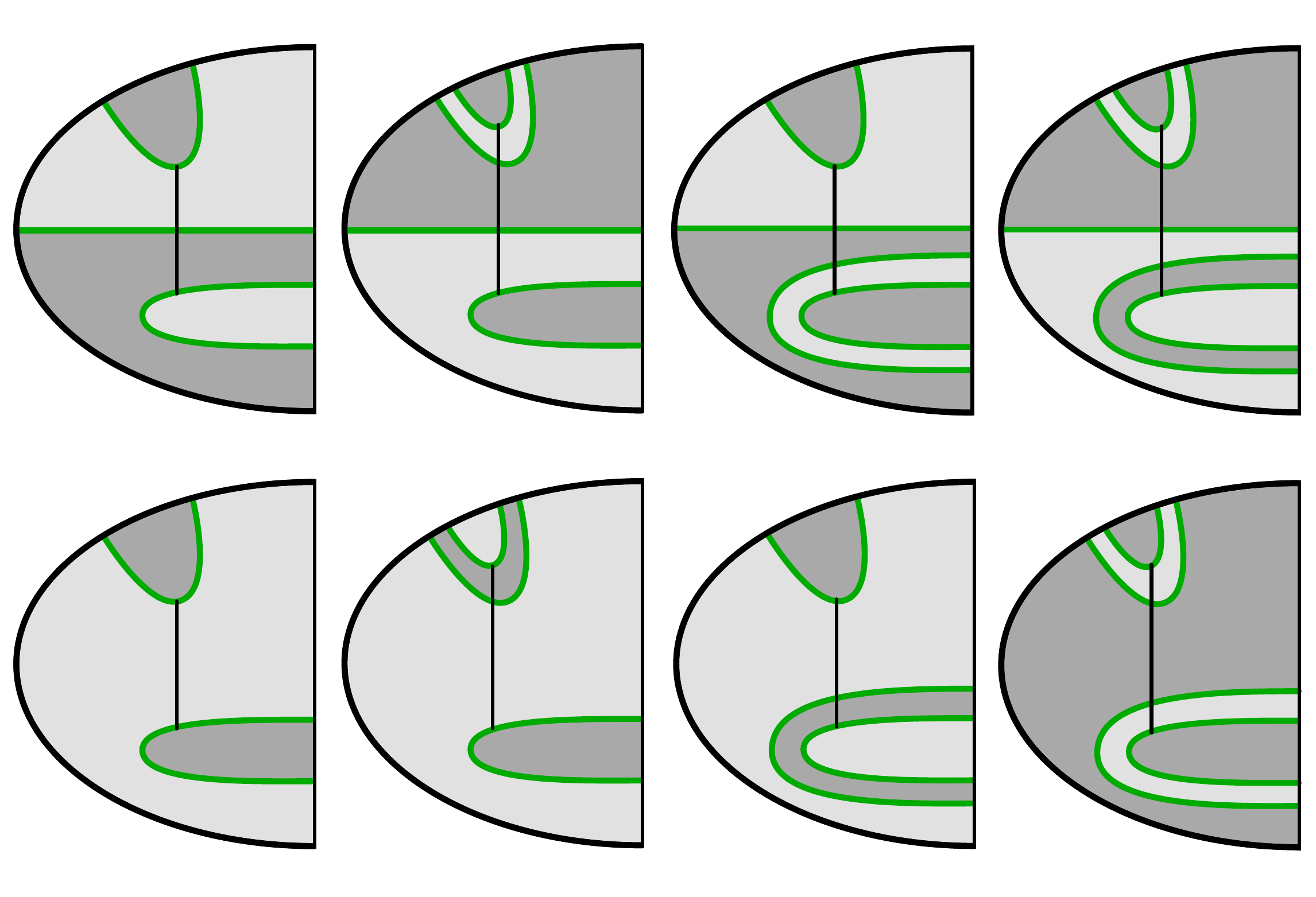}
        \put(-1,50){\small \textcolor{darkgreen}{$p$}}
        \put(25,41){\small \textcolor{darkgreen}{$n$}}
        \put(6,61.5){\small \textcolor{darkgreen}{$m$}}

        \put(18,43.5){\tiny $+\tfrac{1}{2}$}
        \put(11,59){\tiny $-1$}
        \put(6,47){\tiny $-0$}
        \put(18,57){\tiny $+\tfrac{1}{2}$}

        \put(42,43.5){\tiny $-\tfrac{1}{2}$}
        \put(30,53){\tiny $-0$}
        \put(30,47){\tiny $+0$}
        \put(36,60){\tiny $-1$}
        \put(38,57.5){\tiny $+1$}

        \put(68,43.5){\tiny $-\tfrac{1}{2}$}
        \put(61,59){\tiny $-1$}
        \put(56,47){\tiny $-0$}
        \put(68,57){\tiny $+\tfrac{1}{2}$}
        \put(68,47){\tiny $+0$}

        \put(92,43.5){\tiny $+\tfrac{1}{2}$}
        \put(80,47){\tiny $+0$}
        \put(92,57){\tiny $-\tfrac{1}{2}$}
        \put(92,47){\tiny $-0$}
        \put(86,60){\tiny $-1$}
        \put(88.5,57.5){\tiny $+1$}

        \put(11,26){\tiny $-1$}
        \put(18,10.5){\tiny $-\tfrac{1}{2}$}
        \put(18,21){\tiny $+0$}

        \put(42,10.5){\tiny $-\tfrac{1}{2}$}
        \put(42,20.2){\tiny $-0$}
        \put(35.5,26.5){\tiny $+1$}
        \put(38,24){\tiny $-1$}

        \put(61.5,26){\tiny $-1$} 
        \put(68,10.5){\tiny $+\tfrac{1}{2}$}
        \put(68,20.2){\tiny $+0$}
        \put(68,14){\tiny $-0$}

        \put(92,10.5){\tiny $-\tfrac{1}{2}$}
        \put(92,20.2){\tiny $-0$}
        \put(92,13.85){\tiny $+0$}
        \put(85.5,26.5){\tiny $-1$}
        \put(87.75,24){\tiny $+1$}
       
	\end{overpic}
	\caption{Subdisks of Type $\text{\RNum{3}}_{m,n,p}$ split into cases based on parities of $m$, $n$, and $p$.}
	\label{fig:type3}
\end{figure}

\vspace{2mm}
\textsc{Case \RNum{4}:} \textit{One boundary point of $\gamma_i^-$ is contained in a $\Lambda$-arc, and the other is contained in a seam.}
\vspace{2mm}

In this case, $\Gamma \cap D_i$ contains $m$ parallel $\Lambda$-arcs and $p$ parallel seams, yielding Type $\text{\RNum{4}}_{m,p}$, with $m,p > 0$.  Because our computations of $\rot(D_i)$, $I(D_i)$, and $\tb(D_i)$ do not depend on which arcs pass through $\partial \gamma_i^-$, we can use \cref{lem:parallel} and \cref{fig:type3} with $n = 0$.  Observe that $\rot(D_i) \geq -1$ (occurring when $m$ and $p$ are odd); otherwise $\rot(D_i) \geq -\tfrac{1}{2}$ in all other cases.  In addition, $I(D_i) = m-1$ and $\tb(D_i) = -m-\tfrac{p}{2}$.  Thus, if $m$ and $p$ are odd, we have $I(D_i) + \tb(D_i) -\rot(D_i) \leq -\tfrac{p}{2}$, so that $D_i$ is admissible if and only if $p=1$.  Otherwise, we have $I(D_i) + \tb(D_i) -\rot(D_i) \leq -\tfrac{p}{2}-\tfrac{1}{2} < -\tfrac{1}{2}$, and $D_i$ is not admissible.

\vspace{2mm}
\textsc{Case \RNum{5}:} \textit{One boundary point of $\gamma_i^-$ is contained in a $\gamma_i^+$-arc, and the other is contained in a seam.}
\vspace{2mm}

In this case, $\Gamma \cap D_i$ contains $n$ parallel $\gamma_i^+$-arcs and $p$ parallel seams, yielding Type $\text{\RNum{5}}_{n,p}$, with $n,p > 0$.  As above, we can use \cref{lem:parallel} and \cref{fig:type3} with $m = 0$.  We compute $\rot(D_i) \geq -1$, $I(D_i) = -n-1$ and $\tb(D_i) = -\tfrac{p}{2}$.  Thus $I(D_i) + \tb(D_i) -\rot(D_i) \leq -n-\tfrac{p}{2} < -\tfrac{1}{2}$ as $n,p > 0$, and $D_i$ is not admissible.

\vspace{2mm}
\textsc{Case \RNum{6}:} \textit{Both boundary points of $\gamma_i^-$ are contained in distinct seams.}
\vspace{2mm}

In this case, $\Gamma \cap D_i$ contains $p$ parallel seams, yielding Type $\text{\RNum{6}}_{p}$, with $p>1$.  As above, we can use \cref{lem:parallel} and \cref{fig:type3} with $m,n = 0$.  We compute $\rot(D_i) \geq -\tfrac{1}{2}$, $I(D_i) = -1$ and $\tb(D_i) = -\tfrac{p}{2}$.  Thus $I(D_i) + \tb(D_i) -\rot(D_i) \leq -\tfrac{p}{2}-\tfrac{1}{2} < -\tfrac{1}{2}$ as $p > 1$, and $D_i$ is not admissible.

\vspace{2mm}
\textsc{Case \RNum{7}:} \textit{Boundary points of $\gamma_i^-$ are contained in a single $\Lambda$-arc.}
\vspace{2mm}

In this case, there are $p$ parallel $\gamma_i^+$-arcs, with $p-1$ of these disjoint from $\gamma_i^-$, and one extra seam that avoids $\gamma_i^+$ (Type $\text{\RNum{7}}_p$, $p>0$).  Depictions of Type $\text{\RNum{7}}_p$, split based on the parity of $p$, are shown in \cref{fig:type7}.  We have $\rot(D_i) \geq -1$, $I(D_i) = 0$, and $\tb(D_i) = -p-\tfrac{1}{2}$.  Thus, $I(D_i) + \tb(D_i) -\rot(D_i) \leq -p + \tfrac{1}{2}$, and $D_i$ is admissible if and only if $p=1$.

\begin{figure}[ht]
	\centering
	\begin{overpic}[scale=.35]{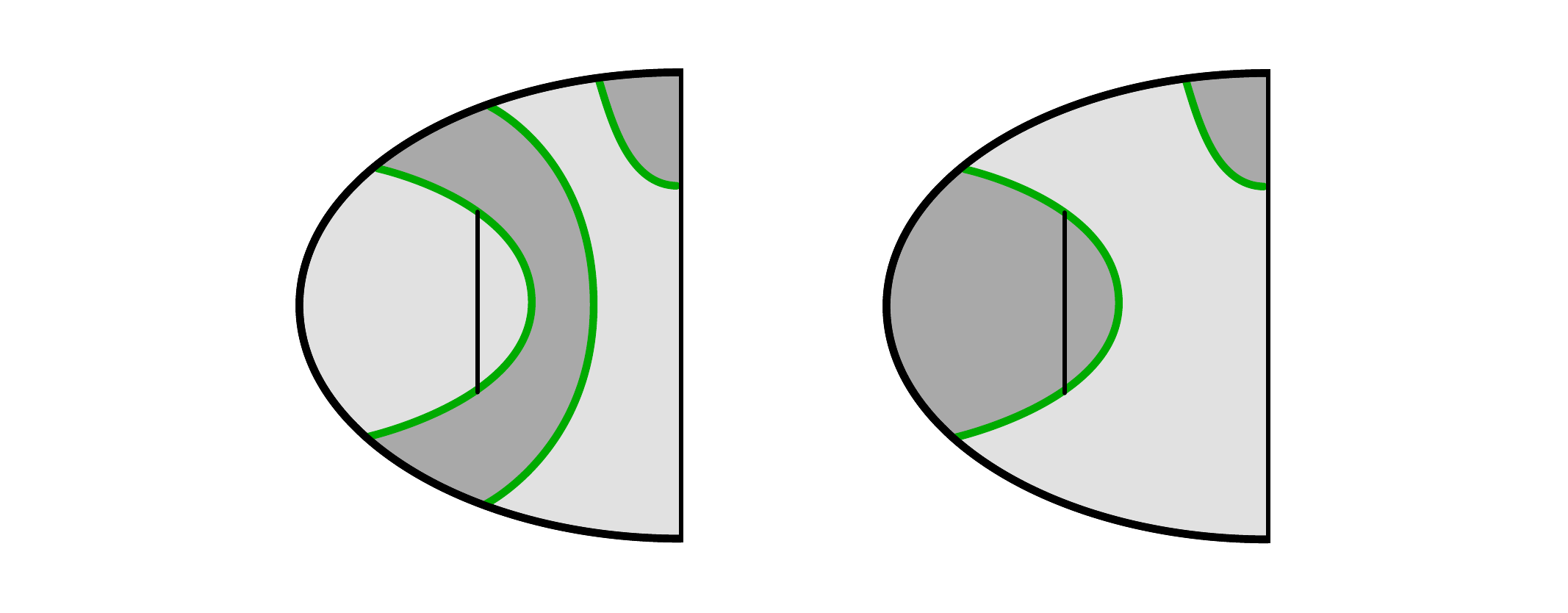}
        
        \put(29,34){\small \textcolor{darkgreen}{$p$}}

        \put(40,31){\tiny $-\tfrac{1}{2}$}
        \put(39.25,19){\tiny $+\tfrac{1}{2}$}
        \put(24,19){\tiny $+1$}
        \put(34.5,19){\tiny $-1$}

        \put(77.5,31){\tiny $-\tfrac{1}{2}$}
        \put(75,19){\tiny $+\tfrac{1}{2}$}
        \put(62,19){\tiny $-1$}

	\end{overpic}
	\caption{Subdisks of Type $\text{\RNum{7}}_{p}$ split into cases based on the parity of $p$.}
	\label{fig:type7}
\end{figure}

\vspace{2mm}
\textsc{Case \RNum{8}:} \textit{Boundary points of $\gamma_i^-$ are contained in a single $\gamma_i^+$-arc.}
\vspace{2mm}

In this case, there are $p$ parallel $\gamma_i^+$-arcs, with $p-1$ of these disjoint from $\gamma_i^-$, and Type $\text{\RNum{8}}_p$ (with $p > 0$) is shown at left and center in \cref{fig:type89}, split into cases based on the parity of $p$.  We compute $\rot(D_i) = -\tfrac{1}{2}$, $I(D_i) = -2p+1$, and $\tb(D_i) = 0$, so that $I(D_i) + \tb(D_i) -\rot(D_i) = -2p+\tfrac{3}{2}$.  It follows that $D_i$ is admissible if and only if $p=1$.

\begin{figure}[ht]
	\centering
	\begin{overpic}[scale=.35]{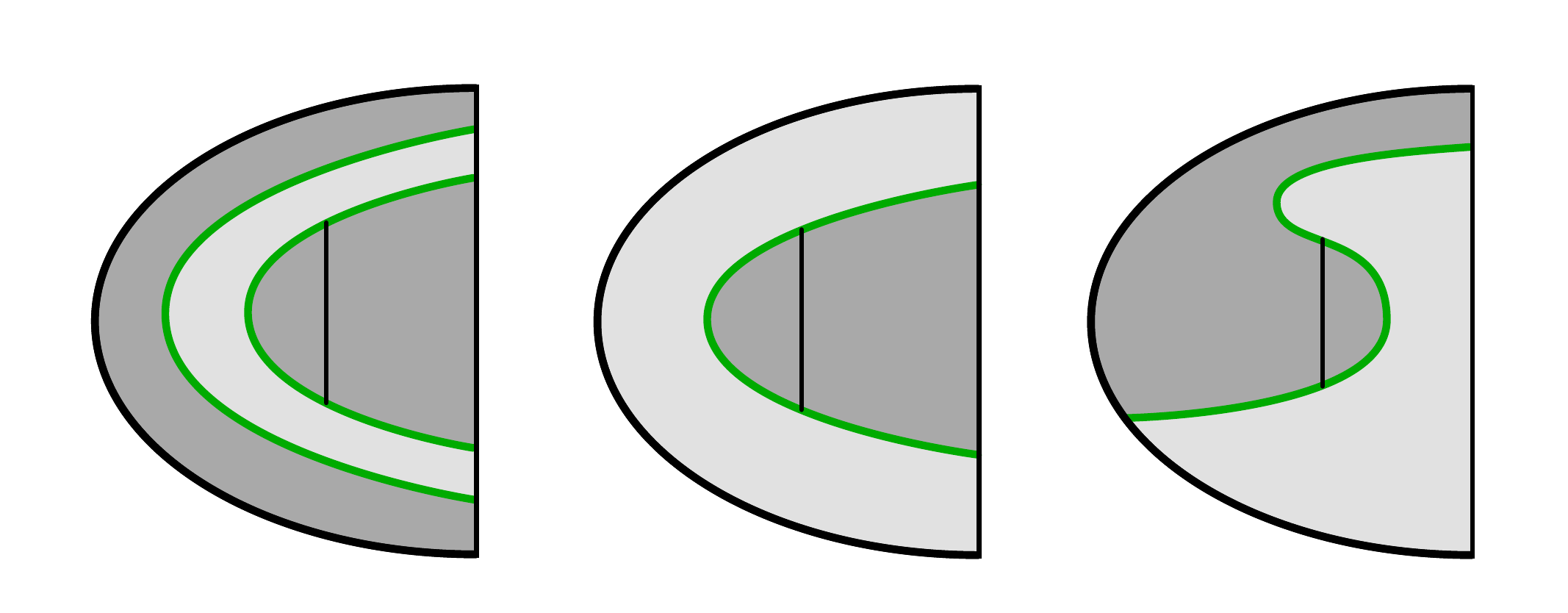}
        
        \put(31,31){\small \textcolor{darkgreen}{$p$}}

        \put(7,19){\tiny $-0$}
        \put(12,19){\tiny $+0$}
        \put(24,19){\tiny $-\tfrac{1}{2}$}

        \put(40,19){\tiny $+0$}
        \put(54,19){\tiny $+\tfrac{1}{2}$}

        \put(75,19){\tiny $-\tfrac{1}{2}$}
        \put(87,9){\tiny $+\tfrac{1}{2}$}

	\end{overpic}
	\caption{Subdisks of Type $\text{\RNum{8}}_p$ (left and center), split into cases based on the parity of $p$, and Type $\text{\RNum{9}}$ (right).}
	\label{fig:type89}
\end{figure}

\vspace{2mm}
\textsc{Case \RNum{9}:} \textit{Boundary points of $\gamma_i^-$ are contained in a single seam.}
\vspace{2mm}

In this case, there are no additional arcs, and Type \RNum{9} is shown at right in \cref{fig:type89}.  We compute $\rot(D_i) = 0$, $I(D_i) = 0$, and $\tb(D_i) = -\tfrac{1}{2}$, so that $I(D_i) + \tb(D_i) -\rot(D_i) = -\tfrac{1}{2}$, and $D_i$ is admissible.

\vspace{2mm}

By exhaustion, we conclude that the inequality \eqref{eq:index} holds, and it follows from the argument at the beginning of the proof that every subdisk $D_i$ must be admissible.  Moreover, as a result of this case-by-case analysis, the set of admissible types is
\[ \{\text{\RNum{1}}'_{m,n,0}, \text{\RNum{3}}_{m,n,0}, \text{\RNum{4}}_{m,1}, \text{\RNum{7}}_1, \text{\RNum{8}}_{1}, \text{\RNum{9}} \},\]
where each undetermined parameter must be odd.  In Types $\text{\RNum{1}}'_{m,n,0}$ and $\text{\RNum{3}}_{m,n,0}$, $\Gamma$ meets the arc $\gamma_i^-$ in a total of $m+n$ points.  In Type $\text{\RNum{4}}_{m,1}$, $\Gamma$ meets $\gamma_i^-$ in $m+1$ points.  In Types $\text{\RNum{7}}_1$, $\text{\RNum{8}}_{1}$, and $\text{\RNum{9}}$, $\Gamma$ meets $\gamma_i^-$ twice.  It follows that in all admissible types, every arc $\gamma_i^-$ meets $\Gamma$ in an even number of points.  By the interleaving condition of \cref{lemma:normalized_disk}, every arc $\gamma_i^+$ meets $\Gamma$ in an odd number of points.  But in Type $\text{\RNum{8}}_{1}$, $\Gamma$ meets $\gamma_i^+$ twice, and in Type $\text{\RNum{3}}_{m,n,0}$, $\Gamma$ meets $\gamma_i^+$ in $2n$ points.  We conclude that neither of these types occur.

Finally, we consider $I(D_i)$ for each of the remaining four types.  For Type $\text{\RNum{1}}'_{m,n,0}$, we have $I(D_i) = m+n-2 \geq 0$, for Type $\text{\RNum{4}}_{m,1}$, we have $I(D_i) = m-1 \geq 0$, and for Types $\text{\RNum{7}}_1$ and \text{\RNum{9}}, $I(D_i) = 0$.  By \cref{lem:I}, $\sum I(D_i) = 0$, from which it follows that $I(D_i) = 0$ for all subdisks $D_i$, which further implies that the only possible types are $\text{\RNum{1}}'_{1,1,0}$, $\text{\RNum{4}}_{1,1}$, $\text{\RNum{7}}_1$, and \text{\RNum{9}}. By inspection, each possible subdisk admits an unknotted, unlinked derivative curve $K_i$ bounding a subdisk $E_i$ of $D_i$ satisfying $\tb(K_i) - \rot(E_i) =  -1$. The transverse pushoff of a Legendrian realization of $K_i$ will then have $\mathrm{sl} = -1$. It follows that $D$ supports a max-sl transverse unlink ribbon derivative.
\end{proof}

\begin{proof}[Proof of \cref{cor:slice-benn-sharp}]
Let $\mathcal{T}$ be a transverse knot with $\mathrm{sl}(\mathcal{T}) = -1$. Assume that $\mathcal{T}$ bounds a star-shaped ribbon disk. Let $\Lambda$ be a Legendrian pushoff of $\mathcal{T}$. After possibly negatively stabilizing and destabilizing $\Lambda$, there is a normalized ribbon-convex star-shaped disk satisfying the hypotheses of \cref{prop:star}. Since negative stabilization and destabilization preserves the transverse isotopy class of the pushoff, the transverse pushoff is $\mathcal{T}$. By the proposition, $\mathcal{T}$ admits a Seifert surface with a transverse unlink derivative, and \cref{thm:BZ1review} then implies that $\mathcal{T}$ is symplectically slice, hence quasipositive.
\end{proof}

\begin{proof}[Proof of \cref{cor:CNSrib2}]
Let $\Lambda$ be a Legendrian knot with $\tb(\Lambda) = -1$ and $\mathrm{rot}(\Lambda) = 0$. Assume that $\Lambda$ bounds a star-shaped ribbon disk. After possibly negatively stabilizing and destabilizing $\Lambda$, there is a normalized ribbon-convex star-shaped disk satisfying the principles assumed in \cref{prop:star}, i.e.\ the bulleted principles at the beginning of the present section. By the proposition, the transverse pushoff of the Legendrian is quasipositive. Since negative stabilization and destabilization preserves the transverse pushoff, we conclude that the transverse pushoff of $\Lambda$ is quasipositive. By \cite[Corollary 2.10]{mark2024fillable}, $\Lambda$ is Lagrangian slice.
\end{proof}

\section{Data supporting the CNS conjecture}\label{sec:data}

Finally, we explore explicit data in support of \cref{conj:CNS} in order to demonstrate the practical efficacy of \cref{{cor:CNSrib2}}.  In~\cite{cornwell2016obstructions}, the authors completed a census of Lagrangian slice knots through $14$ crossings, verifying \cref{conj:CNS} up to crossing number 14, and they narrowed down the collection of all $15$-crossing knots to $48$ potentially Lagangrian slice knots.  In this section, we explore which of these 48 knots fall under the purview of \cref{{cor:CNSrib2}}.

Via direct experimentation using SnapPy~\cite{snappy}, we constructed star-shaped ribbon disks for 11 of these knots:

\begin{equation*}
    \begin{gathered}
        15n_{38594}, 15n_{41697}, \overline{15n_{43982}}, 15n_{46734}, 15n_{77461},\\ 15n_{96161}, 15n_{112479}, 15n_{127852}, \overline{15n_{130682}}, 15n_{132539}, 15n_{135516}.
    \end{gathered}
\end{equation*}

\begin{figure}[ht]
	\centering
	\begin{overpic}[scale=.4]{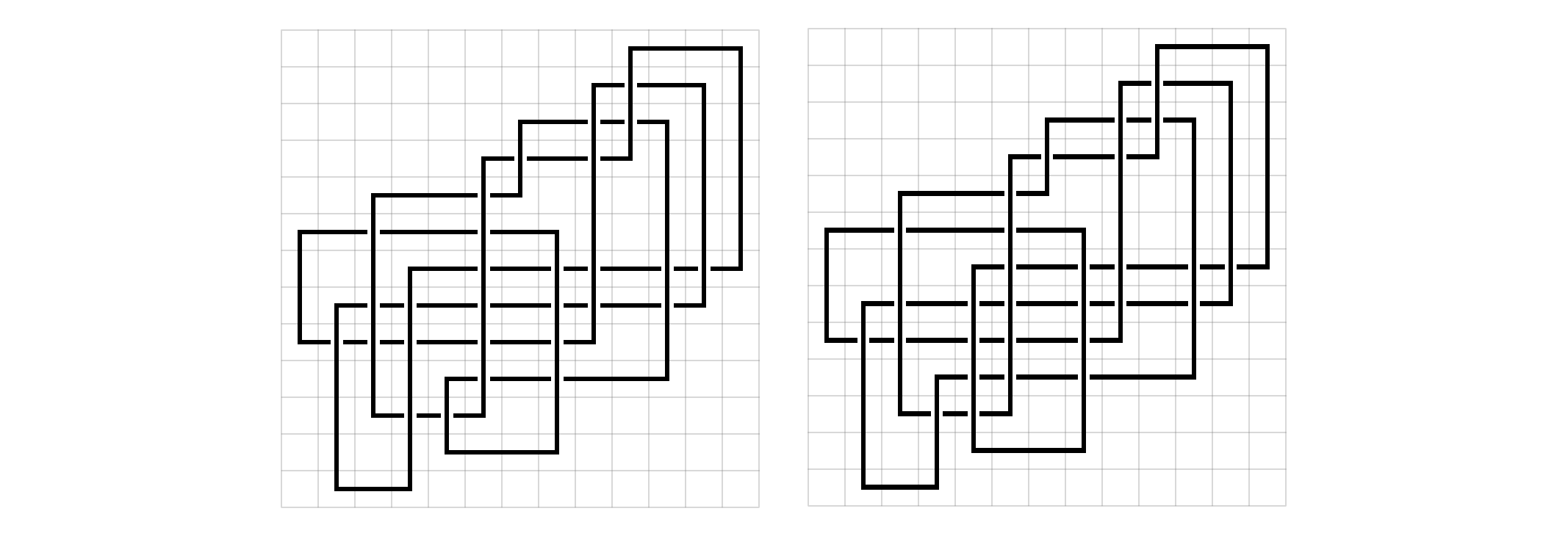}
        \put(30,-0.75){\small $15n_{38594}$}
        \put(63.5,-0.75){\small $15n_{46734}$}
        
	\end{overpic}
	\caption{Minimal grid diagrams for two $15$-crossing knots that can be used to show that they satisfy $\overline{\tb}(K) < -1$.}
	\label{fig:notLS}
\end{figure}

In order to decide whether these 11 knots satisfy the other hypothesis of \cref{{cor:CNSrib2}}, we determine $\overline{\tb(K)}$, either by directly exhibiting a grid diagram with $\tb(\Lambda) = -1$, or by using the grid diagrams appearing in~\cite{arc11} and~\cite{arc12} in conjunction with a result of Dynnikov and Prasolov~\cite{dynnikov-prasolov}, which asserts that $\overline{\tb}(K)$ is realized by any minimal grid diagram $D$ of $K$.  Two of the knots, $15n_{38594}$ and $15n_{46734}$, do not appear in any census of knots with arc index at most 12, implying that they have arc index at least 13.  For these knots, we exhibit minimal grid diagrams with arc index 13 (\cref{fig:notLS}), which verify that $\overline{\tb}(K) < -1$, confirming that they do not satisfy the hypotheses of \cref{{cor:CNSrib2}}, and in addition, they are not Lagrangian slice.

\begin{figure}[ht]
	\centering
	\begin{overpic}[scale=.4]{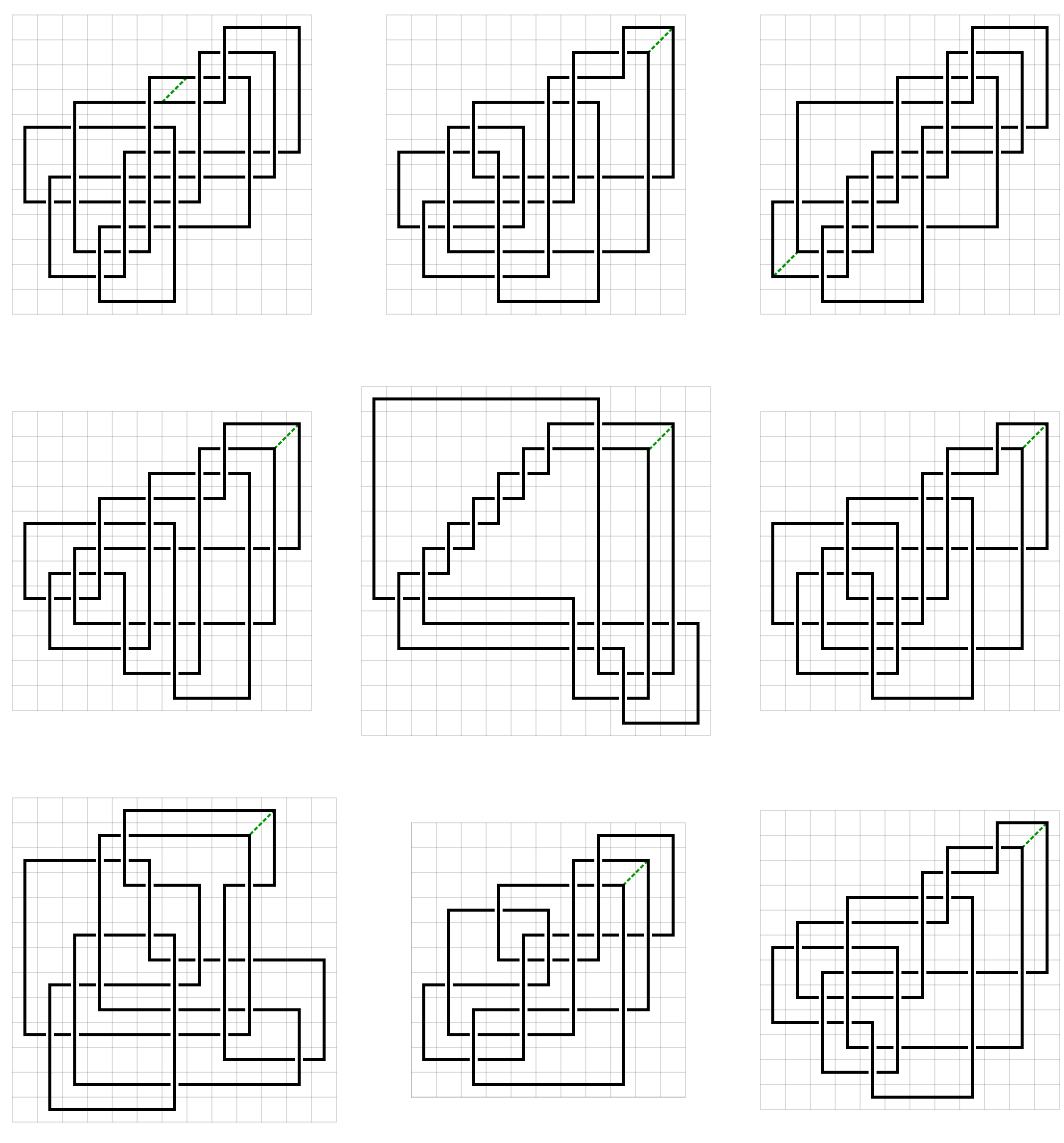}
        \put(11,-1.25){\small $\overline{15n_{130682}}$}
        \put(44.5,1.25){\small $15n_{132539}$}
        \put(75.5,0.5){\small $15n_{135516}$}

        \put(11,35.25){\small $15n_{96161}$}
        \put(44.5,33.25){\small $15n_{112479}$}
        \put(75.5,35.25){\small $15n_{127852}$}

        \put(11,70.25){\small $15n_{41697}$}
         \put(44.5,69.5){\small $\overline{15n_{43982}}$}
          \put(75.5,70.25){\small $15n_{77461}$}
	\end{overpic}
	\caption{Nine $15$-crossing Lagrangian slice knots.}
	\label{fig:grid-new}
\end{figure}

For the remaining nine knots, we show that $\overline{\tb}(K) = -1$, implying that each is Lagrangian slice by \cref{{cor:CNSrib2}}.  We further demonstrate that each knot $K$ is decomposably slice by explicitly exhibiting Legendrian pinch moves converting $K$ to a $2$-component max-tb unlink, as in~\cite{cornwell2016obstructions}.  See~\cite{breen2024regularlysliceimpliesoncestably} for further discussion about the distinction between Lagrangian slicness and decomposable sliceness.  These pinch moves are indicated with dashed lines in \cref{fig:grid-new}.  The associated grid diagrams are imported from~\cite{arc11} for $15n_{132539}$ and from~\cite{arc12} for $15n_{41697}$, $\overline{15n_{43982}}$, $15n_{77461}$, $15n_{96161}$, $15n_{127852}$, and $15n_{135516}$.  For $15n_{112479}$ and $\overline{15n_{130682}}$, we constructed explicit grid diagrams via experimentation in SnapPy.

\bibliography{references}
\bibliographystyle{amsalpha}
\end{document}